\documentclass[article]{IEEEtran}

\ifCLASSINFOpdf
\else
\fi

\usepackage{amssymb}
\usepackage{graphicx}
\usepackage{subfigure}
\usepackage{amsmath}
\usepackage{color}
\usepackage{setspace}
\usepackage{url}
\usepackage{threeparttable}
\usepackage{tabularx}
\usepackage{booktabs} % For formal tables

\begin{document}
\title{R2-based Hypervolume Contribution Approximation}

\author{Ke~Shang,
        Hisao~Ishibuchi,~\IEEEmembership{Fellow,~IEEE,}
        and~Xizi~Ni% <-this % stops a space
\thanks{This work was supported by National Natural Science Foundation of China (Grant No. 61876075), the Program for Guangdong Introducing Innovative and Enterpreneurial Teams (Grant No. 2017ZT07X386), Shenzhen Peacock Plan (Grant No. KQTD2016112514355531), the Science and Technology Innovation Committee Foundation of Shenzhen (Grant No. ZDSYS201703031748284) and the Program for University Key Laboratory of Guangdong Province (Grant No. 2017KSYS008).}
\thanks{K. Shang, H. Ishibuchi, and X. Ni are with Shenzhen Key Laboratory of Computational Intelligence, University Key Laboratory of Evolving Intelligent Systems of Guangdong Province, Department of Computer Science and Engineering, Southern University of Science and Technology, Shenzhen 518055, China (e-mail: kshang@foxmail.com; hisao@sustc.edu.cn).}% <-this % stops a space
\thanks{Corresponding Author: H. Ishibuchi.}
\thanks{Manuscript received XXX; revised XXX.}
}

\markboth{IEEE TRANSACTIONS ON EVOLUTIONARY COMPUTATION, VOL. X, NO. X, MONTH YEAR}%
{Shell \MakeLowercase{\textit{et al.}}: Bare Demo of IEEEtran.cls for IEEE Journals}
\maketitle

\begin{abstract}
In this letter, a new hypervolume contribution approximation method is proposed which is formulated as an R2 indicator. The basic idea of the proposed method is to use different line segments only in the hypervolume contribution region for the hypervolume contribution approximation. Comparing with a traditional method which is based on the R2 indicator to approximate the hypervolume, the new method can directly approximate the hypervolume contribution and will utilize all the direction vectors only in the hypervolume contribution region. The new method, the traditional method and the Monte Carlo sampling method \textcolor{black}{together with two exact methods} are compared through comprehensive experiments. \textcolor{black}{Our results show the advantages of the new method over the other methods. Comparing with the other two approximation methods, the new method achieves the best performance for comparing hypervolume contributions of different solutions and identifying the solution with the smallest hypervolume contribution. Comparing with the exact methods, the new method is computationally efficient in high-dimensional spaces where the exact methods are impractical to use.}
\end{abstract}


\IEEEpeerreviewmaketitle

\section{Introduction}
\label{introduction}
Hypervolume \cite{zitzler1998multiobjective} is a widely used performance indicator in the Evolutionary Multi-objective Optimization (EMO) community due to its well-known unique property (i.e., Pareto compliance) among all existing indicators. The bottleneck of the hypervolume applicability in Evolutionary Multi-objective Optimization Algorithms (EMOA) is the increasing computational burden as the dimensionality of the objective space increases. Whereas several fast hypervolume calculation methods \cite{bradstreet2010fast,While2012A,Russo2012Quick,Guerreiro2018Computing} have been proposed, it has been proved that the exact hypervolume calculation is \#P-hard in the number of dimensions \cite{Bringmann2010Approximating}. Therefore, efforts in the hypervolume approximation have been done to increase the applicability of the hypervolume to  high-dimensional spaces, including the Monte Carlo sampling method \cite{bader2010faster,bader2011hype,Bringmann2012Approximating} and the achievement scalarizing function method \cite{ishibuchi2009hypervolume, ishibuchi2010indicator}. 

\textcolor{black}{In the Monte Carlo sampling method, the sampling space of a solution set is determined first, then a large number of samples are drawn evenly in this sampling space to estimate the hypervolume of the corresponding solution set. A sample is called a hit if it is dominated by the solution set, otherwise it is called a miss. Then the hypervolume is approximated by the ratio of the number of hits to the total number of the samples multiplies the volume of the sampling space.}

In the achievement scalarizing function method, the hypervolume is approximated by a number of achievement scalarizing functions with uniformly distributed weight vectors. Each achievement scalarizing function with a different weight vector is used to measure the distance from the reference point to the attainment surface of the solution set. Then the average distance from the reference point to the attainment surface over a large number of weight vectors is calculated as the hypervolume approximation. 

The achievement scalarizing function method can be formulated as an R2 indicator. In \cite{ma2017tchebycheff}, an R2 indicator is proposed based on a new Tchebycheff function. The proposed R2 indicator shows a clear geometric property for the approximation of the hypervolume. In \cite{shangke}, a new R2 indicator is proposed based on the Divergence theorem and Riemann sum approximation for better hypervolume approximation. The new R2 indicator significantly improved the approximation quality for the hypervolume compared with the achievement scalarizing function method.

In the hypervolume-based EMOAs (e.g., SMS-EMOA \cite{Emmerich2005An,beume2007sms}, FV-MOEA \cite{jiang2015simple}), the hypervolume contribution is used to evaluate the fitness value of each individual. In SMS-EMOA, with a steady-state $(\mu+1)$ ES-type generation update mechanism, the worst individual with the smallest hypervolume contribution is eliminated from the population. \textcolor{black}{In FV-MOEA, with a $(\mu+\lambda)$ ES-type generation update mechanism, the individual with the smallest hypervolume contribution is removed one by one from the population.} Thus, the hypervolume contribution plays an important role in the hypervolume-based EMOAs. 

\textcolor{black}{There are some works devoted to calculate the hypervolume contribution exactly such as IHSO \cite{bradstreet2008fast}, IWFG \cite{while2012applying} and exQHV \cite{russo2016extending}. However, the exact hypervolume contribution calculation is still \#P-hard \cite{Bringmann2012Approximating}, so there is a need to develop an approximation method for the hypervolume contribution calculation especially for high-dimensional spaces.}

This letter investigates the hypervolume contribution approximation methods based on the R2 indicator. The main contribution of this letter is that a new method which is formulated as an R2 indicator is proposed for the hypervolume contribution approximation. \textcolor{black}{We conduct comprehensive experiments to test the performance of the new method and show its advantages over the other methods. Some interesting observations and insights are obtained from our results.}

%The rest of the letter is organized as follows. Section \ref{preliminaries} gives the preliminaries of the letter. R2-based hypervolume contribution approximation methods are investigated in Section \ref{newr2forhv}. Numerical studies are given in Section \ref{numstudy}. We conclude the letter in Section \ref{conclude}.

\section{Preliminaries}
\label{preliminaries}
\subsection{Hypervolume and hypervolume contribution}
\label{section:hv}
Given a reference point $\mathbf{r}$ and an approximation solution set $A$, the hypervolume of the set $A$ is defined as:
\begin{equation}
HV(A,\mathbf{r}) = \mathcal{L}\left(\bigcup_{\mathbf{a}\in A}\left\{\mathbf{b}|\mathbf{a}\succ \mathbf{b}\succ \mathbf{r}\right\}\right),
\end{equation}
where $\mathcal{L}(.)$ is the Lebesgue measure of a set, and $\mathbf{a}\succ \mathbf{b}$ means $\mathbf{a}$ dominates $\mathbf{b}$.

For a solution $\mathbf{s}\in A$, its hypervolume contribution is defined as:
\begin{equation}
C_{HV}(\mathbf{s},A,\mathbf{r}) = HV(A,\mathbf{r}) - HV(A\setminus \{\mathbf{s}\},\mathbf{r}).
\end{equation}

\begin{figure}[!htbp]
\centering                                           %居中
\subfigure[Hypervolume]{      
\hspace{-0.5cm}              %第一张子图
\begin{minipage}{0.45\columnwidth}\centering                                                          %子图居中
\includegraphics[scale=0.35]{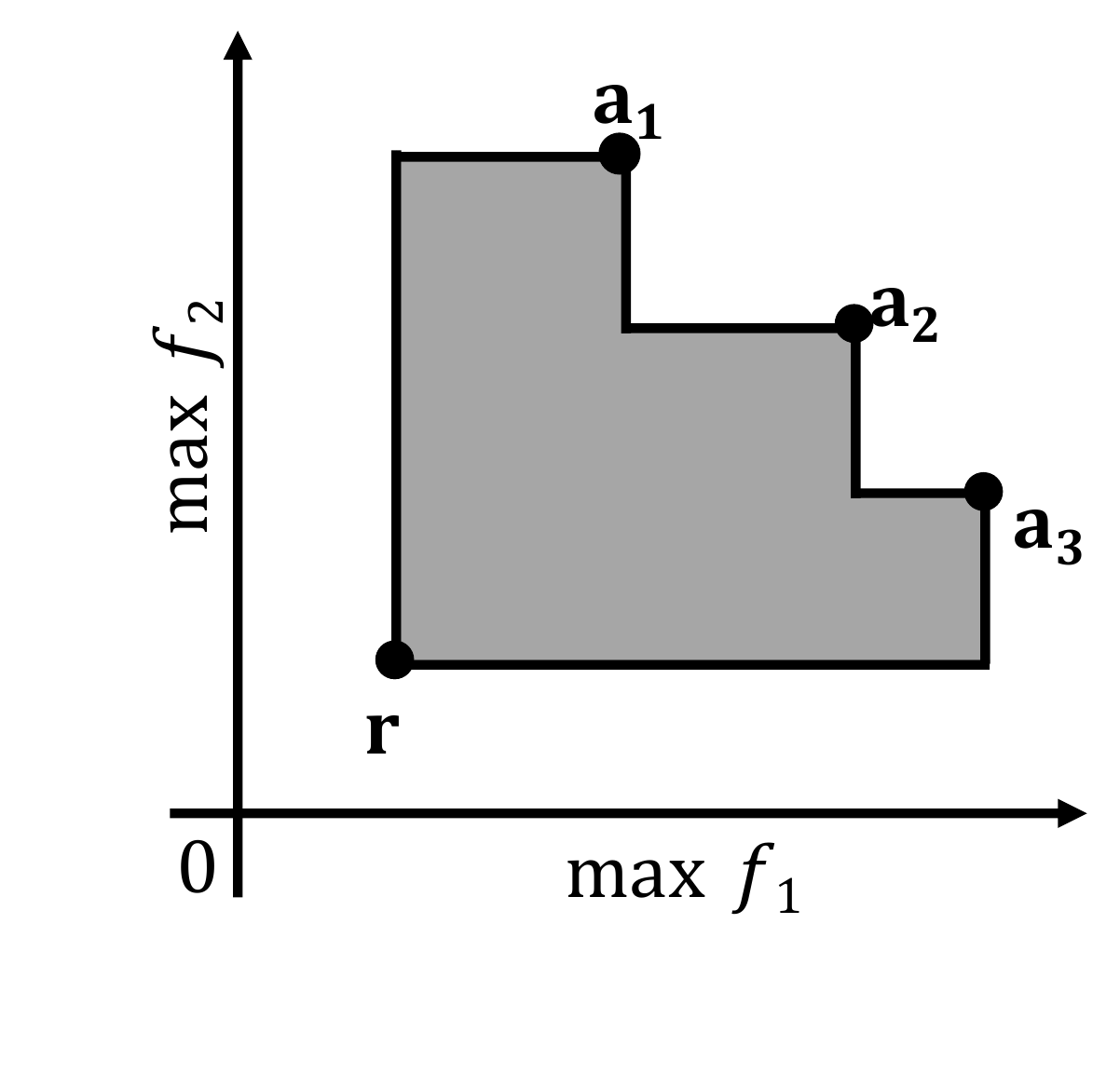}               %以pic.jpg的0.5倍大小输出
\end{minipage}}
\subfigure[Hypervolume contribution]{                    %第二张子图
\begin{minipage}{0.45\columnwidth}\centering                                                          %子图居中
\includegraphics[scale=0.35]{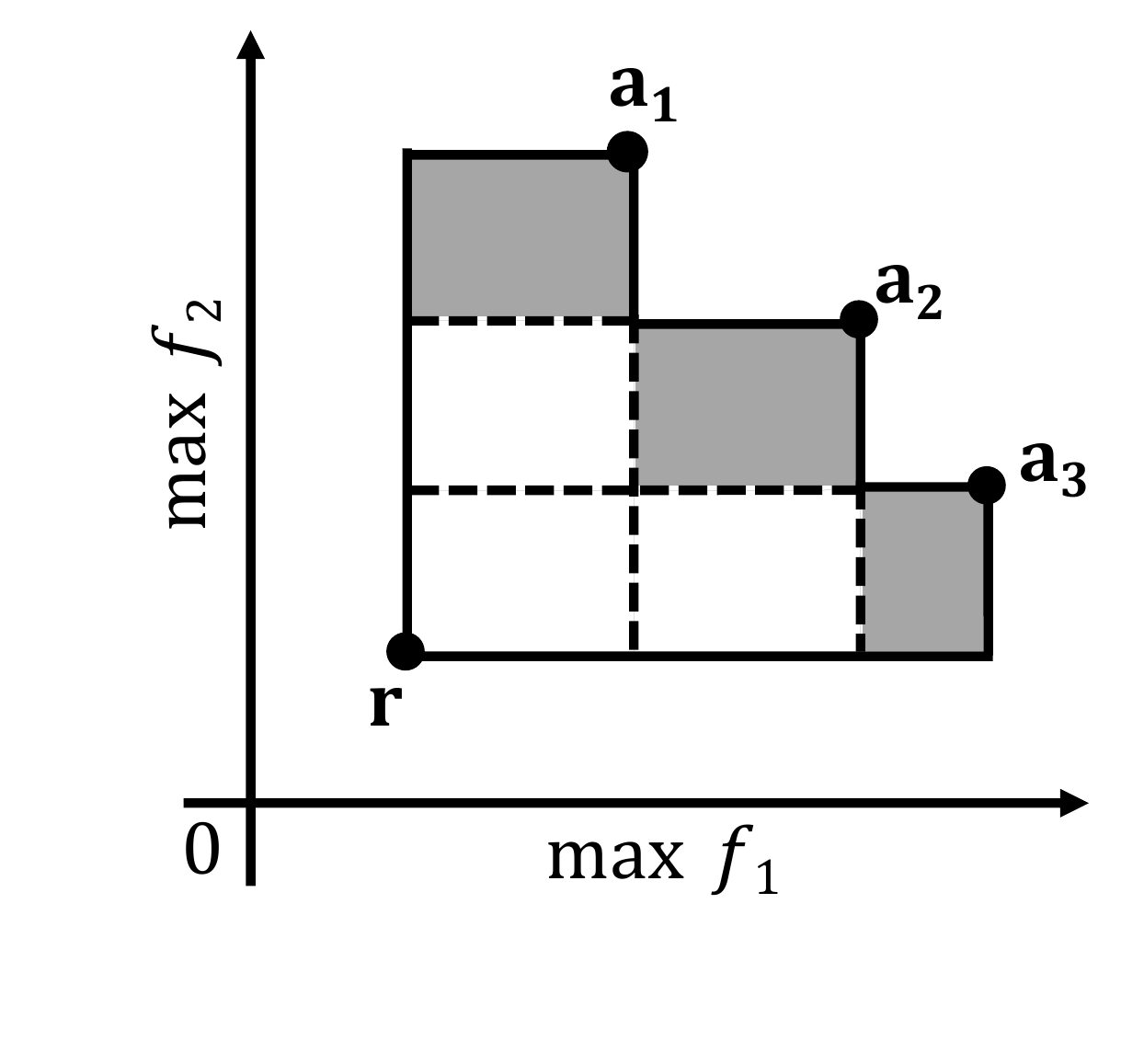}                %以pic.jpg的0.5倍大小输出
\end{minipage}}
\caption{An illustration of the hypervolume and the hypervolume contribution} %                         %大图名称
\label{hv}                                                        %图片引用标记
\end{figure}
Fig. \ref{hv} gives an illustration of the hypervolume of a solution set and the hypervolume contribution of each solution\footnote{\textcolor{black}{Throughout of the letter, maximization of all objectives is assumed.}}. The hypervolume is the shaded area of the enclosed polygon in Fig. \ref{hv} (a). The hypervolume contribution of each solution is the shaded area which is uniquely dominated by the corresponding solution in Fig. \ref{hv} (b).

\subsection{R2 indicator for hypervolume approximation}

Given a solution set $A$, a set of direction vectors $\Lambda$, and a utopian point $\mathbf{r}^*$, the R2 indicator based on the 2-Tch function \cite{ma2017tchebycheff} is defined for an $m$-objective problem as:
\begin{equation}
\label{r2tch}
R^{\textrm{2tch}}_2(A,\Lambda,\mathbf{r}^*)=\frac{1}{|\Lambda|}\sum_{\lambda\in \Lambda}\min_{\mathbf{a}\in A}\left\{g^{\textrm{2tch}}(\mathbf{a}|\lambda,\mathbf{r}^*)\right\},\end{equation}
where the 2-Tch function is defined as follows: 
\begin{equation}
\label{g2tch}
g^{\textrm{2tch}}(\mathbf{a}|\lambda,\mathbf{r}^*)=\max_{j\in\{1,...,m\}}\left\{\frac{|r^*_j-a_j|}{\lambda_j}\right\}.
\end{equation}
In Eq.~\eqref{g2tch}, $\lambda = (\lambda_1,\lambda_2,...,\lambda_m)$ is a given direction vector with $\left \| \lambda \right \|_2 = 1$ and $\lambda_i\geq 0$, $i = 1,...,m$.

\begin{figure}[!htbp]
\centering                                                          %居中
\subfigure[$R^{\textrm{2tch}}_2$]{        
\hspace{-0.5cm}            %第一张子图
\begin{minipage}{0.45\columnwidth}\centering                                                          %子图居中
\includegraphics[scale=0.35]{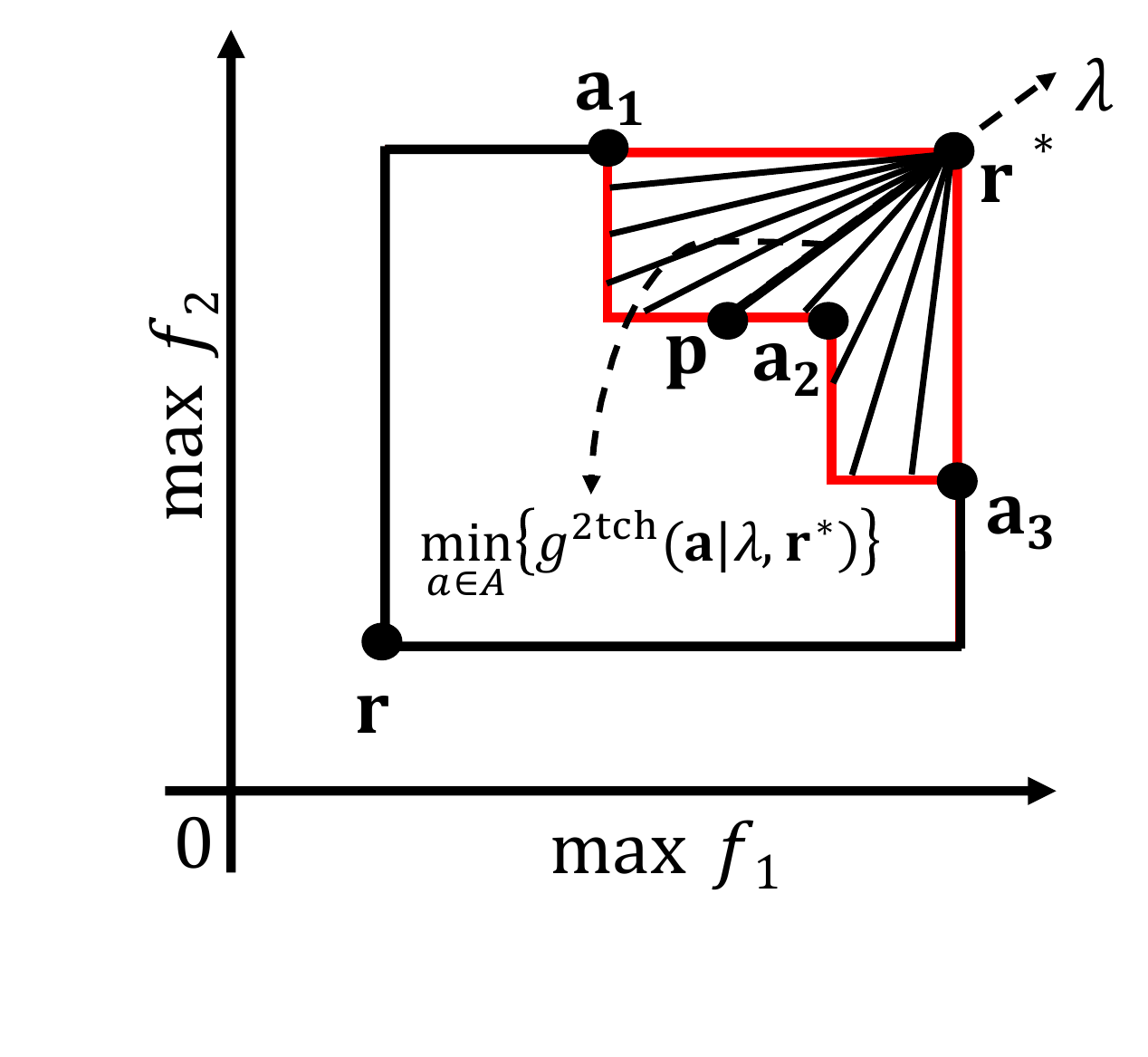}               %以pic.jpg的0.5倍大小输出
\end{minipage}}
\subfigure[$R^{\textrm{mtch}}_2$]{                    %第二张子图
\begin{minipage}{0.45\columnwidth}\centering                                                          %子图居中
\includegraphics[scale=0.35]{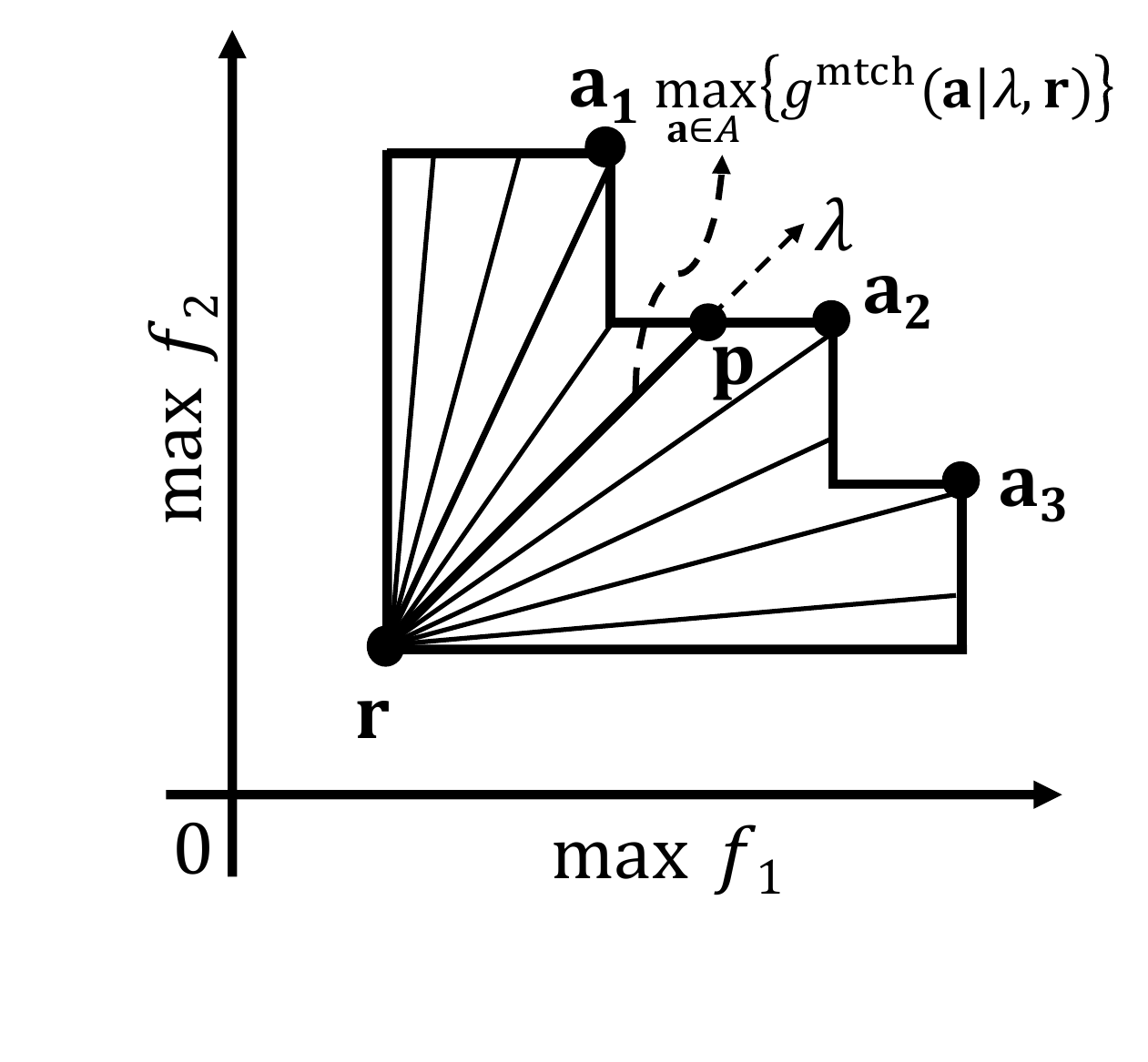}                %以pic.jpg的0.5倍大小输出
\end{minipage}}
\caption{An illustration of the geometric property of $R^{\textrm{2tch}}_2$ and $R^{\textrm{mtch}}_2$.}\label{r2}                                                        
\end{figure}

Fig. \ref{r2} (a) shows the geometric property of $R^{\textrm{2tch}}_2$. In this figure, $\mathbf{r}^*$ is the utopian point and $\mathbf{r}$ is the reference point for the hypervolume calculation. Suppose a line follows the direction $\lambda$, passes through $\mathbf{r}^*$ and intersects with the attainment surface of the solution set $A$ at $\mathbf{p}$, then the length of the line segment with the end points $\mathbf{r}^*$ and $\mathbf{p}$ is determined by $\min_{\mathbf{a}\in A}\left\{g^{\textrm{2tch}}(\mathbf{a}|\lambda,\mathbf{r}^*)\right\}$. $R^{\textrm{2tch}}_2$ in Eq.~\eqref{r2tch} is the average length of the line segments with different directions as shown in Fig. \ref{r2} (a).

The idea of using different line segments starting from a reference point to the attainment surface of the solution sets for the hypervolume approximation was firstly proposed in \cite{ishibuchi2009hypervolume} as shown in Fig.~\ref{r2} (b). Intuitively, the average length of the line segments in Fig. \ref{r2} (b) is closely related to the hypervolume.

The $R^{\textrm{2tch}}_2$ indicator cannot be directly used for the hypervolume approximation as shown in Fig. \ref{r2} (a). Its modified version which is able to directly approximate the hypervolume is defined as follows:
\begin{equation}
\label{rmtch}
R^{\textrm{mtch}}_2(A,\Lambda,\mathbf{r})=\frac{1}{|\Lambda|}\sum_{\lambda\in \Lambda}\max_{\mathbf{a}\in A}\left\{g^{\textrm{mtch}}(\mathbf{a}|\lambda,\mathbf{r})\right\},\end{equation}
where $g^{\textrm{mtch}}$ function is defined as follows:
\begin{equation}
\label{gmtch}
g^{\textrm{mtch}}(\mathbf{a}|\lambda,\mathbf{r})=\min_{j\in\{1,...,m\}}\left\{\frac{|r_j-a_j|}{\lambda_j}\right\}.
\end{equation}
In Eq.~\eqref{gmtch}, $\mathbf{r}$ is the reference point for the hypervolume calculation, $\lambda$ is a given direction vector which is the same as in $g^{\textrm{2tch}}$ in Eq. \eqref{g2tch}. \textcolor{black}{Comparing $R^{\textrm{2tch}}_2$ and $R^{\textrm{mtch}}_2$, we can see that $R^{\textrm{mtch}}_2$ is obtained from $R^{\textrm{2tch}}_2$ by changing the min (max) sign to max (min) and the utopian point $\mathbf{r}^*$ to the reference point $\mathbf{r}$.}

It is easy to show that the $R^{\textrm{mtch}}_2$ indicator has a geometric property as shown in Fig. \ref{r2} (b), which can be used directly for the hypervolume approximation.

In our previous work \cite{shangke}, a new R2 indicator is proposed for better hypervolume approximation. The new R2 indicator is derived based on the Divergence theorem and Riemann sum approximation. It is defined as follows:
\begin{equation}
\label{newr2}
R^{\textrm{HV}}_2(A,\Lambda,\mathbf{r},m)=\frac{1}{|\Lambda|}\sum_{\lambda\in \Lambda}\max_{\mathbf{a}\in A}\left\{g^{\textrm{mtch}}(\mathbf{a}|\lambda,\mathbf{r})\right\}^m.
\end{equation}

Comparing Eq. \eqref{newr2} with Eq. \eqref{rmtch}, we can see that the only difference between $R^{\textrm{mtch}}_2$ and $R^{\textrm{HV}}_2$ is the added exponential $m$ in $R^{\textrm{HV}}_2$. This small change in  $R^{\textrm{HV}}_2$ significantly improves the approximation quality of the R2 indicator for the hypervolume. For detailed explanations of the new R2 indicator, please refer to our previous work \cite{shangke}.

\section{R2-based hypervolume contribution approximation}
\label{newr2forhv}
\subsection{Traditional method and its drawbacks}
\label{traditionmethod}
In order to use the R2 indicator (e.g., $R^{\textrm{mtch}}_2$ and $R^{\textrm{HV}}_2$) to approximate the hypervolume contribution, the simplest and straightforward method is to use the R2 contribution to approximate the hypervolume contribution.

For a solution $\mathbf{s}\in A$, its R2 contribution is defined as follows:
\begin{equation}
\label{r2c}
C_{R_2}(\mathbf{s},A,\Lambda,\mathbf{r}) = R_2(A,\Lambda,\mathbf{r}) - R_2(A\setminus \{\mathbf{s}\},\Lambda,\mathbf{r}),
\end{equation}
where $R_2$ can be $R^{\textrm{mtch}}_2$ and $R^{\textrm{HV}}_2$ for the hypervolume approximation.

\begin{figure}[!htb]
\centering                                                          %居中
\subfigure[Step 1]{      
\hspace{-0.5cm}              %第一张子图
\begin{minipage}{0.3\columnwidth}\centering                                                          %子图居中
\includegraphics[scale=0.25]{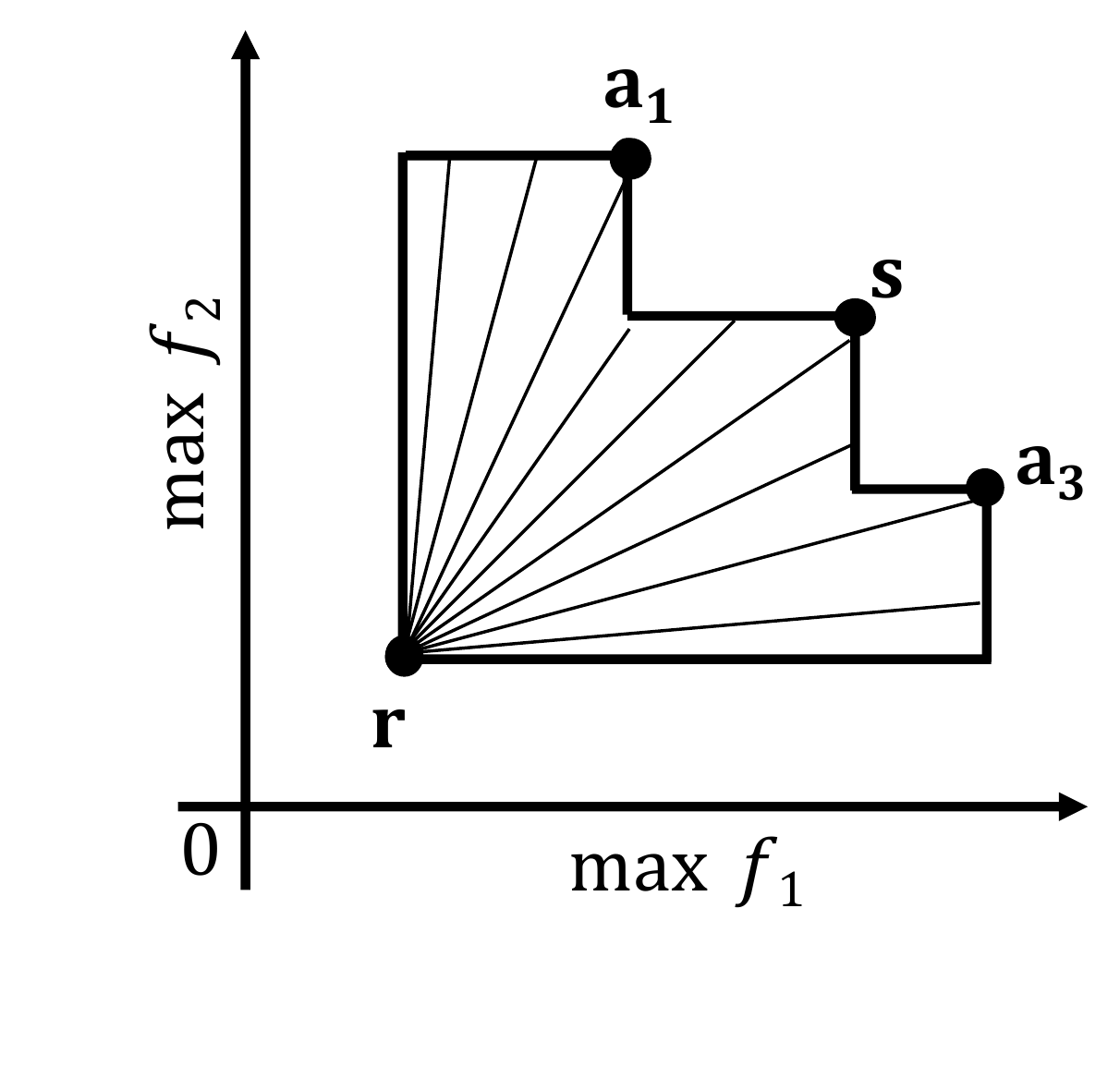}               %以pic.jpg的0.5倍大小输出
\end{minipage}}
\subfigure[Step 2]{                    %第二张子图
\begin{minipage}{0.3\columnwidth}\centering                                                          %子图居中
\includegraphics[scale=0.25]{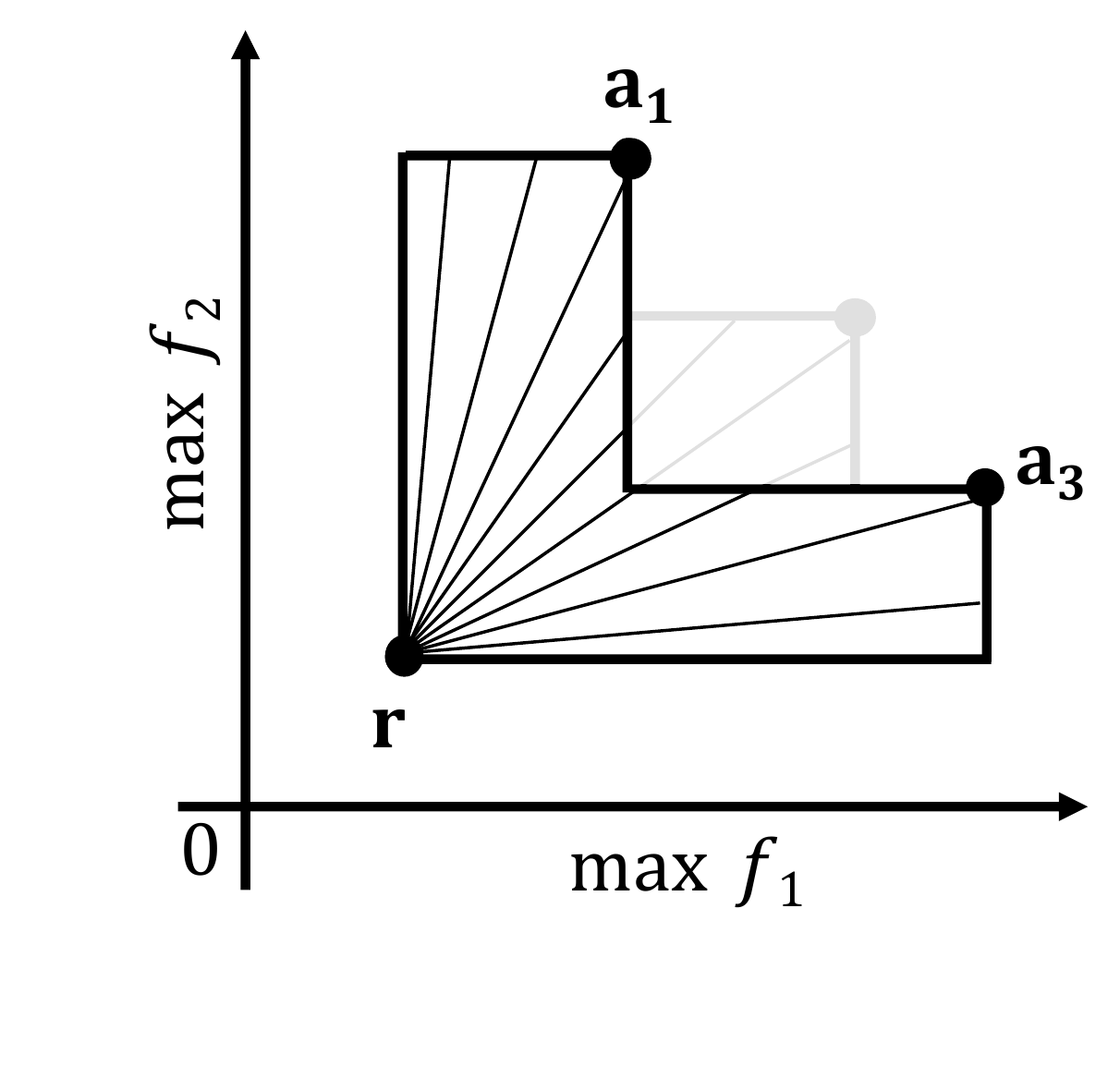}                %以pic.jpg的0.5倍大小输出
\end{minipage}}
\subfigure[Step 3]{                    %第一张子图
\begin{minipage}{0.3\columnwidth}\centering                                                          %子图居中
\includegraphics[scale=0.25]{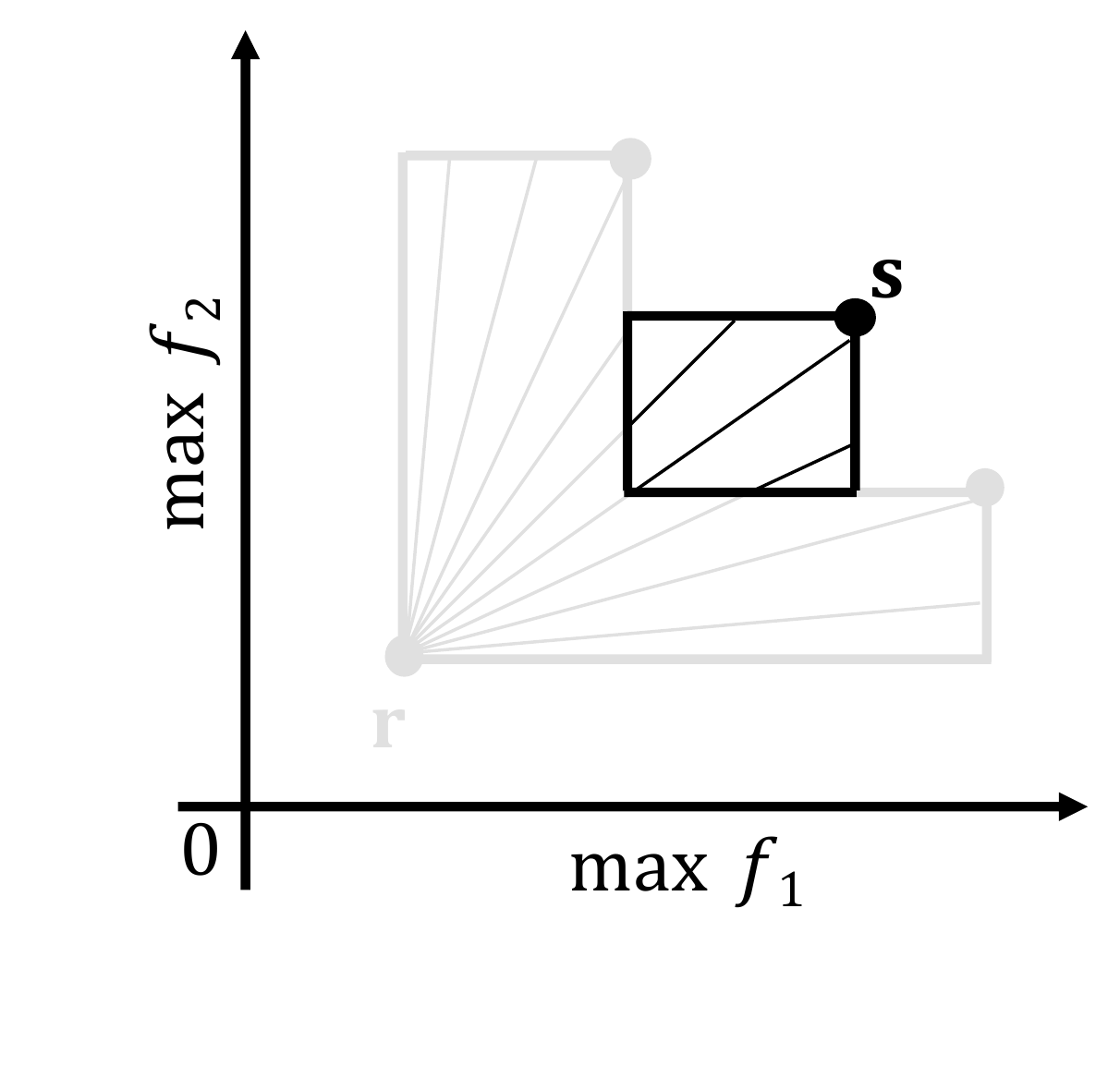}               %以pic.jpg的0.5倍大小输出
\end{minipage}}
\caption{An illustration of the traditional method for the hypervolume contribution approximation. The hypervolume contribution of solution $\mathbf{s}$ is approximated by the difference between the R2 values of the solution sets $\{\mathbf{a}_1,\mathbf{s},\mathbf{a}_3\}$ and $\{\mathbf{a}_1,\mathbf{a}_3\}$. 
}\label{r2appro}                                                        %图片引用标记
\end{figure}

According to the definition of the R2 contribution, the traditional method for the hypervolume contribution approximation (as illustrated in Fig.~\ref{r2appro}) can be described by the following three steps:
\begin{itemize}
\item Step 1: Calculate the R2 value of the solution set $A$.
\item Step 2: Calculate the R2 value of the solution set $A\setminus \{\mathbf{s}\}$.
\item Step 3: The hypervolume contribution approximation of the solution $\mathbf{s}$ is obtained by calculating the difference between the above two R2 values.
\end{itemize}

The drawbacks of the traditional method are summarized as follows:
\begin{enumerate}
\item In order to obtain the hypervolume contribution approximation of a solution, we need to calculate R2 values for two solution sets. The hypervolume contribution cannot be approximated directly.
\item Usually errors exist in the R2 approximation for the hypervolume. So the error of the hypervolume contribution approximation could be amplified by the errors of the two R2 values. For this reason, the approximation accuracy of the traditional method may be low.
\item In order to improve the approximation accuracy, a large number of vectors are needed for calculating the R2 indicator. Thus, the amount of computation in the traditional method could be large so as to achieve a high approximation accuracy.
\end{enumerate}
\subsection{A new method for hypervolume contribution approximation}
\label{newmethod1}
\begin{figure}[!htb]
\centering                                                          %居中
\subfigure[Situation 1]{        
\hspace{-0.5cm}            %第一张子图
\begin{minipage}{0.45\columnwidth}\centering                                                          %子图居中
\includegraphics[scale=0.35]{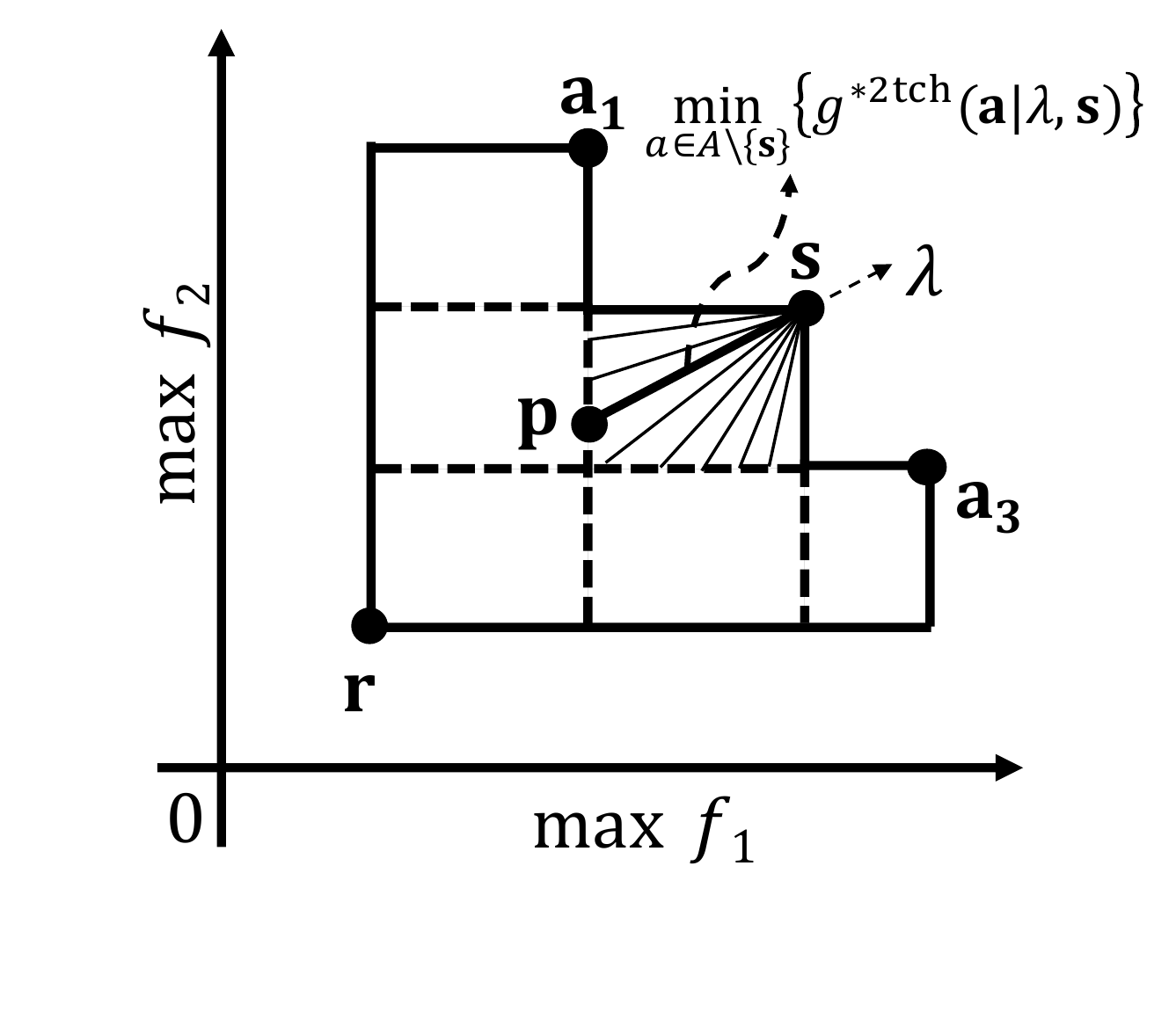}               %以pic.jpg的0.5倍大小输出
\end{minipage}}
\subfigure[Contour lines of $g^{*\textrm{2tch}}(\mathbf{a}|\lambda,\mathbf{s})$]{                    %第二张子图
\begin{minipage}{0.45\columnwidth}\centering                                                          %子图居中
\includegraphics[scale=0.35]{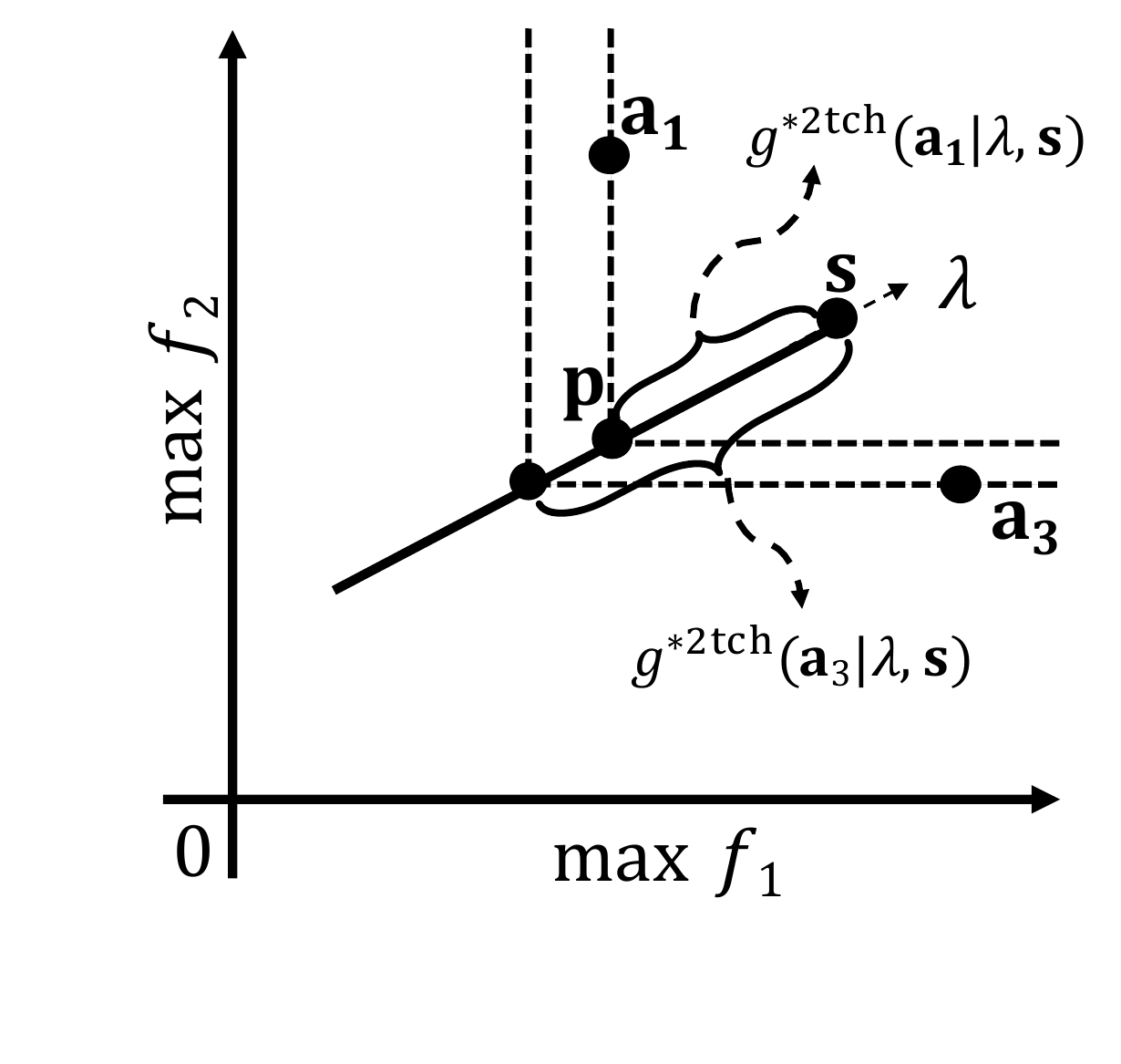}                %以pic.jpg的0.5倍大小输出
\end{minipage}}
\caption{An illustration of the proposed idea (Situation 1).}\label{situation1}                                                        %图片引用标记
\end{figure}

In this subsection, we propose a new method for the hypervolume contribution approximation. The proposed idea is illustrated in Fig.~\ref{situation1} (a). In Fig.~\ref{situation1} (a), the hypervolume contribution region of a solution $\mathbf{s}$ is occupied by the line segments with different directions which start from the solution $\mathbf{s}$ and end on the boundaries of the hypervolume contribution region of the solution $\mathbf{s}$. Then we can utilize all the line segments in the hypervolume contribution region to approximate the corresponding hypervolume contribution, which is similar to $R^{\textrm{mtch}}_2$ or $R^{\textrm{HV}}_2$ for the hypervolume approximation.

First, let us consider the situation in Fig.~\ref{situation1} (a). In this situation, the line segments in the hypervolume contribution region of the solution $\mathbf{s}$ only intersect with the attainment surface of the solution set $A\setminus \{\mathbf{s}\}$, i.e., the lengths of the line segments are only determined by the solution set $A\setminus \{\mathbf{s}\}$. This situation is similar to $R^{\textrm{2tch}}_2$ as illustrated in Fig.~\ref{r2} (a).

Given a direction vector $\lambda$, a solution $\mathbf{s}$ and a solution set $A\setminus \{\mathbf{s}\}$, the length of the line segment in Fig.~\ref{situation1} (a) can be calculated by the following formula:
\begin{equation}
\label{length1}
L\left(\mathbf{s},A,\lambda\right) =\min_{\mathbf{a}\in A\setminus \{\mathbf{s}\}}\left\{g^{*\textrm{2tch}}(\mathbf{a}|\lambda,\mathbf{s})\right\}.
\end{equation}
In Eq.~\eqref{length1}, the $g^{*\textrm{2tch}}(\mathbf{a}|\lambda,\mathbf{s})$ function is defined for maximization problems as 
\begin{equation}
\label{g2tchnew}
g^{*\textrm{2tch}}(\mathbf{a}|\lambda,\mathbf{s})=
\max_{j\in\{1,...,m\}}\left\{\frac{s_j-a_j}{\lambda_j}\right\},
\end{equation}
For minimization problems, it is defined as\footnote{\textcolor{black}{Even though we do not consider minimization case in this letter, we still give the formulation for minimization case to make it more comprehensive. If not explicitly stated, the formulations given in this letter are applicable to both maximization and minimization cases.}}
\begin{equation}
g^{*\textrm{2tch}}(\mathbf{a}|\lambda,\mathbf{s})=
\max_{j\in\{1,...,m\}}\left\{\frac{a_j-s_j}{\lambda_j}\right\}.
\end{equation}

In the maximization case which is considered in Fig.~\ref{situation1} (a), notice that $g^{*\textrm{2tch}}(\mathbf{a}|\lambda,\mathbf{s})$  in Eq.~\eqref{g2tchnew} is slightly different from $g^{\textrm{2tch}}(\mathbf{a}|\lambda,\mathbf{r}^*)$ in Eq.~\eqref{g2tch} \textcolor{black}{in that} there is no absolute value sign in $g^{*\textrm{2tch}}(\mathbf{a}|\lambda,\mathbf{s})$. The reason can be clearly explained by Fig.~\ref{situation1} (b), which shows the contour lines of $g^{*\textrm{2tch}}(\mathbf{a}|\lambda,\mathbf{s})$ in the maximization case. \textcolor{black}{In order to correctly calculate the length of the line segment $\mathbf{s}\mathbf{p}$, $\mathbf{a}_1$ and $\mathbf{p}$ should be on the same contour line of $g^{*\textrm{2tch}}(\mathbf{a}|\lambda,\mathbf{s})$. It is clear that from $\mathbf{p}$ to $\mathbf{a}_1$, only $f_2$ value increases and $f_1$ value does not change, so $g^{*\textrm{2tch}}(\mathbf{a}_1|\lambda,\mathbf{s}) = g^{*\textrm{2tch}}(\mathbf{p}|\lambda,\mathbf{s})$ holds, which means that $\mathbf{a}_1$ and $\mathbf{p}$ are on the same contour line. However, if we use $g^{\textrm{2tch}}(\mathbf{a}|\lambda,\mathbf{s})$ which is $g^{*\textrm{2tch}}(\mathbf{a}|\lambda,\mathbf{s})$ with the absolute value sign, $\mathbf{a}_1$ and $\mathbf{p}$ may not be on the same contour line because if $f_2$ value of $\mathbf{a}_1$ is large enough and then $g^{\textrm{2tch}}(\mathbf{a}_1|\lambda,\mathbf{s}) > g^{\textrm{2tch}}(\mathbf{p}|\lambda,\mathbf{s})$ might hold, thus $\mathbf{a}_1$ and $\mathbf{p}$ will not be on the same contour line anymore.} 

\begin{figure}[!htb]
\centering                                                          %居中
\subfigure[Situation 2]{      
\hspace{-0.5cm}              %第一张子图
\begin{minipage}{0.45\columnwidth}\centering                                                          %子图居中
\includegraphics[scale=0.35]{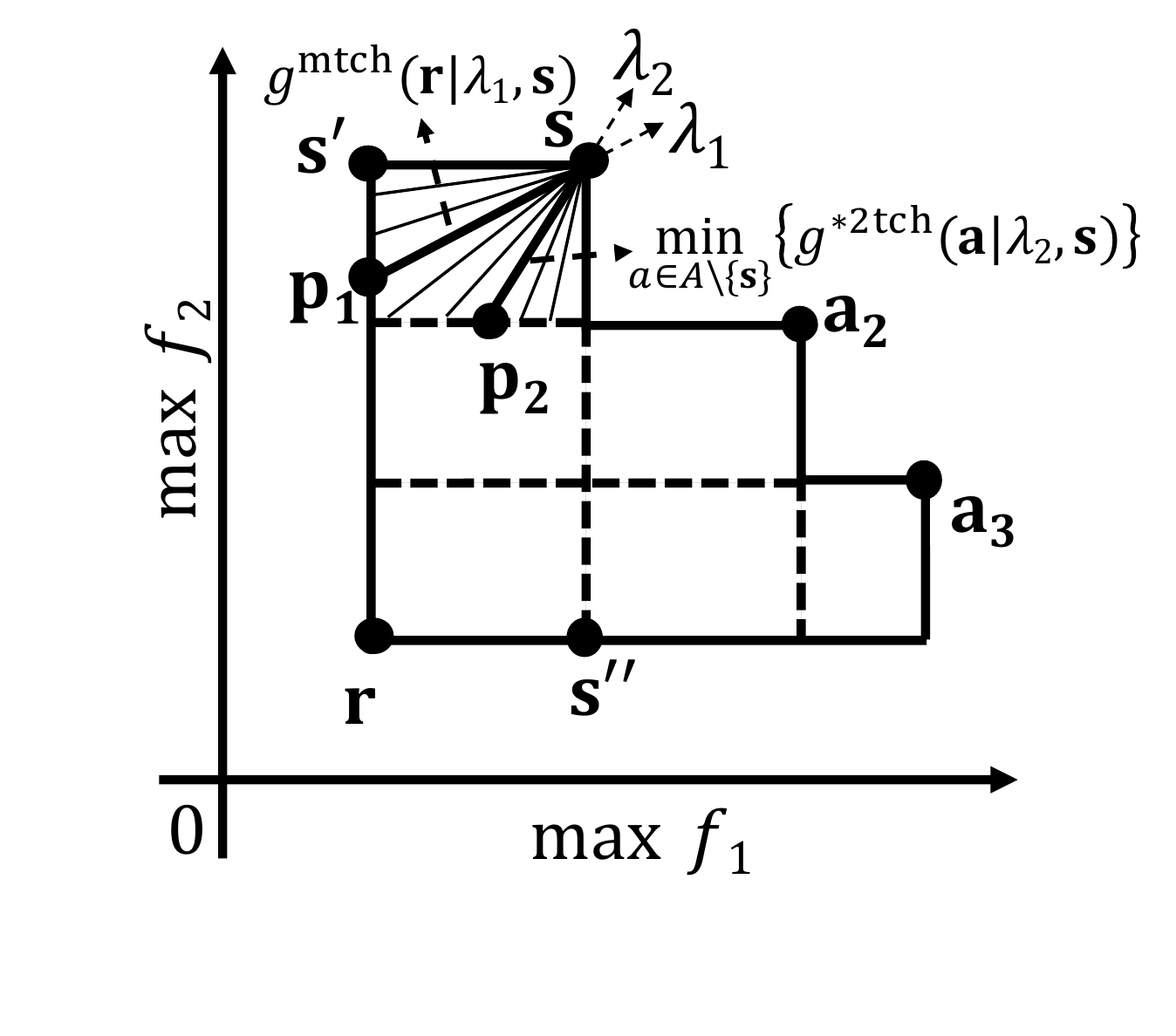}               %以pic.jpg的0.5倍大小输出
\end{minipage}}\\
\subfigure[\scriptsize Contour lines of $g^{*\textrm{2tch}}(\mathbf{a}|\lambda_1,\mathbf{s})$ and $g^{\textrm{mtch}}(\mathbf{r}|\lambda_1,\mathbf{s})$]{       
\hspace{-0.5cm}             %第二张子图
\begin{minipage}{0.45\columnwidth}\centering                                                          %子图居中
\includegraphics[scale=0.35]{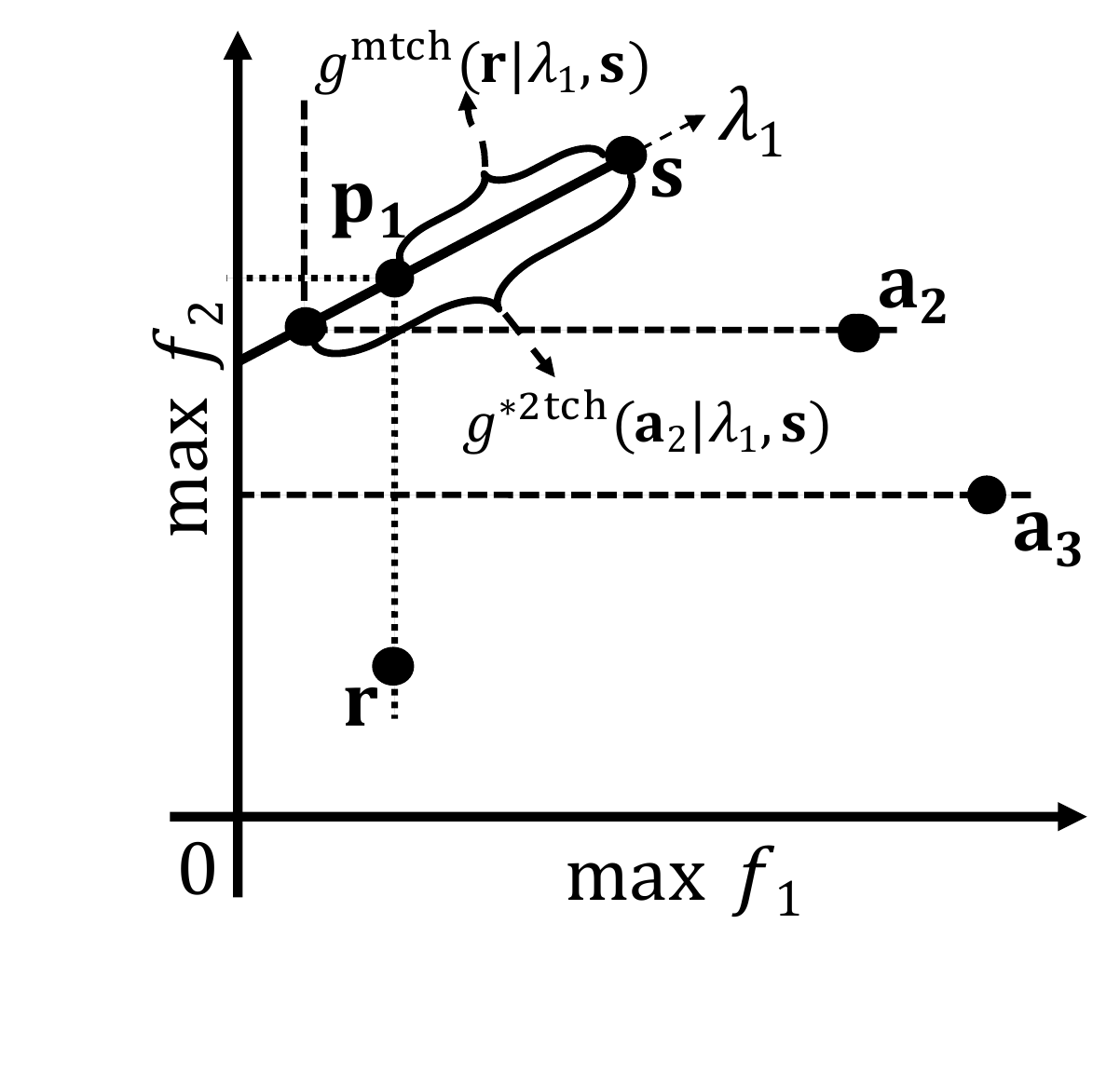}                %以pic.jpg的0.5倍大小输出
\end{minipage}}
\subfigure[\scriptsize Contour lines  of $g^{*\textrm{2tch}}(\mathbf{a}|\lambda_2,\mathbf{s})$ and $g^{\textrm{mtch}}(\mathbf{r}|\lambda_2,\mathbf{s})$]{                    %第二张子图
\begin{minipage}{0.45\columnwidth}\centering                                                          %子图居中
\includegraphics[scale=0.35]{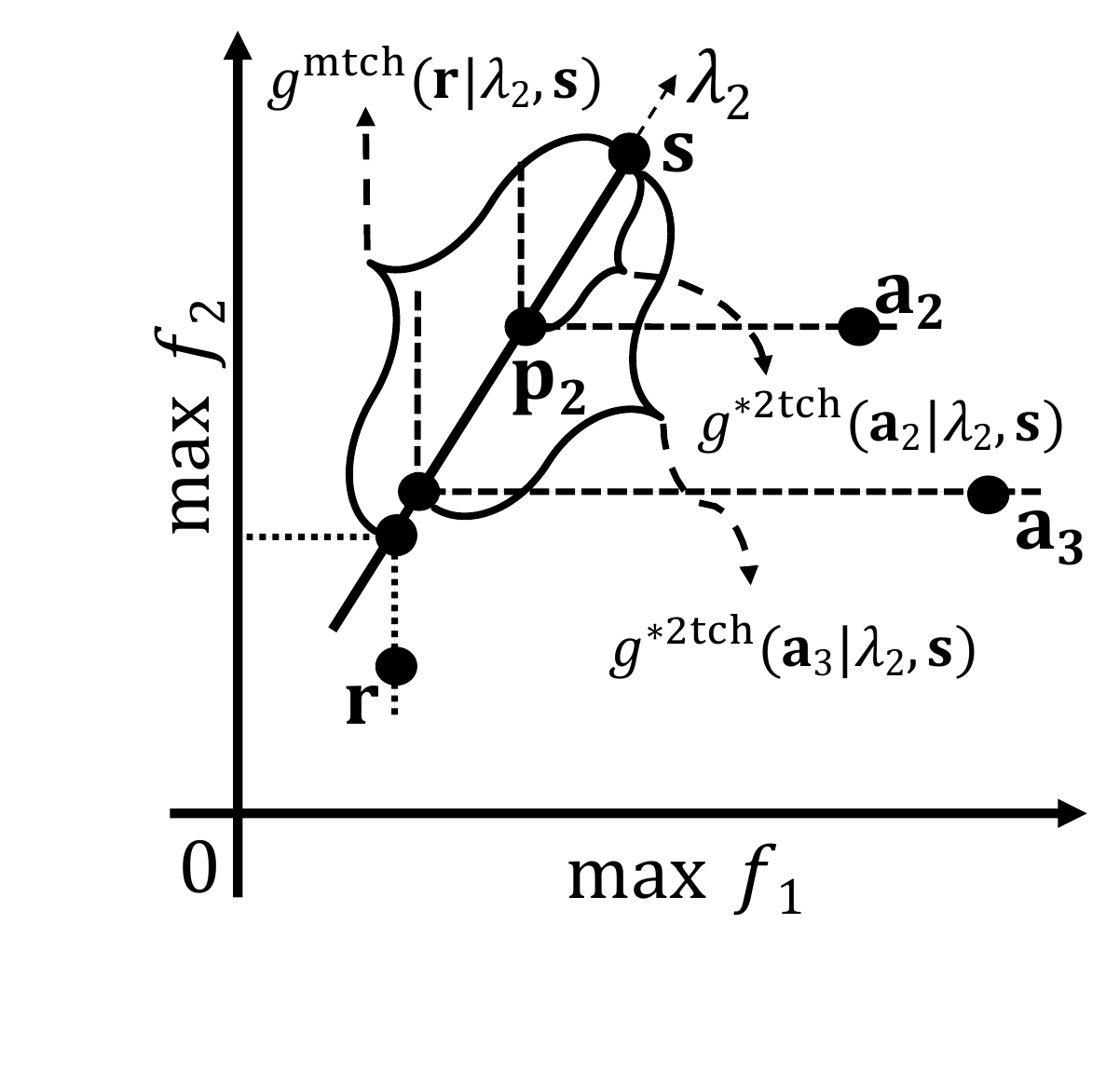}                %以pic.jpg的0.5倍大小输出
\end{minipage}}
\caption{An illustration of the proposed idea (Situation 2).}\label{situation2}                                                        %图片引用标记
\end{figure}

Next, let us consider the situation in Fig.~\ref{situation2} (a). In this situation, the line segments in the hypervolume contribution region of the solution $\mathbf{s}$ not only intersect with the attainment surface of the solution set $A\setminus \{\mathbf{s}\}$ but also intersect with the hypervolume boundary of the solution set $A$ associated with the reference point $\mathbf{r}$. For the line segments intersecting with the attainment surface of the solution set $A\setminus \{\mathbf{s}\}$, the lengths can be calculated by  Eq.~\eqref{length1}. For the line segments intersecting with the hypervolume boundary associated with the reference point $\mathbf{r}$, the lengths are calculated as follows:
\begin{equation}
L\left(\mathbf{s}, \mathbf{r},\lambda\right) =g^{\textrm{mtch}}(\mathbf{r}|\lambda,\mathbf{s}),
\end{equation}
where $g^{\textrm{mtch}}(\mathbf{r}|\lambda,\mathbf{s})$ is defined as:
\begin{equation}
\label{gmtchnew}
g^{\textrm{mtch}}(\mathbf{r}|\lambda,\mathbf{s})=\min_{j\in\{1,...,m\}}\left\{\frac{|s_j-r_j|}{\lambda_j}\right\}.
\end{equation}

Fig.~\ref{situation2} (b) and (c) show the contour lines of $g^{\textrm{mtch}}(\mathbf{r}|\lambda,\mathbf{s})$ in Eq.~\eqref{gmtchnew} and  $g^{*\textrm{2tch}}(\mathbf{a}|\lambda,\mathbf{s})$ in Eq.~\eqref{g2tchnew} respectively. From Fig.~\ref{situation2} (b) for a line segment intersecting with the hypervolume boundary associated with the reference point $\mathbf{r}$, its length is $L\left(\mathbf{s}, \mathbf{r},\lambda_1\right)$ where $L\left(\mathbf{s}, \mathbf{r},\lambda_1\right)<L\left(\mathbf{s},A,\lambda_1\right)$ holds. From Fig.~\ref{situation2} (c) for a line segment intersecting with the attainment surface of the solution set $A\setminus \{\mathbf{s}\}$, its length is $L\left(\mathbf{s},A,\lambda_2\right)$ where $L\left(\mathbf{s},A,\lambda_2\right) < L\left(\mathbf{s}, \mathbf{r},\lambda_2\right)$ holds.

Wrapping up the above two situations, the length of any line segment along the direction $\lambda$ starting from $\mathbf{s}$ and intersecting with the boundaries of the hypervolume contribution region is calculated as follows:
\begin{equation}
\label{thisidea}
L\left(\mathbf{s},A,\mathbf{r},\lambda\right) =\min\{L\left(\mathbf{s},A,\lambda\right),L\left(\mathbf{s},\mathbf{r},\lambda\right)\}.
\end{equation}

\textcolor{black}{Another idea of tackling situation 2 without introducing $g^{\textrm{mtch}}(\mathbf{r}|\lambda,\mathbf{s})$ function is as follows: we can introduce $m$ points which are almost copies of $\mathbf{s}$, except that in one coordinate they actually contain the value in $\mathbf{r}$. Then we can use these $m$ points instead of $\mathbf{r}$ to calculate the length of the line segments using Eq.~\eqref{length1}. This idea can be illustrated in Fig.~\ref{situation2} (a) where two points $\mathbf{s}'$ and $\mathbf{s}''$ can be introduced and treated as solutions. In this manner, we can use Eq.~\eqref{length1} directly to calculate the length of the line segments as:
$L\left(\mathbf{s},A\bigcup \{\mathbf{s}',\mathbf{s}''\},\lambda\right) =\min_{\mathbf{a}\in A\bigcup\{\mathbf{s}',\mathbf{s}''\} \setminus \{\mathbf{s}\}}\left\{g^{*\textrm{2tch}}(\mathbf{a}|\lambda,\mathbf{s})\right\}.$
}

\textcolor{black}{The above idea can treat the calculation in a unified manner. 
However, we need to introduce $m$ points and remove one point (i.e., reference point $\mathbf{r}$) to do the calculation. In general, we have $(m-1)$ more points compared with the original point set. This will lead to extra computational load. So based on the computational load consideration, we will use Eq.~\eqref{thisidea} to calculate the length of the line segments in this letter.}

Based on the lengths of the line segments in the hypervolume contribution region, we can define the following R2 indicator to approximate the hypervolume contribution of a solution $\mathbf{s}$ in solution set $A$ as:
\begin{equation}
\begin{aligned}
\label{newmethod}
&R^{\textrm{HVC}}_2(\mathbf{s},A,\Lambda,\mathbf{r},\alpha)=\frac{1}{|\Lambda|}\sum_{\lambda\in \Lambda}L\left(\mathbf{s},A,\mathbf{r},\lambda\right)^\alpha\\
&=\frac{1}{|\Lambda|}\sum_{\lambda\in \Lambda}\min\left\{\min_{\mathbf{a}\in A\setminus \{\mathbf{s}\}}\left\{g^{*\textrm{2tch}}(\mathbf{a}|\lambda,\mathbf{s})\right\},g^{\textrm{mtch}}(\mathbf{r}|\lambda,\mathbf{s})\right\}^\alpha ,
\end{aligned}
\end{equation}
where $\alpha=1$ if we want to use the average length of the line segments for the approximation, or $\alpha=m$ if we want to use the average $m$ powered length of the line segments for the approximation.

Comparing with the traditional method described in the previous subsection, the proposed new method for the hypervolume contribution approximation has the following advantages:
\begin{enumerate}
\item The hypervolume contribution approximation of a solution can be directly calculated by Eq.~\eqref{newmethod}.
\item The new method will utilize all direction vectors only in the hypervolume contribution region, while the traditional method utilizes them across the entire hypervolume region. For this reason, the approximation accuracy of the new method could be much higher than the traditional method when the same number of direction vectors are used in each method.
\end{enumerate}

\section{Numerical studies}
\label{numstudy}
\subsection{Experiment settings}
\subsubsection{Solution sets generation}
In our experiments, we examine six types of Pareto front (PF): triangular PF (including linear, convex and concave) and inverted triangular PF (including linear, convex and concave). For the triangular PF, the linear PF is expressed as $f_1+f_2+...+f_m = 1$ and $f_i\geq 0$ for $i = 1,...,m$, the concave PF is expressed as $f_1^2+f_2^2+...+f_m^2=1$ and $f_i\geq 0$ for $i = 1,...,m$, the convex PF is expressed as $\sqrt{f_1}+\sqrt{f_2}+...+\sqrt{f_m}=1$ and $f_i\geq 0$ for  $i = 1,...,m$. For the inverted triangular PF, it is obtained by transforming each point of the triangular PF by $1-f_i$ for $i = 1,...,m$.

We examine 5- and 10-dimension cases (i.e., $m=5,10$). \textcolor{black}{For each case, $N$ solutions are randomly sampled from each PF to form solution sets with different PF shapes. This sampling procedure is repeated to generate 100 solution sets for each PF shape. In order to examine the effect of the number of solutions on the performance of each method, we choose $N=100,200,...,500$.}
\textcolor{black}{
\subsubsection{Reference point specification}
For the reference point $\mathbf{r}=(r,...,r)$ (i.e., each element in $\mathbf{r}$ has the same value, denoted as $r$), we examine 5 different specifications: $r=0.0,-0.1,...,-0.4$.}
\subsubsection{Compared methods}
Three methods for the hypervolume contribution approximation are compared. The first method is the traditional method described in Section \ref{traditionmethod}. The second method is the new method proposed in Section \ref{newmethod1}. The third method is the Monte Carlo sampling method proposed in \cite{bader2010faster}. 

For the traditional method, we use $R^{\textrm{HV}}_2$ (Eq.~\eqref{newr2}) as the R2 indicator and approximate the hypervolume contribution according to the R2 contribution (Eq.~\eqref{r2c}). For the new method, we use $R^{\textrm{HVC}}_2$ (Eq.~\eqref{newmethod}) for the approximation and choose $\alpha=m$. The direction vectors used in the R2 indicator are generated by uniformly sampling points on the unit hypersphere\footnote{First we randomly sample points $\mathbf{x}$ according to the normal distribution $\mathcal{N}_m(0,I_m)$, then the corresponding direction vectors are obtained by $\lambda = |\mathbf{x}|/\left \| \mathbf{x} \right \|_2$.}, while the sampling points used in the Monte Carlo method are uniformly sampled in the sampling space.

In order to investigate the effect of the number of the direction vectors and the sampling points on the performance of the three methods, we examine ten settings: $100,200,...,1000$ direction vectors and sampling points.

\textcolor{black}{In addition to the three approximation methods mentioned above, we also compare two exact hypervolume contribution calculation methods: IWFG \cite{while2012applying} and exQHV \cite{russo2016extending}. IWFG is a method to identify the solution with the smallest hypervolume contribution in a solution set, while exQHV is a method to calculate the hypervolume contributions of all solutions in a solution set.}
\textcolor{black}{
\subsubsection{Platforms and Implementations}
We conduct the experiments on a PC equipped with Intel Core i7-8700K CPU@3.70GHz, 16GB RAM and Windows 10 Operating System. The three approximation methods are implemented by ourselves in MATLAB R2018b. The two exact methods are based on their available source code\footnote{IWFG from \url{http://www.wfg.csse.uwa.edu.au/hypervolume/#code}, exQHV from \url{http://web.tecnico.ulisboa.pt/luis.russo/QHV/#down}.} which are both written in C. We compile the source code with \texttt{gcc} 7.4.0 in Cygwin Version 2.11.2-1.}

\subsection{Performance metrics}
\textcolor{black}{Three} performance metrics are used to evaluate the performance of the compared methods for the hypervolume contribution approximation \textcolor{black}{and calculation}. 

The first metric is the consistency rate in the order of the solution pairs in the solution set between their true hypervolume contributions and their approximated hypervolume contributions. Given a solution set $A$, there are a total number of $\binom{|A|}{2}$ solution pairs. For two solutions $\mathbf{s}_1$ and $\mathbf{s}_2$ from $A$, if the orders of their true hypervolume contributions and their approximated hypervolume contributions are the same, then the solution pair $(\mathbf{s}_1,\mathbf{s}_2)$ is called a consistent pair. \textcolor{black}{The consistency rate of a solution set $A$ is the ratio of the total number of the consistent pairs with $\binom{|A|}{2}$}. 

The second metric is the correct identification rate of the worst solution with the smallest hypervolume contribution in each solution set over all solution sets. For a solution set, if the solution with the smallest approximated hypervolume contribution also has the smallest true hypervolume contribution, then we call it as a correct identification for this solution set. The correct identification rate is the ratio of the total number of the correct identifications to the total number of the solution sets. 

\textcolor{black}{Our choice of these two metrics is based on practical considerations. As discussed in Section \ref{introduction}, in hypervolume-based EMOAs the solution with the smallest hypervolume contribution is removed from the population. So the relationship between the hypervolume contributions of two solutions is more important than their values. The two metrics are able to evaluate such a relationship, especially the second metric is able to evaluate the ability of a method to identify the solution with the smallest hypervolume contribution in a solution set.}

\begin{figure}[!htbp]
\centering                                                          %居中
\subfigure{       
\hspace{-0.5cm}             %第一张子图
\begin{minipage}{0.45\columnwidth}\centering                                                          %子图居中
\includegraphics[scale=0.23]{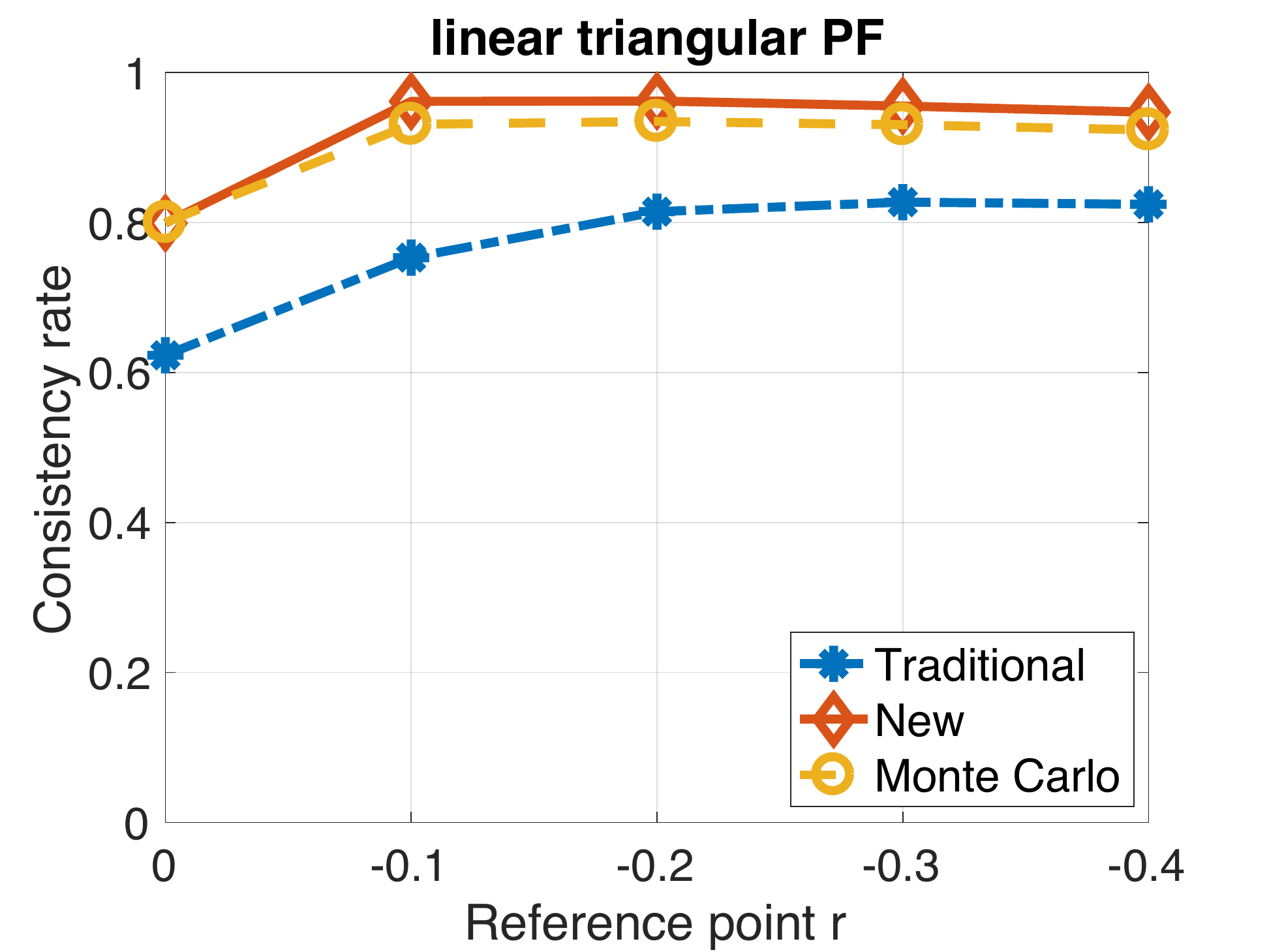}               %以pic.jpg的0.5倍大小输出
\end{minipage}}
\subfigure{                    %第二张子图
\begin{minipage}{0.45\columnwidth}\centering                                                          %子图居中
\includegraphics[scale=0.23]{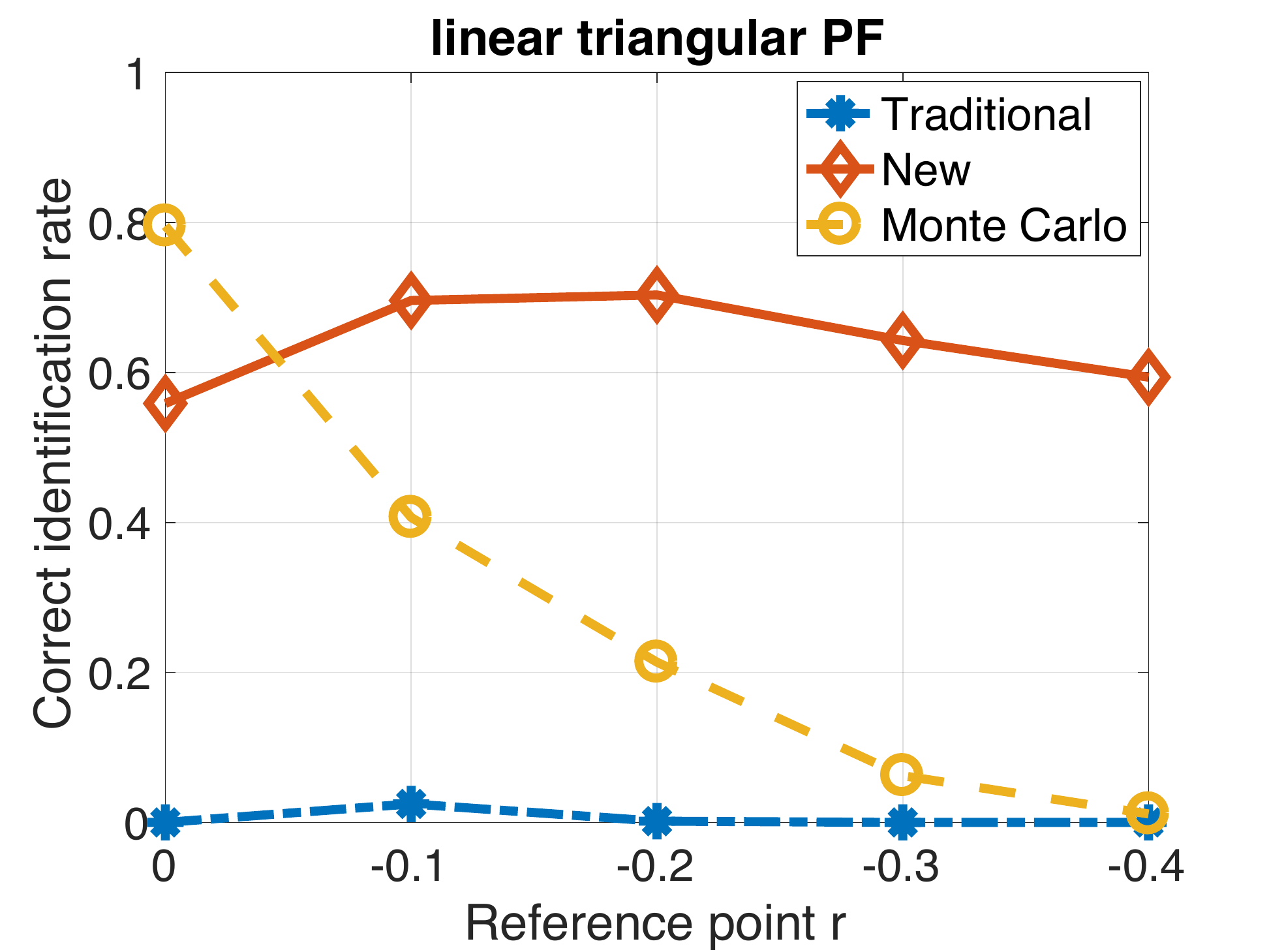}                %以pic.jpg的0.5倍大小输出
\end{minipage}}
\subfigure{      
\hspace{-0.5cm}              %第一张子图
\begin{minipage}{0.45\columnwidth}\centering                                                          %子图居中
\includegraphics[scale=0.23]{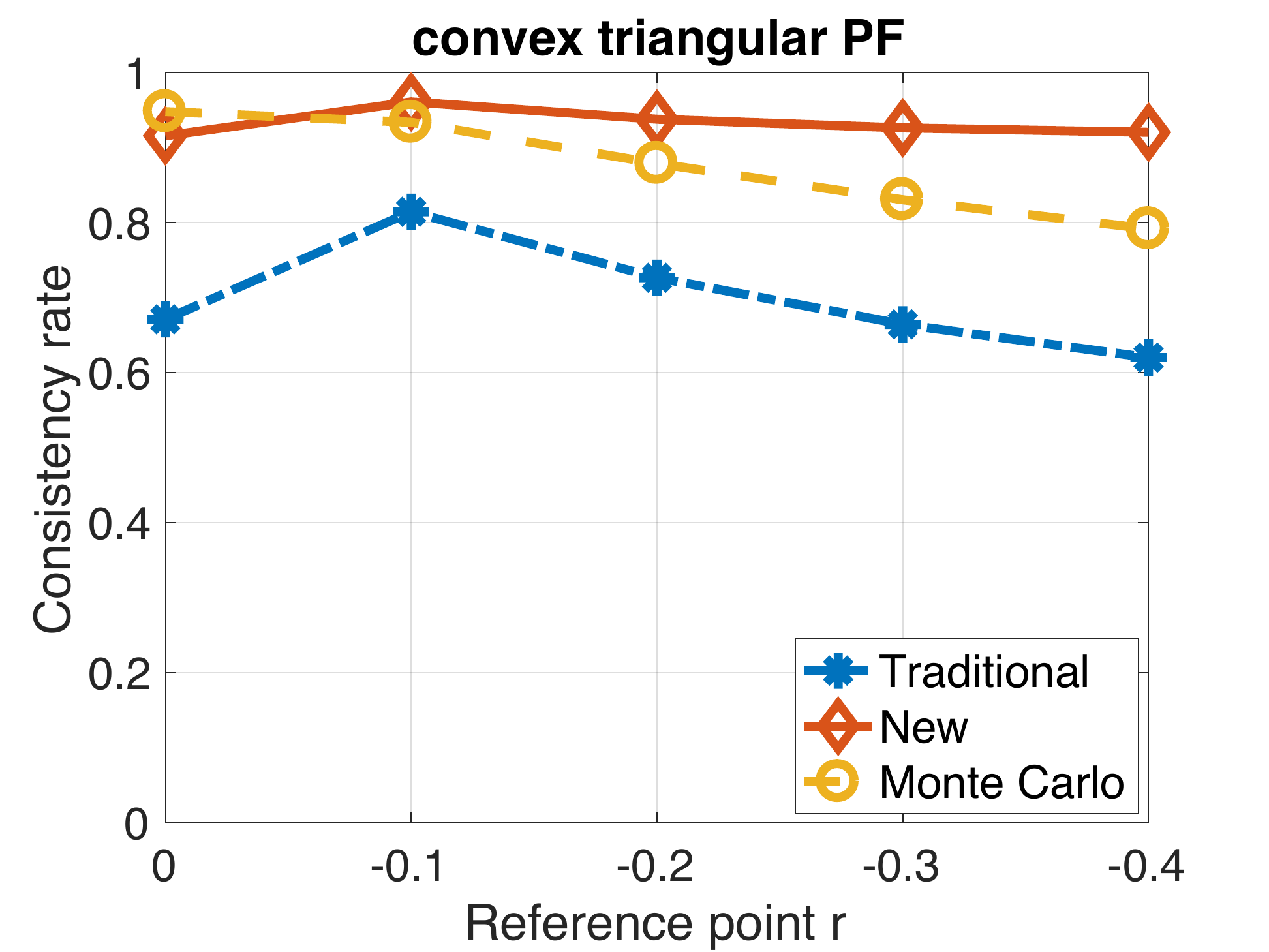}               %以pic.jpg的0.5倍大小输出
\end{minipage}}
\subfigure{                    %第二张子图
\begin{minipage}{0.45\columnwidth}\centering                                                          %子图居中
\includegraphics[scale=0.23]{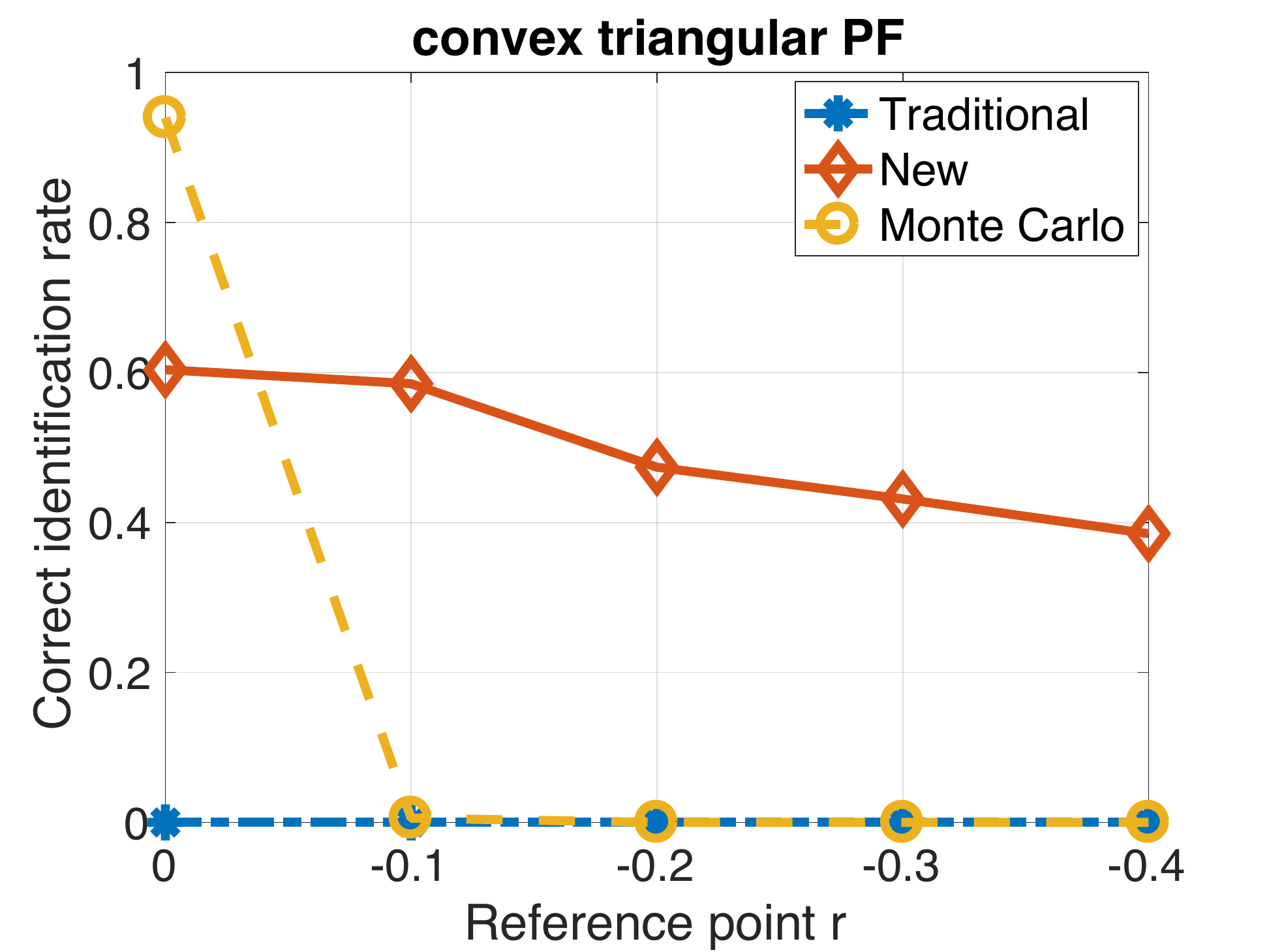}                %以pic.jpg的0.5倍大小输出
\end{minipage}}
\subfigure{      
\hspace{-0.5cm}              %第一张子图
\begin{minipage}{0.45\columnwidth}\centering                                                          %子图居中
\includegraphics[scale=0.23]{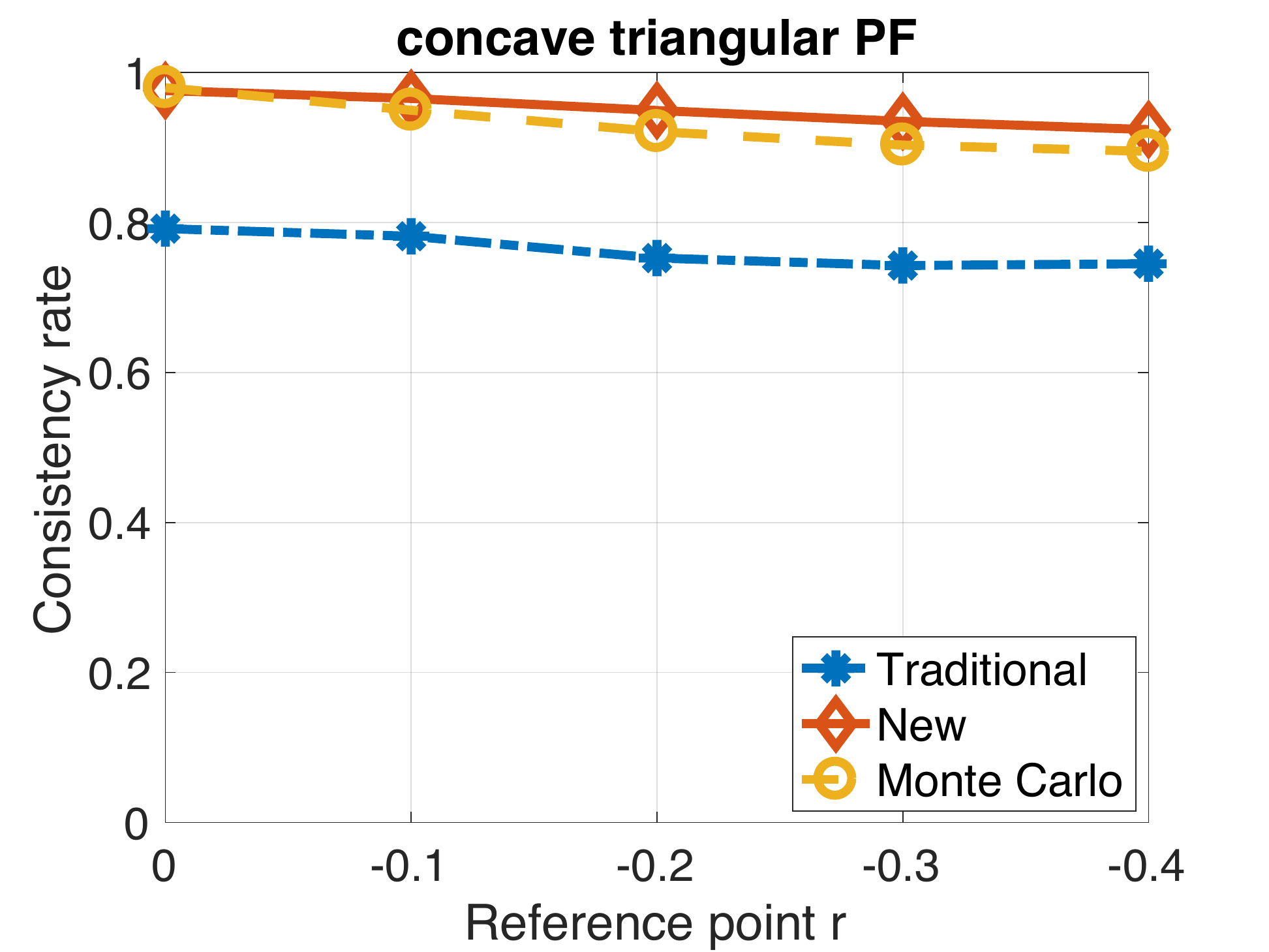}               %以pic.jpg的0.5倍大小输出
\end{minipage}}
\subfigure{                    %第二张子图
\begin{minipage}{0.45\columnwidth}\centering                                                          %子图居中
\includegraphics[scale=0.23]{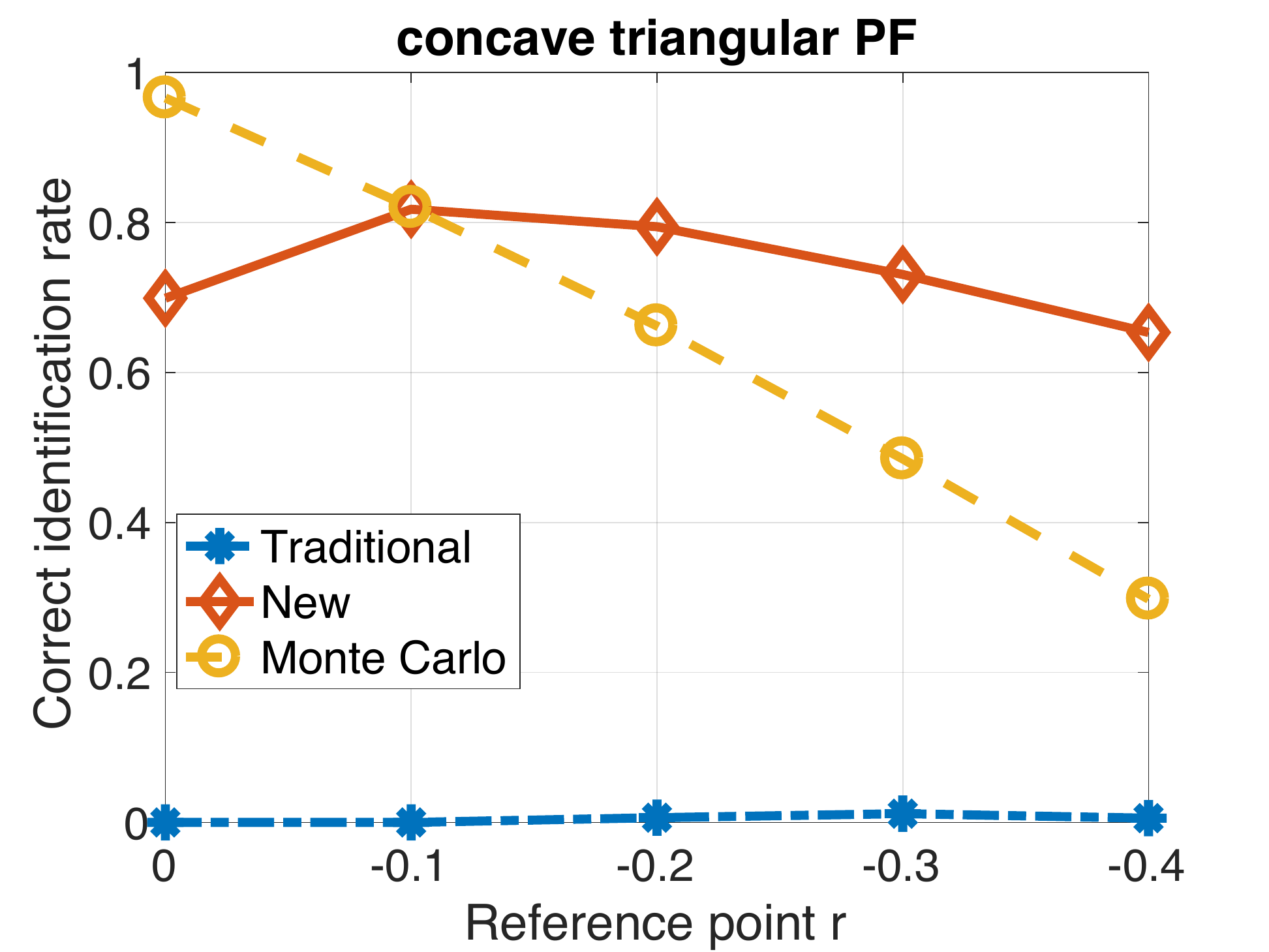}                %以pic.jpg的0.5倍大小输出
\end{minipage}}
\caption{The approximation accuracy with respect to the reference point on the triangular PF solution sets.} %                         %大图名称
\label{ref1}                                                        %图片引用标记
\end{figure}

\begin{figure}[!htbp]
\centering                                                          %居中
\subfigure{   
\hspace{-0.5cm}                 %第一张子图
\begin{minipage}{0.45\columnwidth}\centering                                                          %子图居中
\includegraphics[scale=0.23]{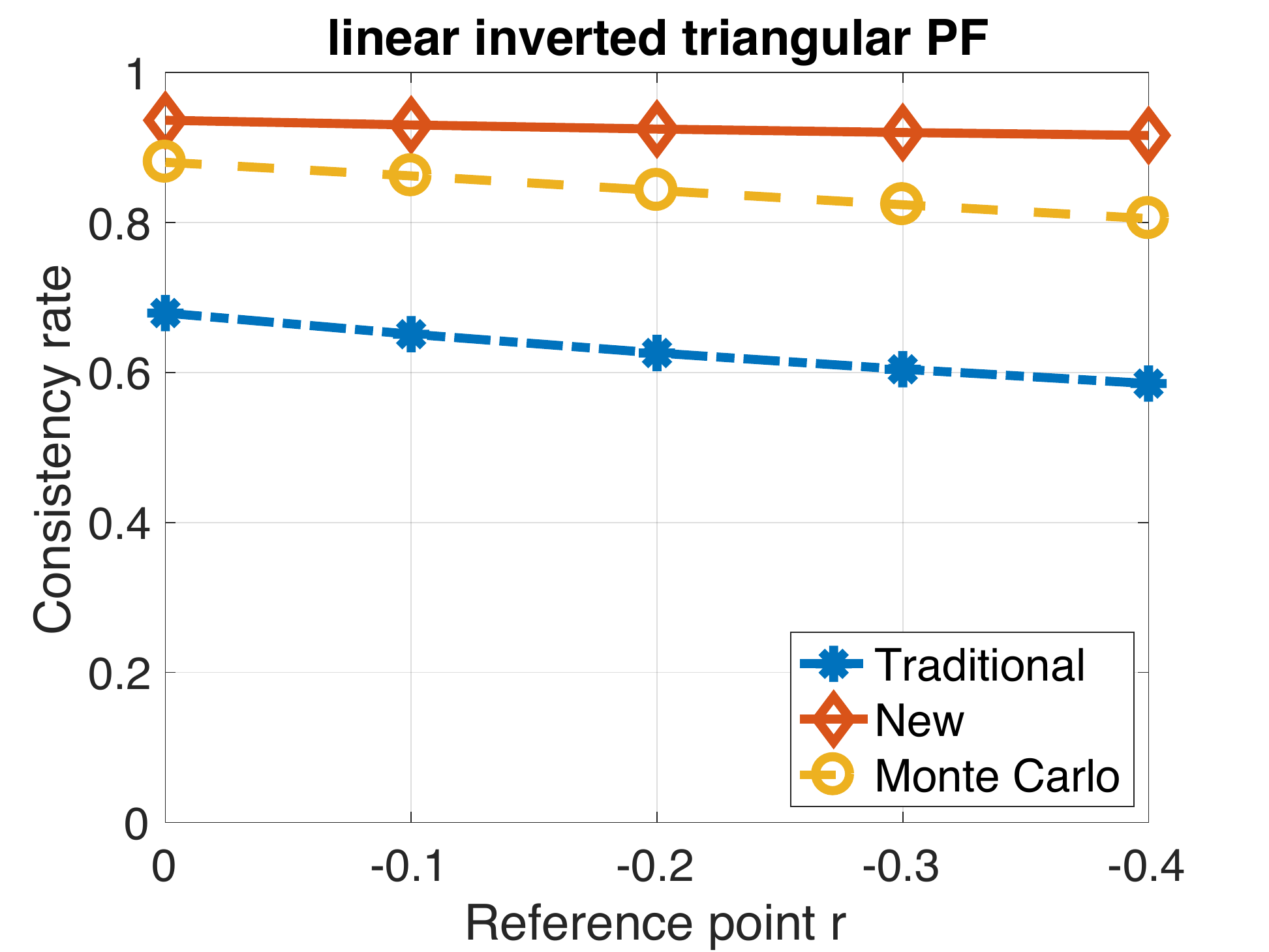}               %以pic.jpg的0.5倍大小输出
\end{minipage}}
\subfigure{                    %第二张子图
\begin{minipage}{0.45\columnwidth}\centering                                                          %子图居中
\includegraphics[scale=0.23]{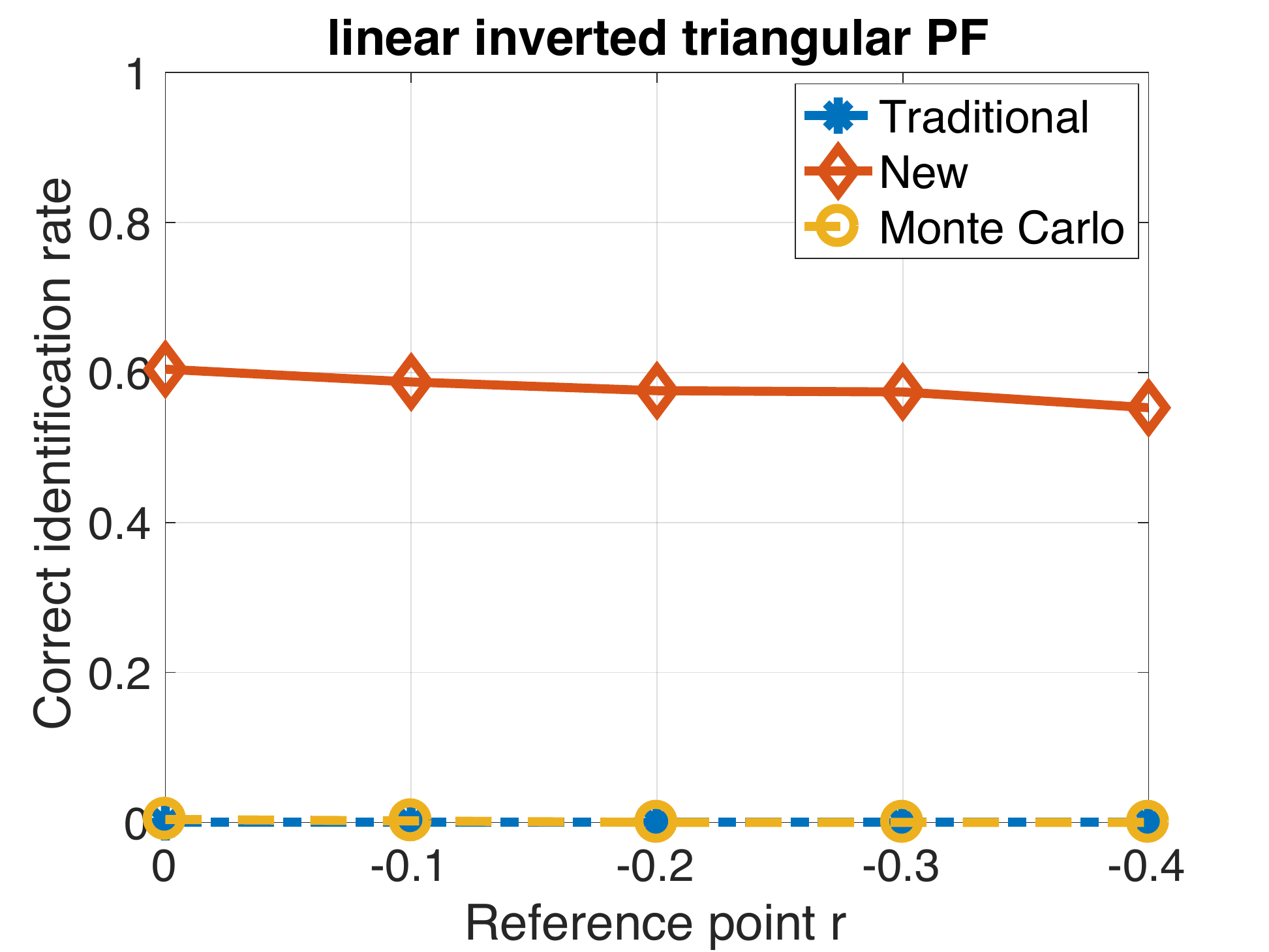}                %以pic.jpg的0.5倍大小输出
\end{minipage}}
\subfigure{     
\hspace{-0.5cm}               %第一张子图
\begin{minipage}{0.45\columnwidth}\centering                                                          %子图居中
\includegraphics[scale=0.23]{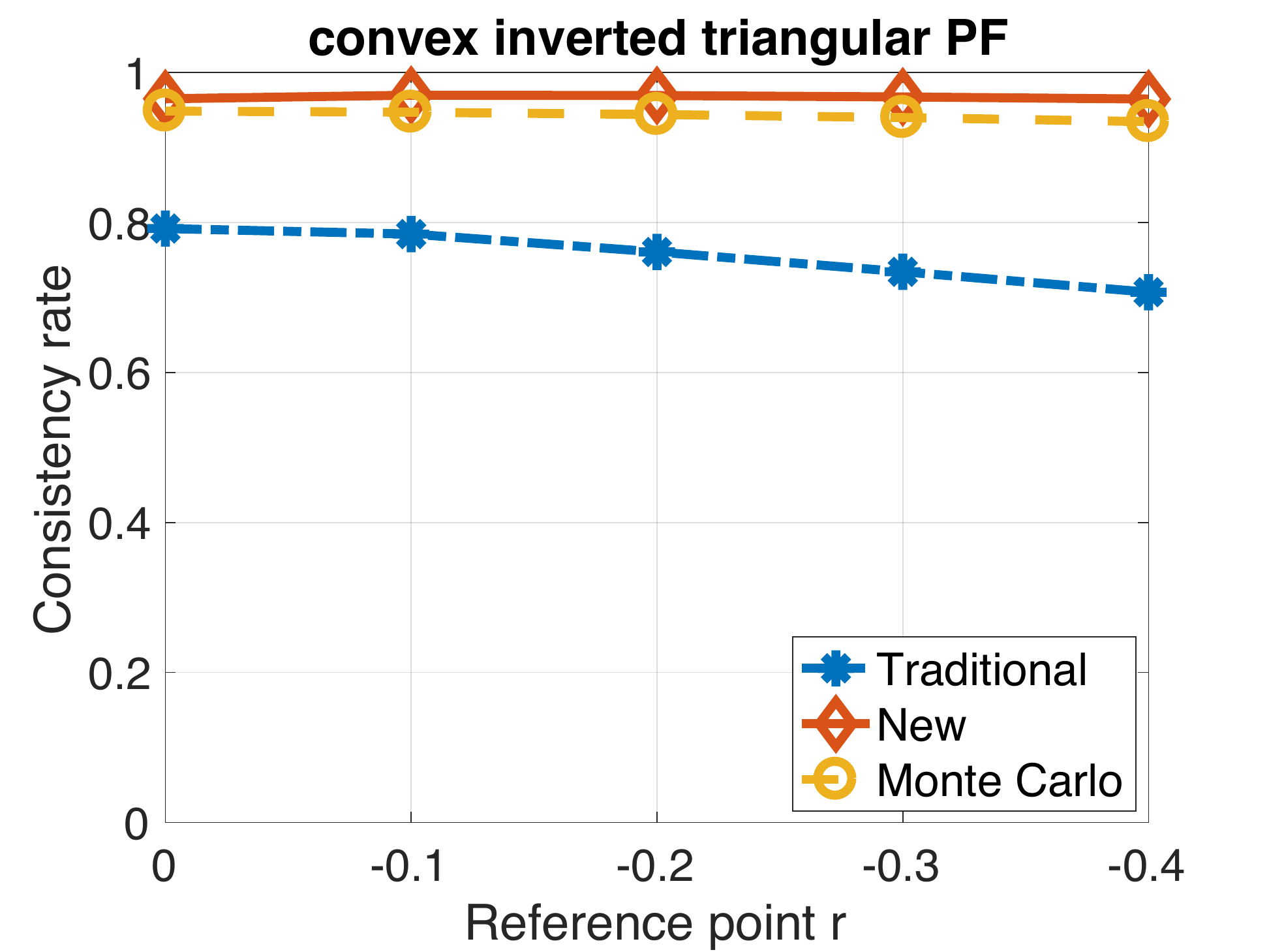}               %以pic.jpg的0.5倍大小输出
\end{minipage}}
\subfigure{                    %第二张子图
\begin{minipage}{0.45\columnwidth}\centering                                                          %子图居中
\includegraphics[scale=0.23]{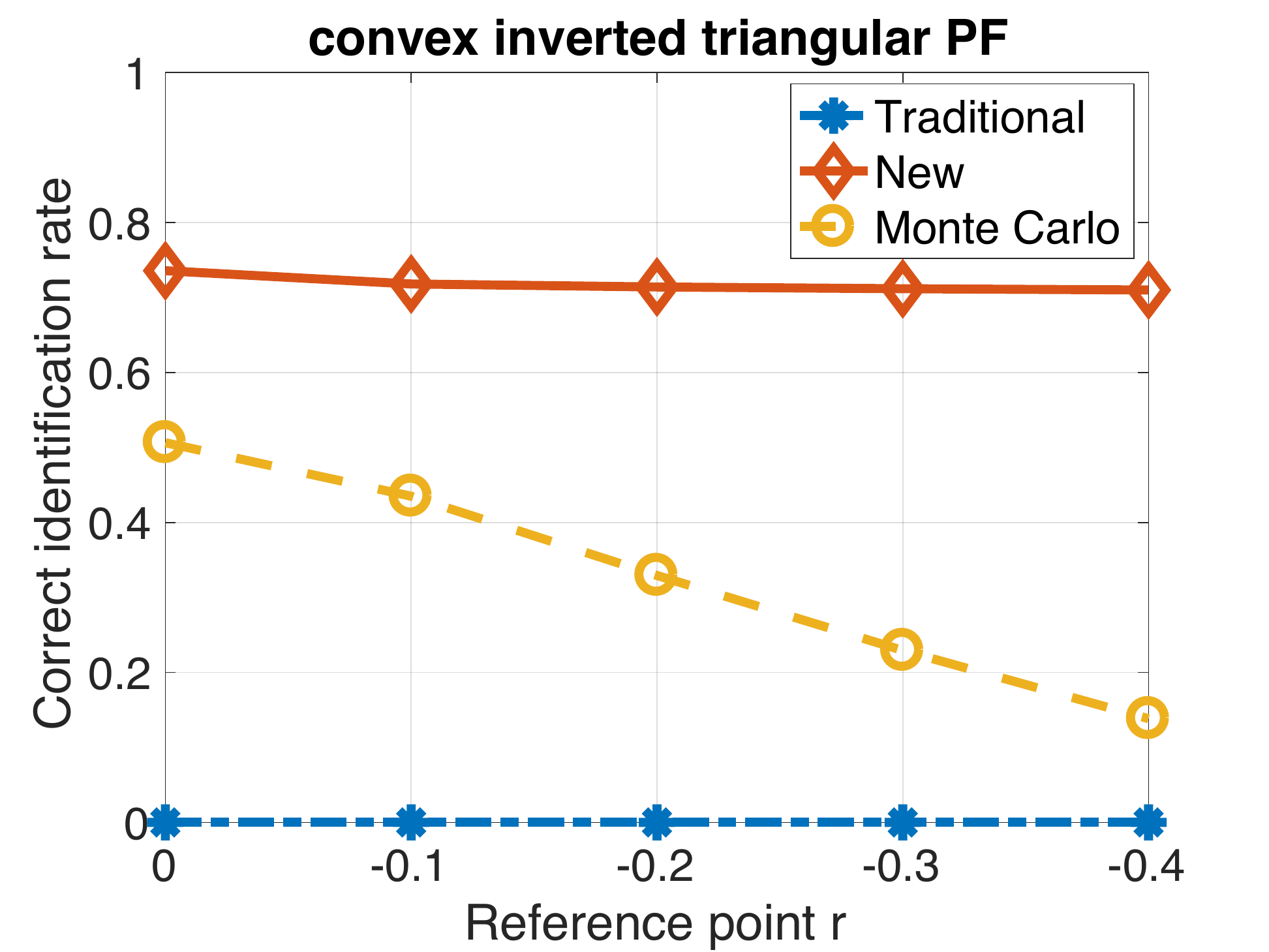}                %以pic.jpg的0.5倍大小输出
\end{minipage}}
\subfigure{       
\hspace{-0.5cm}             %第一张子图
\begin{minipage}{0.45\columnwidth}\centering                                                          %子图居中
\includegraphics[scale=0.23]{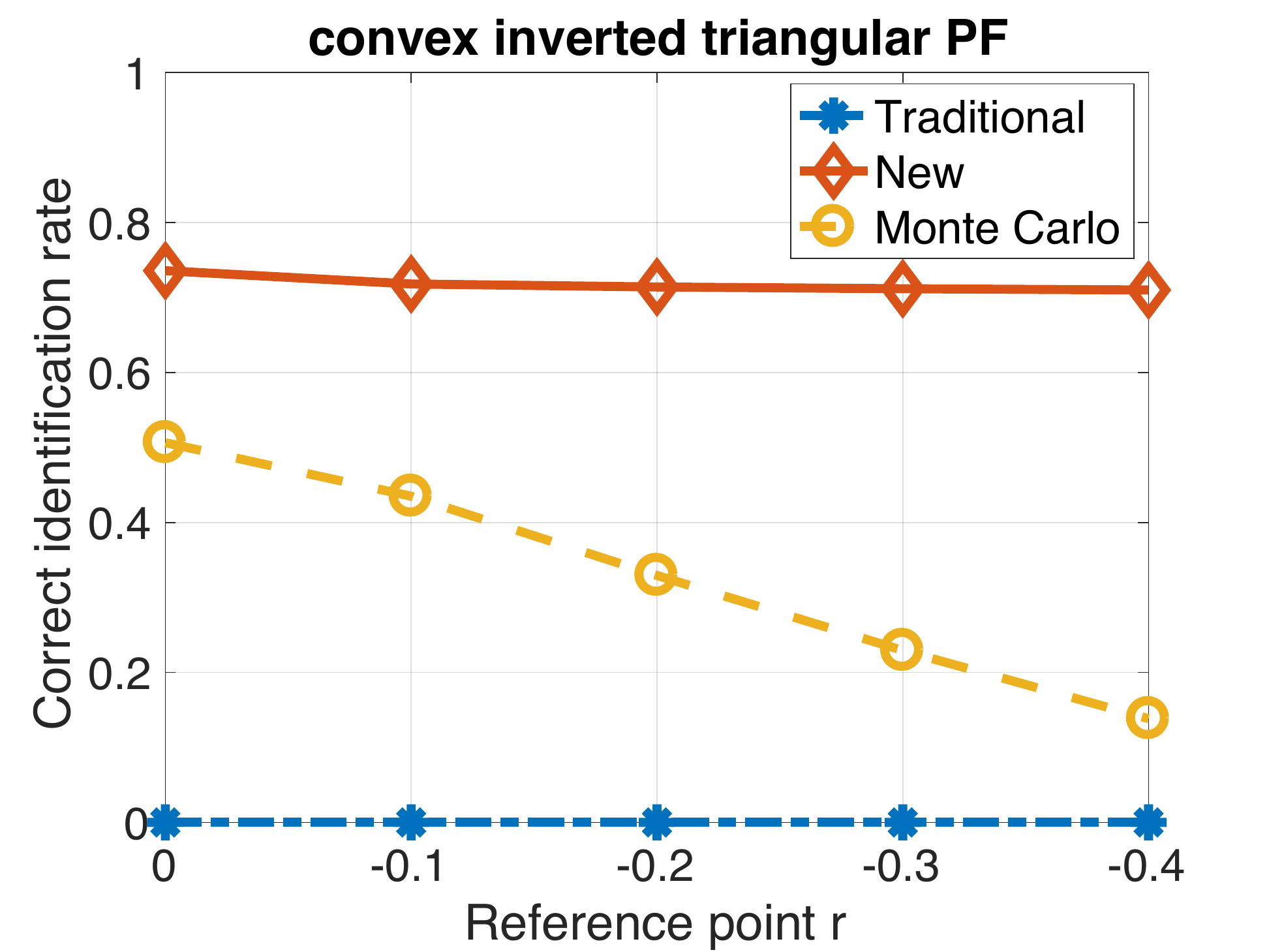}               %以pic.jpg的0.5倍大小输出
\end{minipage}}
\subfigure{                    %第二张子图
\begin{minipage}{0.45\columnwidth}\centering                                                          %子图居中
\includegraphics[scale=0.23]{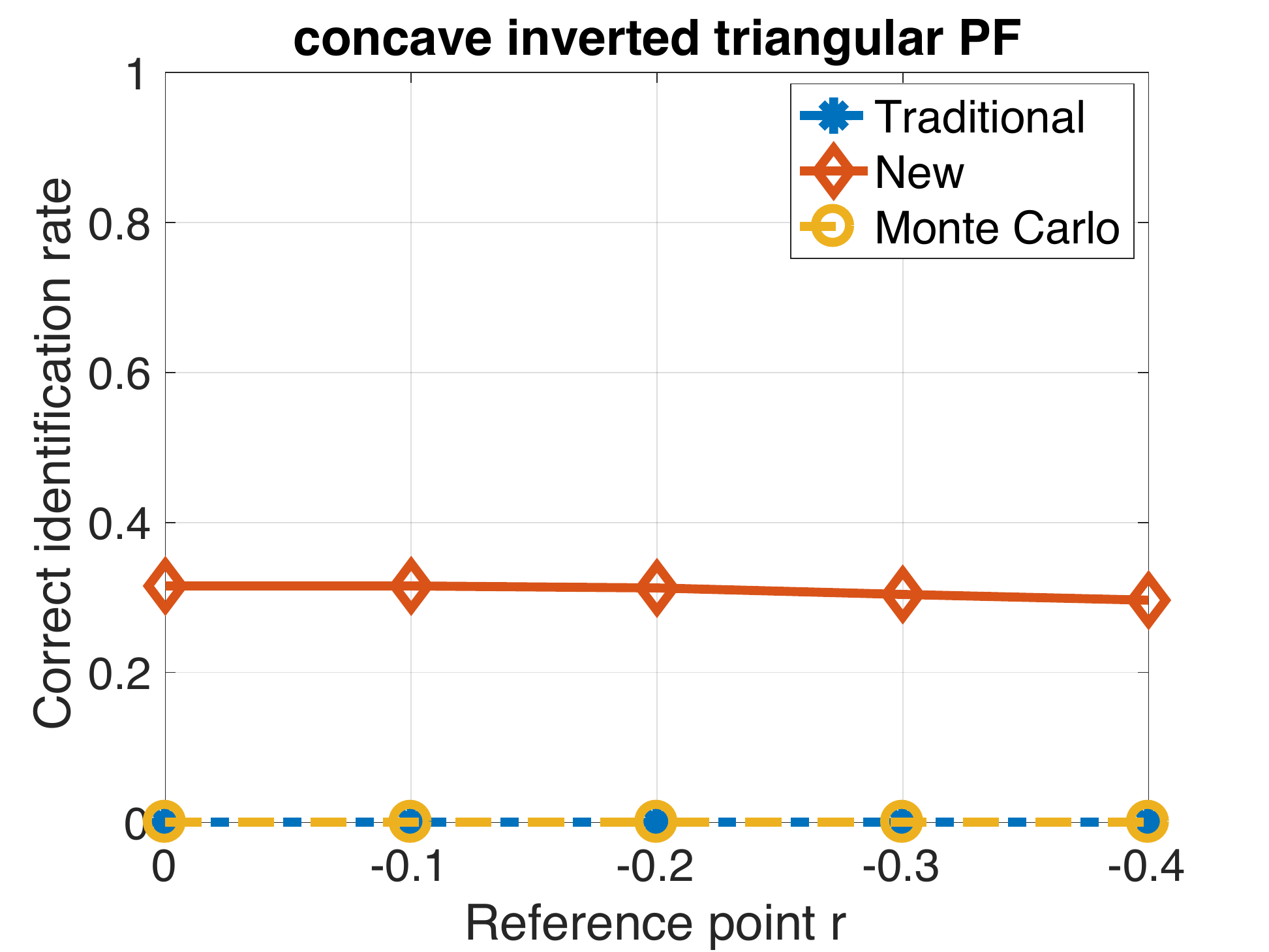}                %以pic.jpg的0.5倍大小输出
\end{minipage}}
\caption{The approximation accuracy with respect to the reference point on the inverted triangular PF solution sets.} %                         %大图名称
\label{ref2}                                                        %图片引用标记
\end{figure}

\textcolor{black}{The above two metrics are only used for evaluating the approximation accuracy of the three approximation methods. In addition to these two metrics, we also compare the runtime of all methods to evaluate their computational efficiency for the hypervolume contribution approximation and calculation.}
\textcolor{black}{
\subsection{Approximation accuracy comparison}
\subsubsection{The effect of the reference point}
First let us examine the effect of the reference point on the performance of the three approximation methods. We fix the number of the direction vectors and the sampling points to 500, the number of solutions in each solution set $N=100$. The results on 5-dimension solution sets are shown in Fig.~\ref{ref1}-\ref{ref2}. All the results shown in the figures are the average of 30 independent runs with different seeds in the generation of the direction vectors and the sampling points\footnote{This applies to all the experiments in this letter.}. }

For the triangular PF solution sets (see Fig.~\ref{ref1}), we can observe that the new method always outperforms the traditional method in terms of both performance metrics. When the reference point $r=0.0$, the Monte Carlo method shows the best performance, while its performance deteriorates dramatically as the reference point changes from $0.0$ to $-0.4$. The new method shows a robust performance with respect to the specification of the reference point. Its performance is worse than the Monte Carlo method only when $r=0.0$, while it outperforms the Monte Carlo method with other reference point specifications. 

For the inverted triangular PF solution sets (see Fig.~\ref{ref2}), we can observe that the new method always outperforms the traditional method and the Monte Carlo method in terms of both performance metrics. We can also see that the new method achieves a robust performance with respect to the specification of the reference point. 

The reason why the Monte Carlo method achieves a good performance on the triangular PF solution sets when $r=0.0$ (while its performance deteriorates dramatically with other reference point specifications) can be explained as follows: If $r=0.0$, \textcolor{black}{which means the reference point $\mathbf{r}$ is the nadir point of the PF}, all solutions on the PF boundaries will have zero hypervolume contributions \cite{Ishibuchi:2017}. So the solution with the smallest hypervolume contribution tends to lie close to the PF boundaries. In this case, the solution with the smallest hypervolume contribution tends to have a small sampling space compared with other solutions, \textcolor{black}{which makes it easy for the Monte Carlo method to identify it}. If $r=-0.1,...,-0.4$, \textcolor{black}{which means the reference point $\mathbf{r}$ is worse than the nadir point of the PF,} all solutions on the PF boundaries will have nonzero hypervolume contributions and their hypervolume contributions will increase as the reference point moves further \cite{Ishibuchi:2017}. As a result, the solution with the smallest hypervolume contribution tends to lie inside of the PF. In this case, the sampling space for the solution with the smallest hypervolume contribution tends to have similar size to other solutions, \textcolor{black}{which makes it difficult for the Monte Carlo method to identify it}. We give some detailed examples to clearly explain this phenomenon in the supplementary material.

\begin{figure}[!htbp]
\centering                                                          %居中
\subfigure{        
\hspace{-0.5cm}            %第一张子图
\begin{minipage}{0.45\columnwidth}\centering                                                          %子图居中
\includegraphics[scale=0.23]{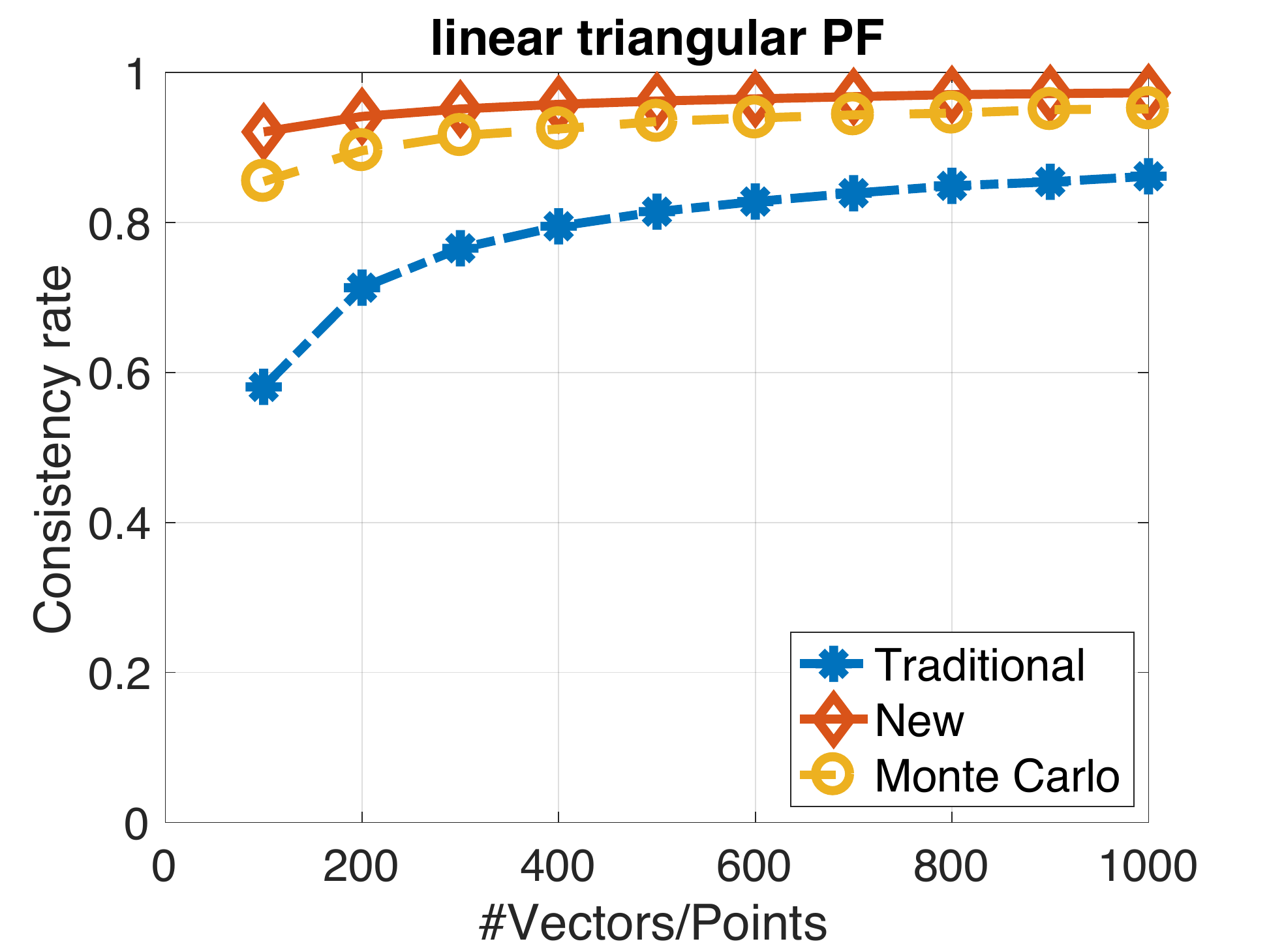}               %以pic.jpg的0.5倍大小输出
\end{minipage}}
\subfigure{                    %第二张子图
\begin{minipage}{0.45\columnwidth}\centering                                                          %子图居中
\includegraphics[scale=0.23]{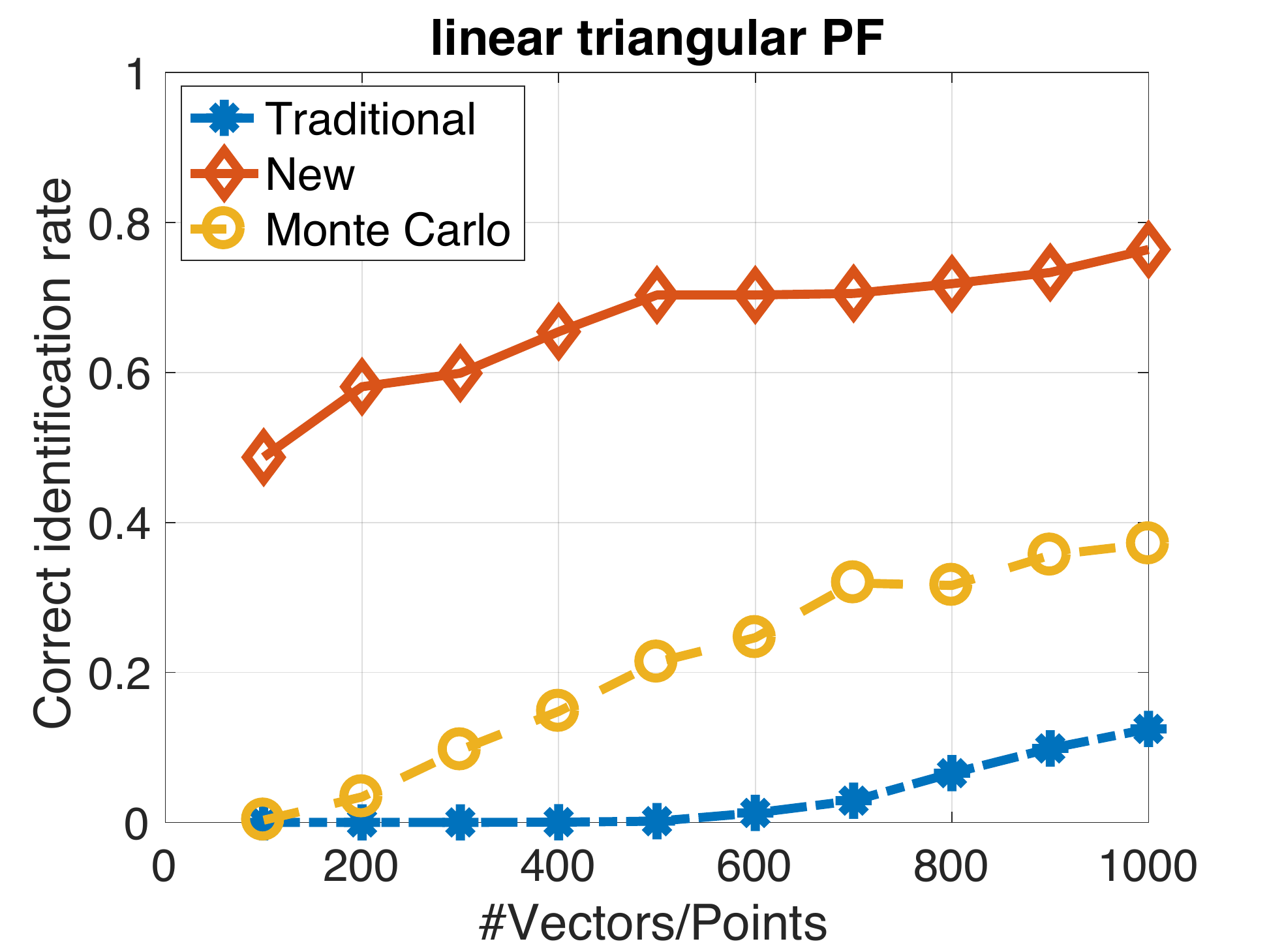}                %以pic.jpg的0.5倍大小输出
\end{minipage}}
\subfigure{        
\hspace{-0.5cm}            %第一张子图
\begin{minipage}{0.45\columnwidth}\centering                                                          %子图居中
\includegraphics[scale=0.23]{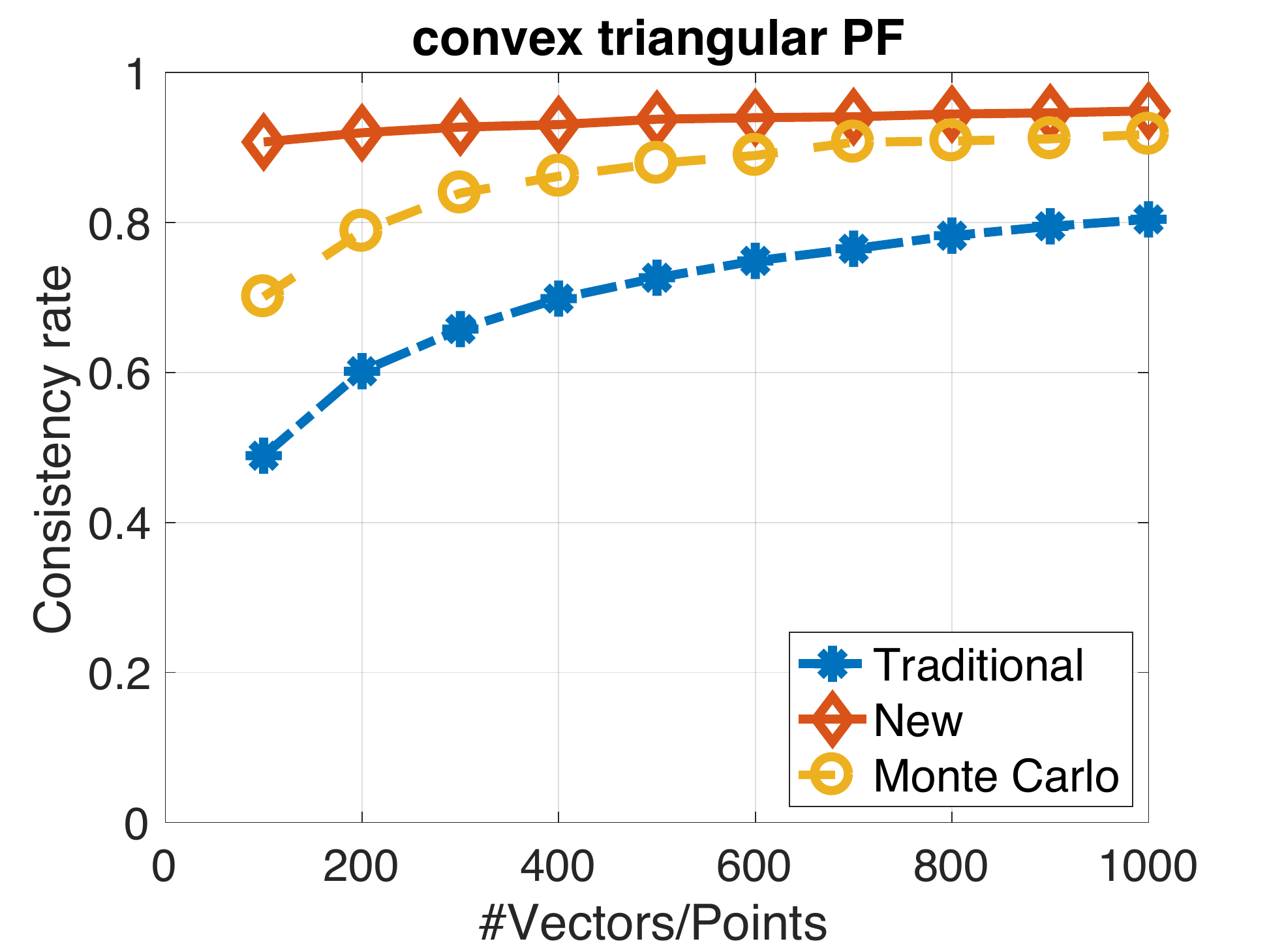}               %以pic.jpg的0.5倍大小输出
\end{minipage}}
\subfigure{                    %第二张子图
\begin{minipage}{0.45\columnwidth}\centering                                                          %子图居中
\includegraphics[scale=0.23]{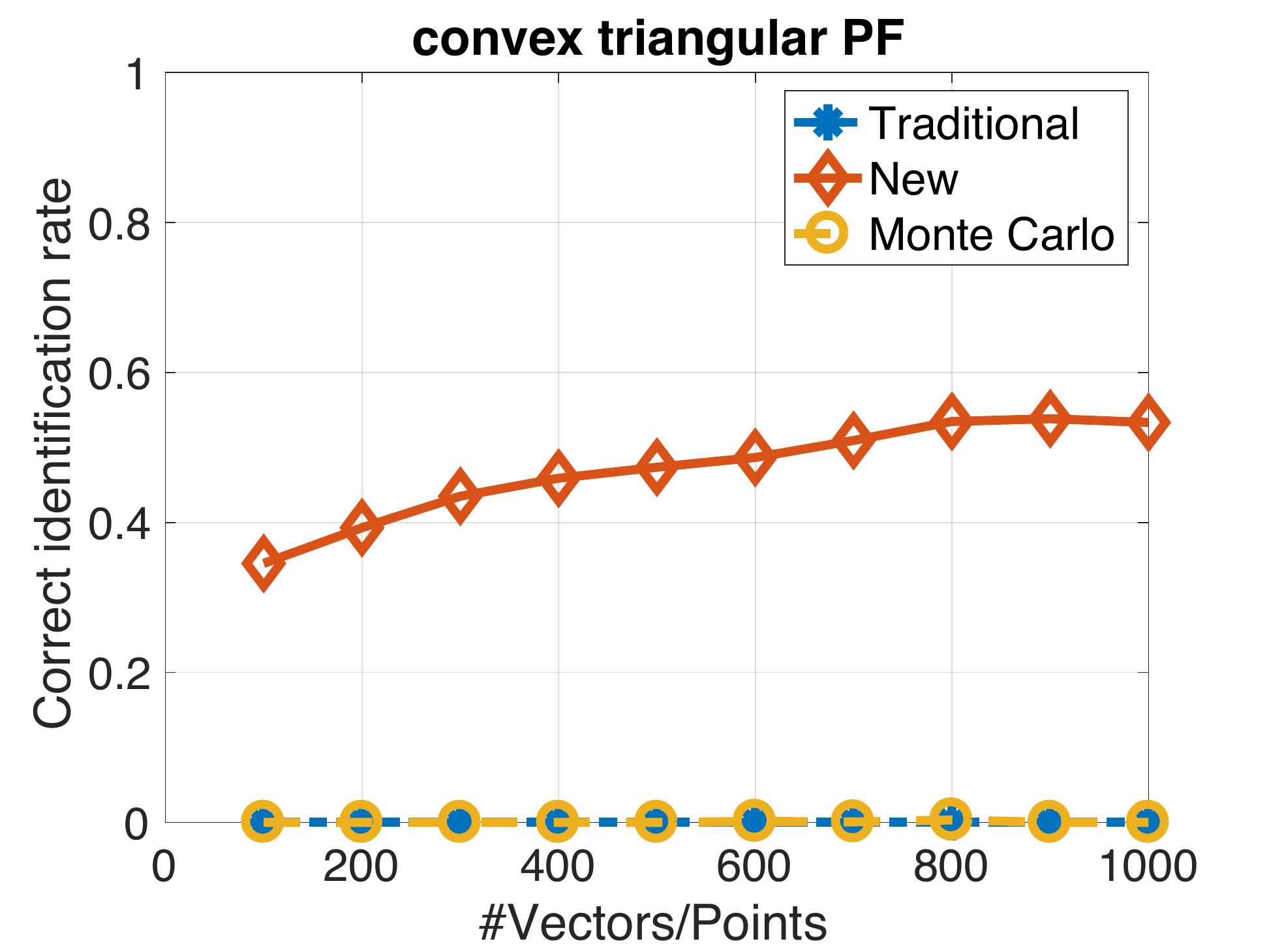}                %以pic.jpg的0.5倍大小输出
\end{minipage}}
\subfigure{ 
\hspace{-0.5cm}                   %第一张子图
\begin{minipage}{0.45\columnwidth}\centering                                                          %子图居中
\includegraphics[scale=0.23]{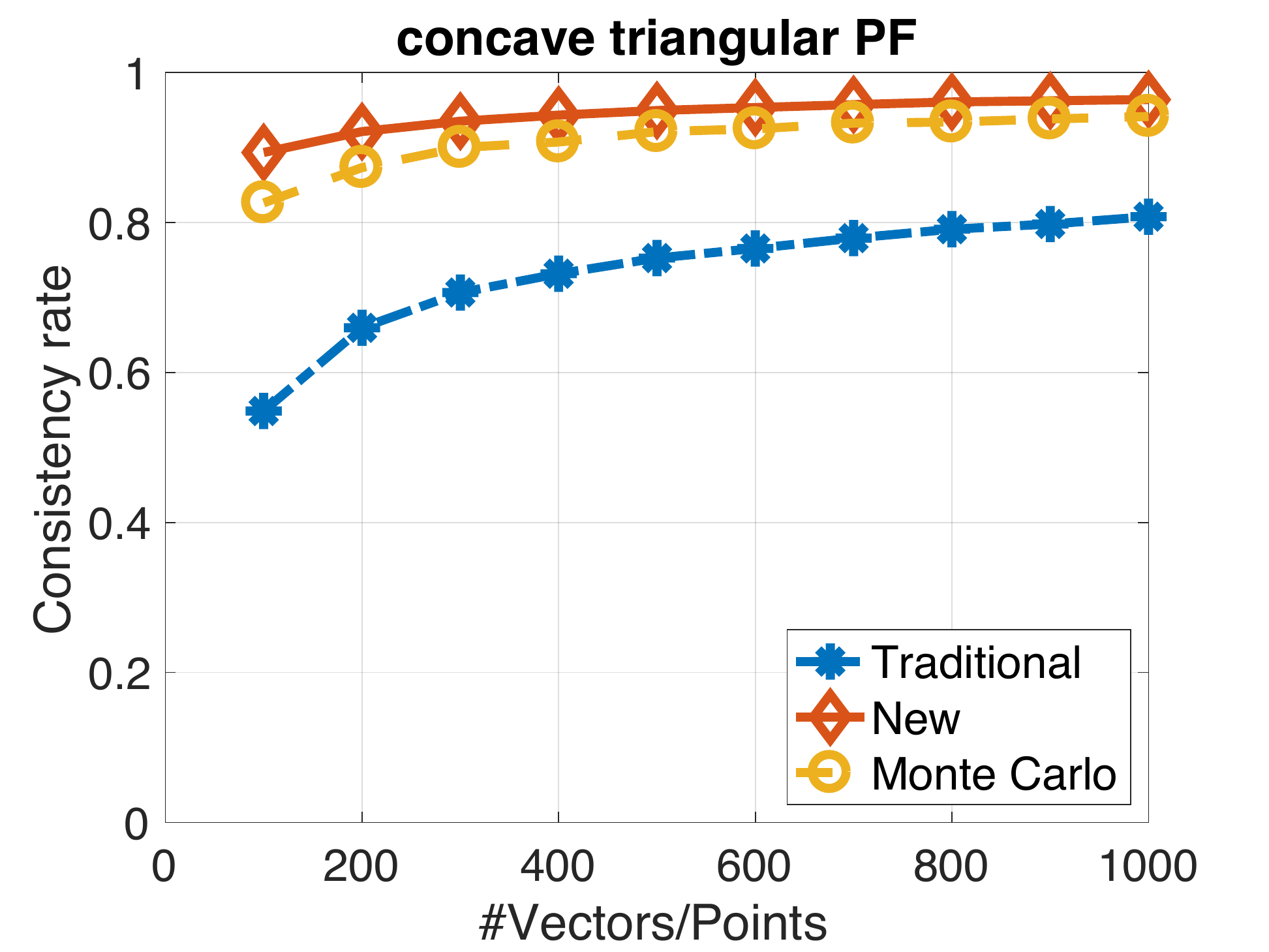}               %以pic.jpg的0.5倍大小输出
\end{minipage}}
\subfigure{                    %第二张子图
\begin{minipage}{0.45\columnwidth}\centering                                                          %子图居中
\includegraphics[scale=0.23]{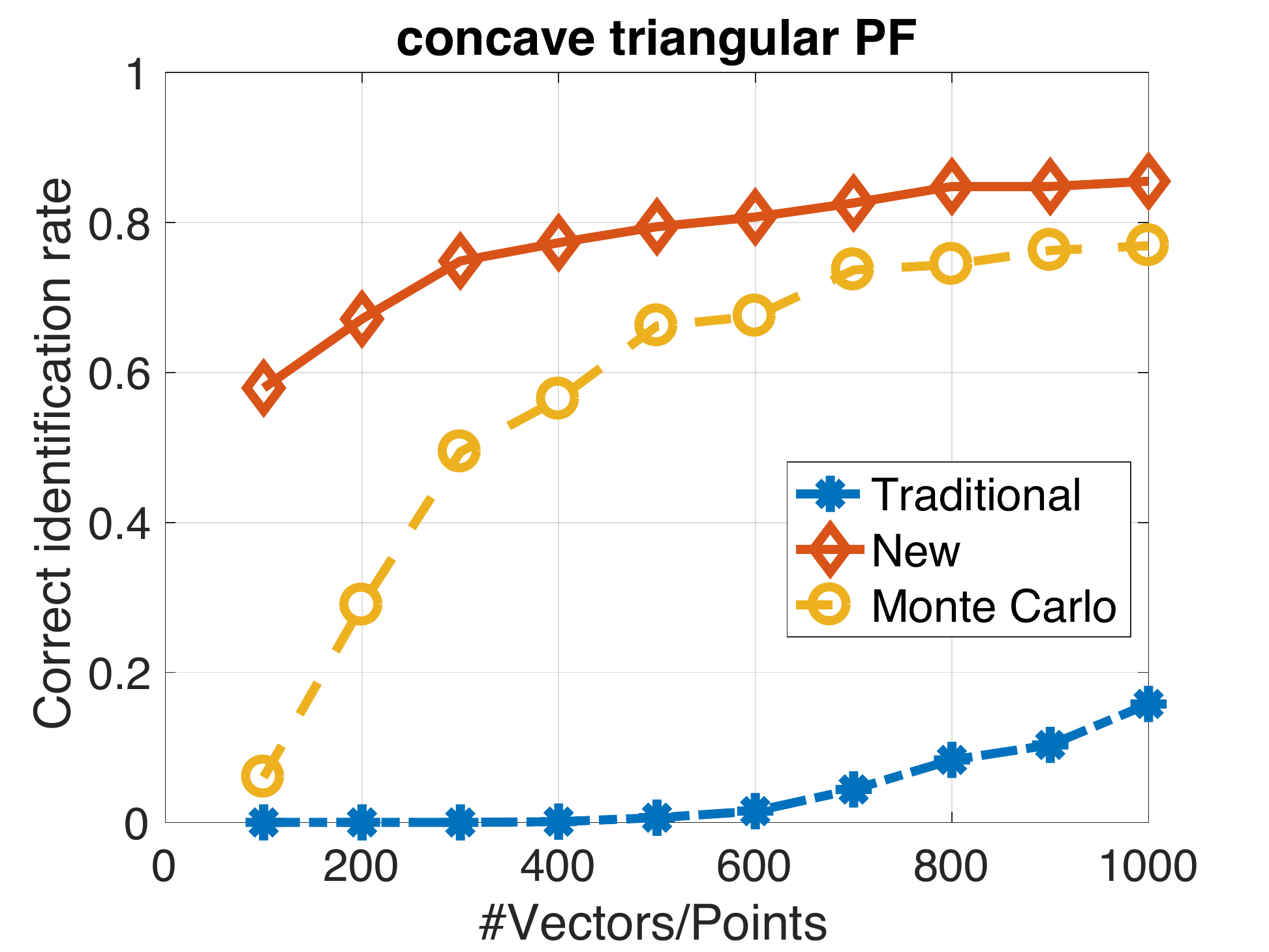}                %以pic.jpg的0.5倍大小输出
\end{minipage}}
\caption{The approximation accuracy with respect to the number of the direction vectors and the sampling points on the triangular PF solution sets.} %                         %大图名称
\label{numvec1}                                                        %图片引用标记
\end{figure}

\begin{figure}[!htbp]
\centering                                                          %居中
\subfigure{                    %第一张子图
\hspace{-0.5cm}
\begin{minipage}{0.45\columnwidth}\centering                                                          %子图居中
\includegraphics[scale=0.23]{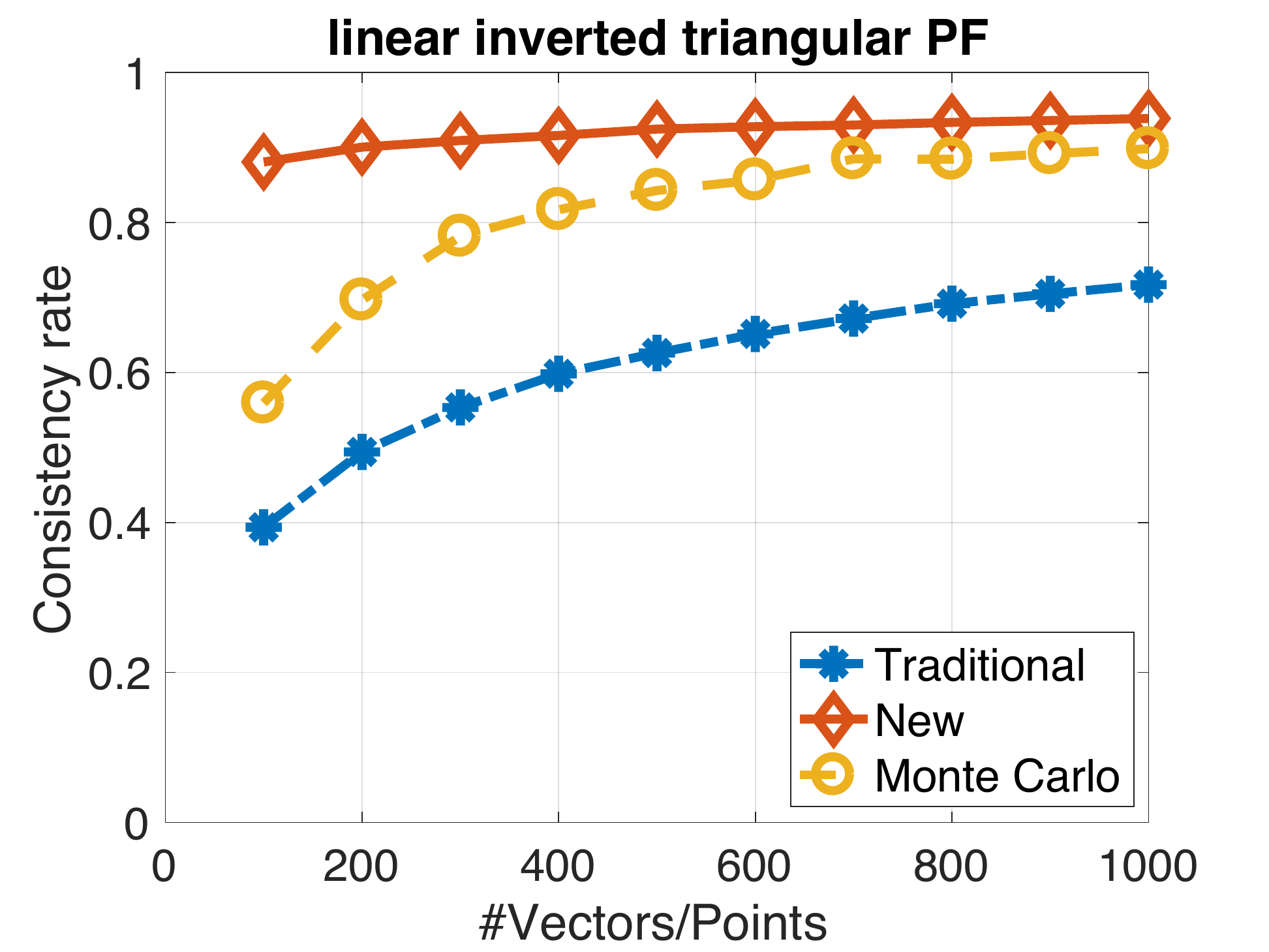}               %以pic.jpg的0.5倍大小输出
\end{minipage}}
\subfigure{                    %第二张子图
\begin{minipage}{0.45\columnwidth}\centering                                                          %子图居中
\includegraphics[scale=0.23]{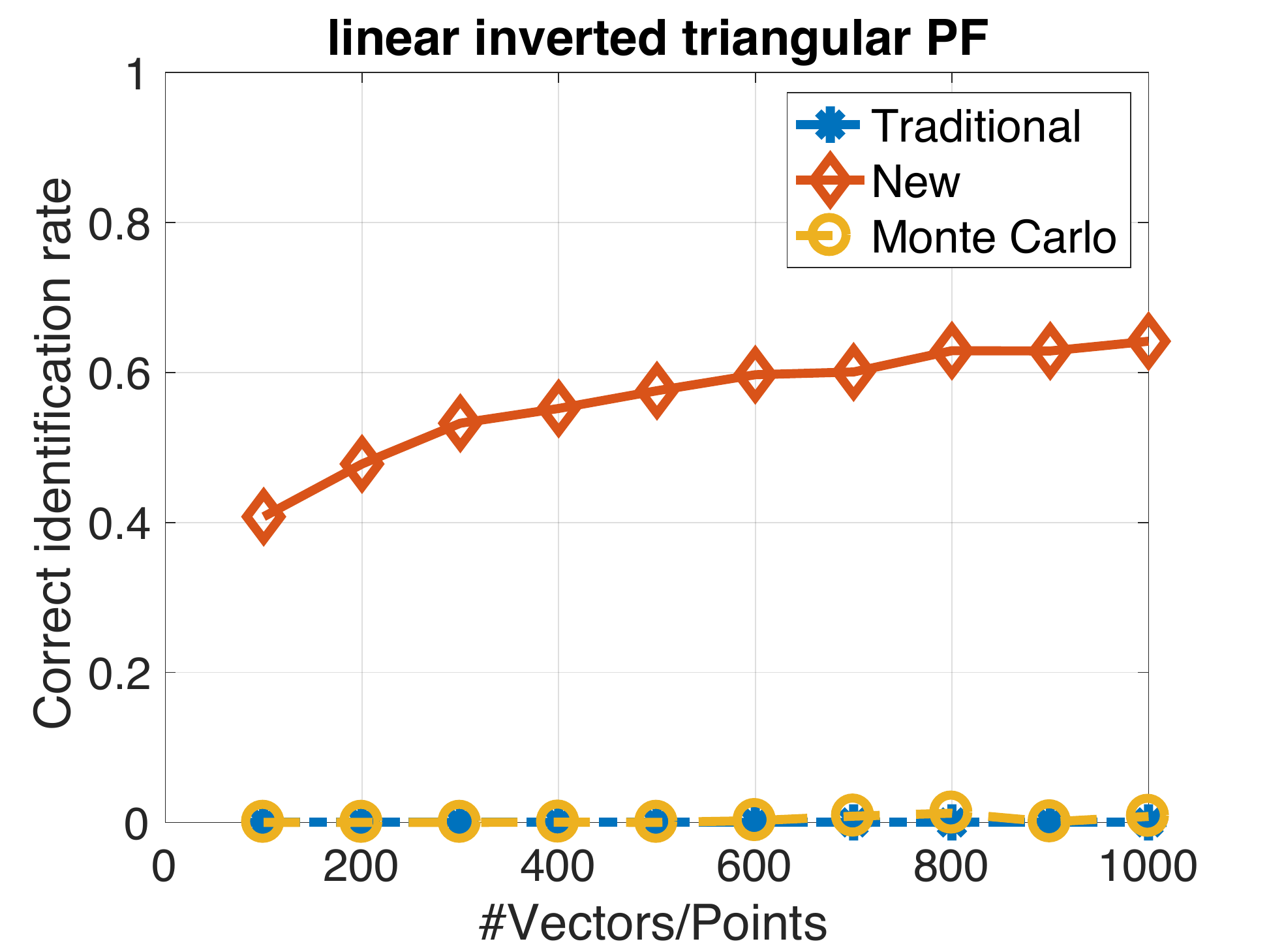}                %以pic.jpg的0.5倍大小输出
\end{minipage}}
\subfigure{      
\hspace{-0.5cm}              %第一张子图
\begin{minipage}{0.45\columnwidth}\centering                                                          %子图居中
\includegraphics[scale=0.23]{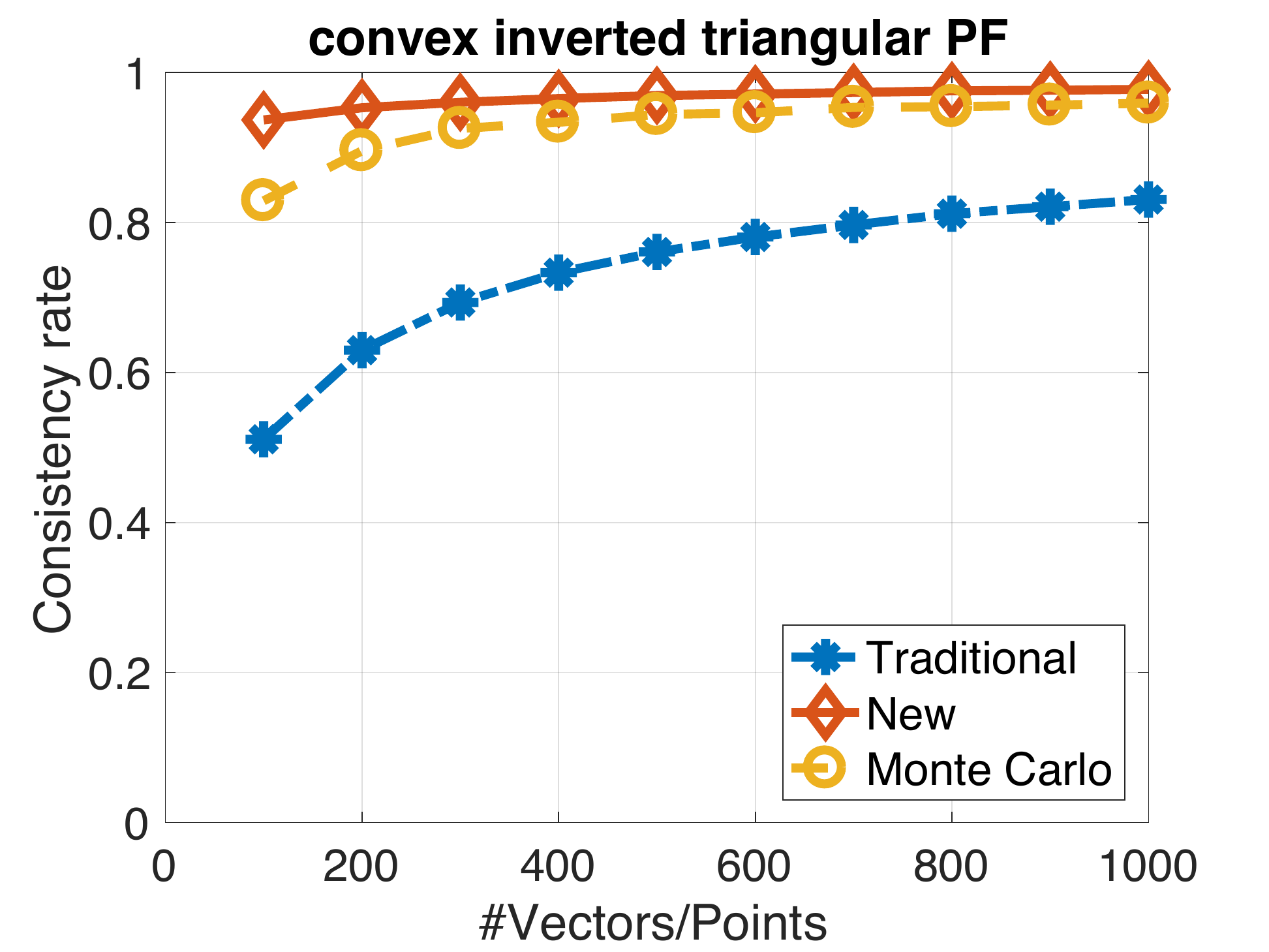}               %以pic.jpg的0.5倍大小输出
\end{minipage}}
\subfigure{                    %第二张子图
\begin{minipage}{0.45\columnwidth}\centering                                                          %子图居中
\includegraphics[scale=0.23]{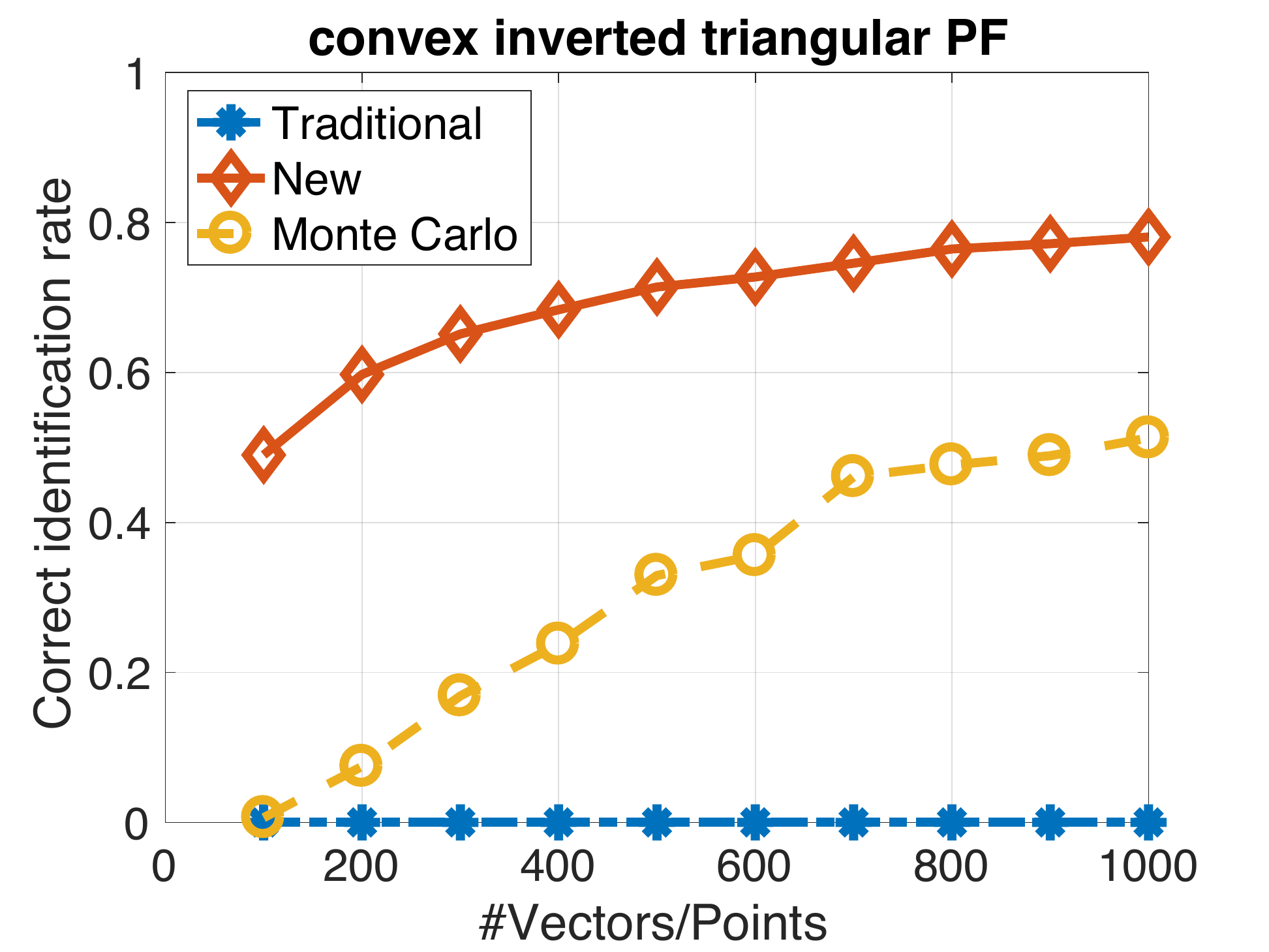}                %以pic.jpg的0.5倍大小输出
\end{minipage}}
\subfigure{      
\hspace{-0.5cm}              %第一张子图
\begin{minipage}{0.45\columnwidth}\centering                                                          %子图居中
\includegraphics[scale=0.23]{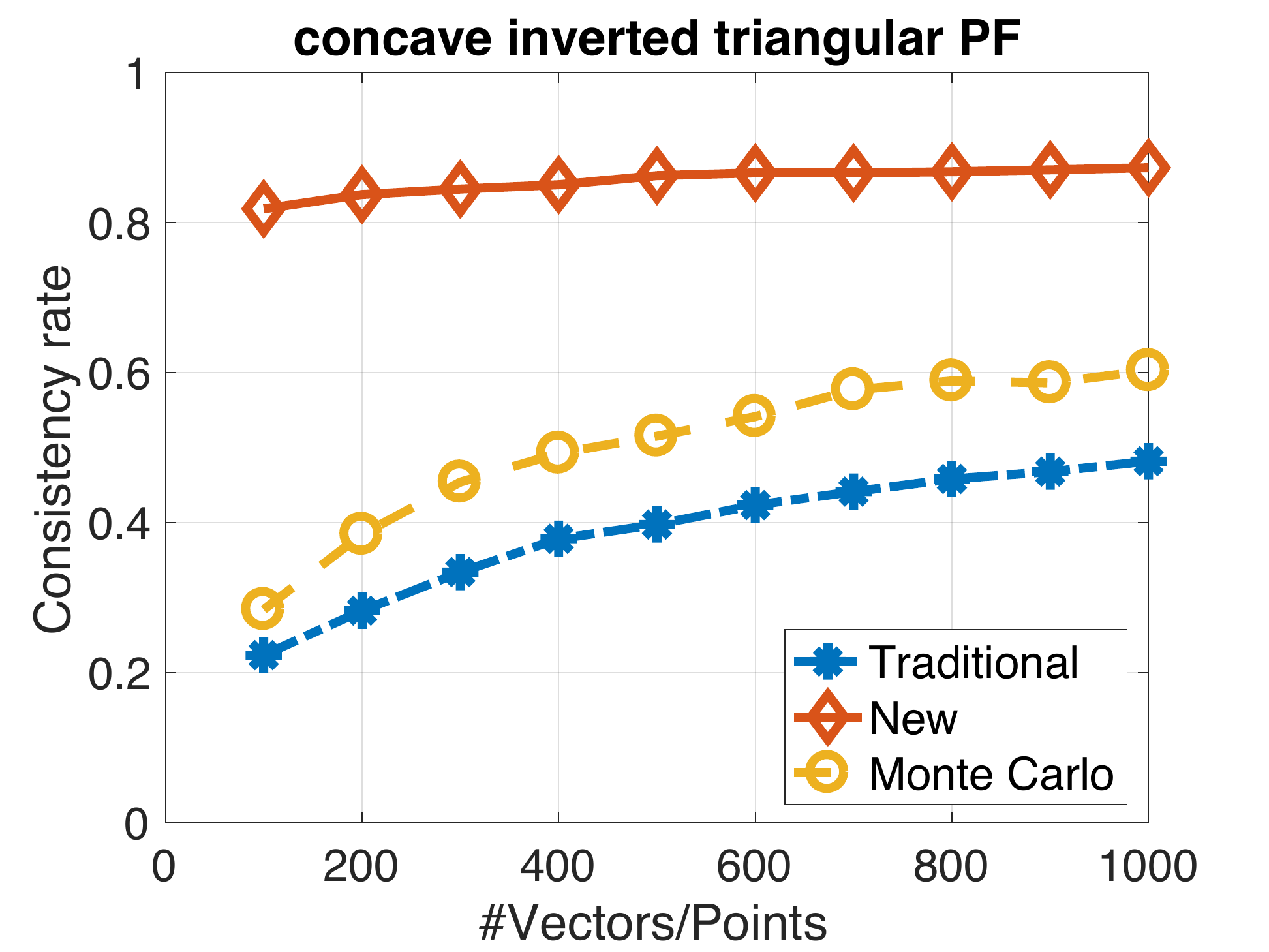}               %以pic.jpg的0.5倍大小输出
\end{minipage}}
\subfigure{                    %第二张子图
\begin{minipage}{0.45\columnwidth}\centering                                                          %子图居中
\includegraphics[scale=0.23]{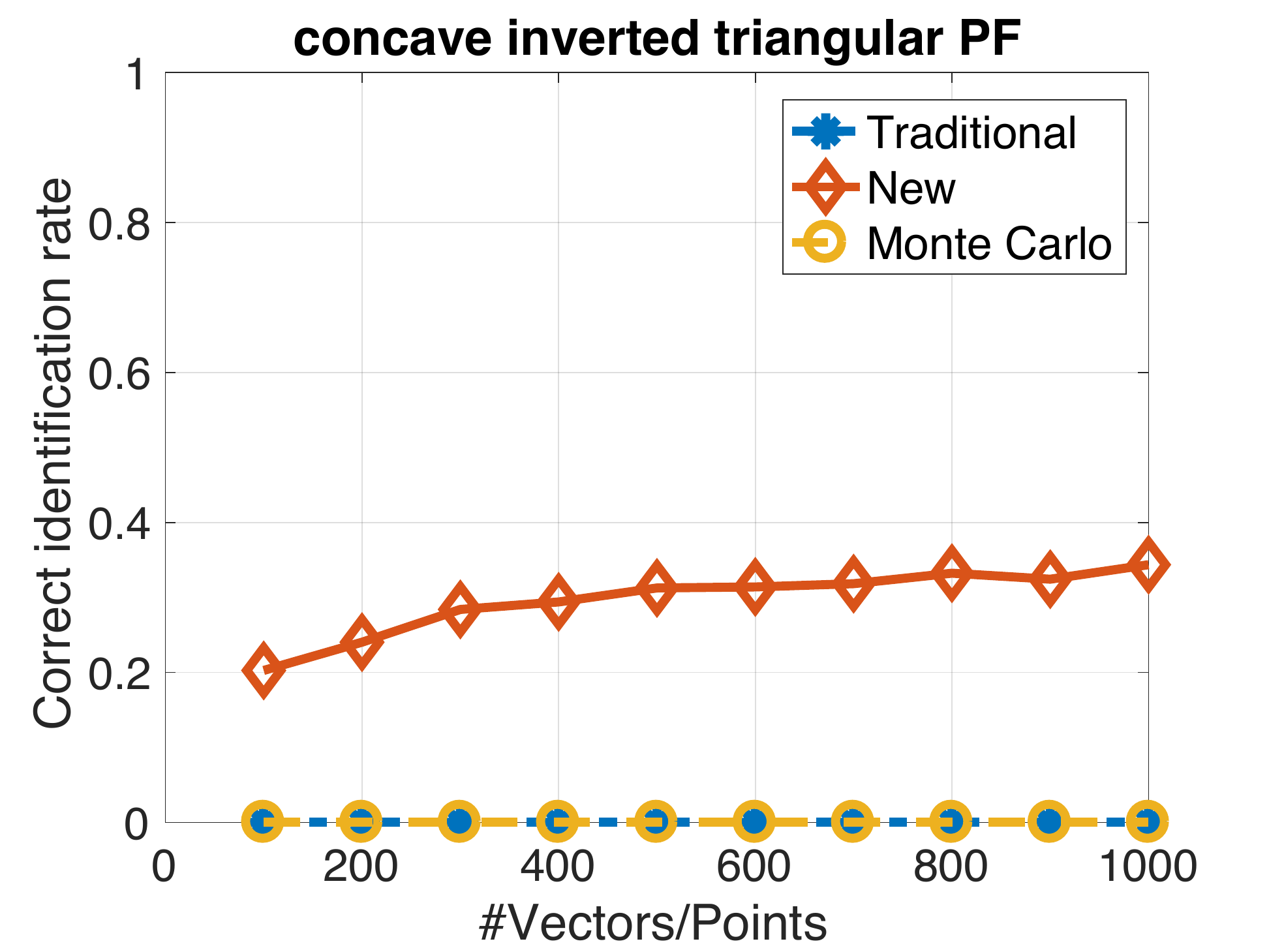}                %以pic.jpg的0.5倍大小输出
\end{minipage}}
\caption{The approximation accuracy with respect to the number of the direction vectors and the sampling points on the inverted triangular PF solution sets.} %                         %大图名称
\label{numvec2}                                                        %图片引用标记
\end{figure}

\textcolor{black}{As studied in \cite{Ishibuchi:2017}, appropriate reference point specification is essential in the hypervolume-based EMOAs for fair performance comparison. The Monte Carlo method is not practical because it only performs well on the triangular PF solution sets with $r=0.0$ (i.e., when the reference point is the nadir point of the PF), and generally this is not a good choice for the reference point specification. Usually a point worse than (dominated by) the nadir point is set as the reference point. In the following experiments, we will fix the reference point $r=-0.2$ because this is a reasonable reference point specification according to \cite{Ishibuchi:2017}.}
\textcolor{black}{
\subsubsection{The effect of the number of direction vectors and sampling points}
\label{section4c2}
Next we examine the effect of the number of the direction vectors and the sampling points on the performance of the three approximation methods. We fix the number of solutions in each solution set $N=100$. The results on 5-dimension solution sets are shown in Fig.~\ref{numvec1}-\ref{numvec2}. The transverse axis represents the number of the direction vectors for the traditional and the new methods, or the number of the sampling points for the Monte Carlo sampling method.}

\textcolor{black}{
From Fig.~\ref{numvec1}-\ref{numvec2} we can clearly see that the new method outperforms the traditional method and the Monte Carlo sampling method. With the number of the direction vectors and the sampling points increases, the performance of the three approximation methods can be improved. However, this is not so obvious to the Monte Carlo method and the traditional method on the convex triangular, linear inverted triangular and concave inverted triangular PF solution sets because the correct identification rate does not increase as the number of the direction vectors and the sampling points increases in these cases.}

\begin{figure}[!htbp]
\centering                                                          %居中
\subfigure{        
\hspace{-0.5cm}            %第一张子图
\begin{minipage}{0.45\columnwidth}\centering                                                          %子图居中
\includegraphics[scale=0.23]{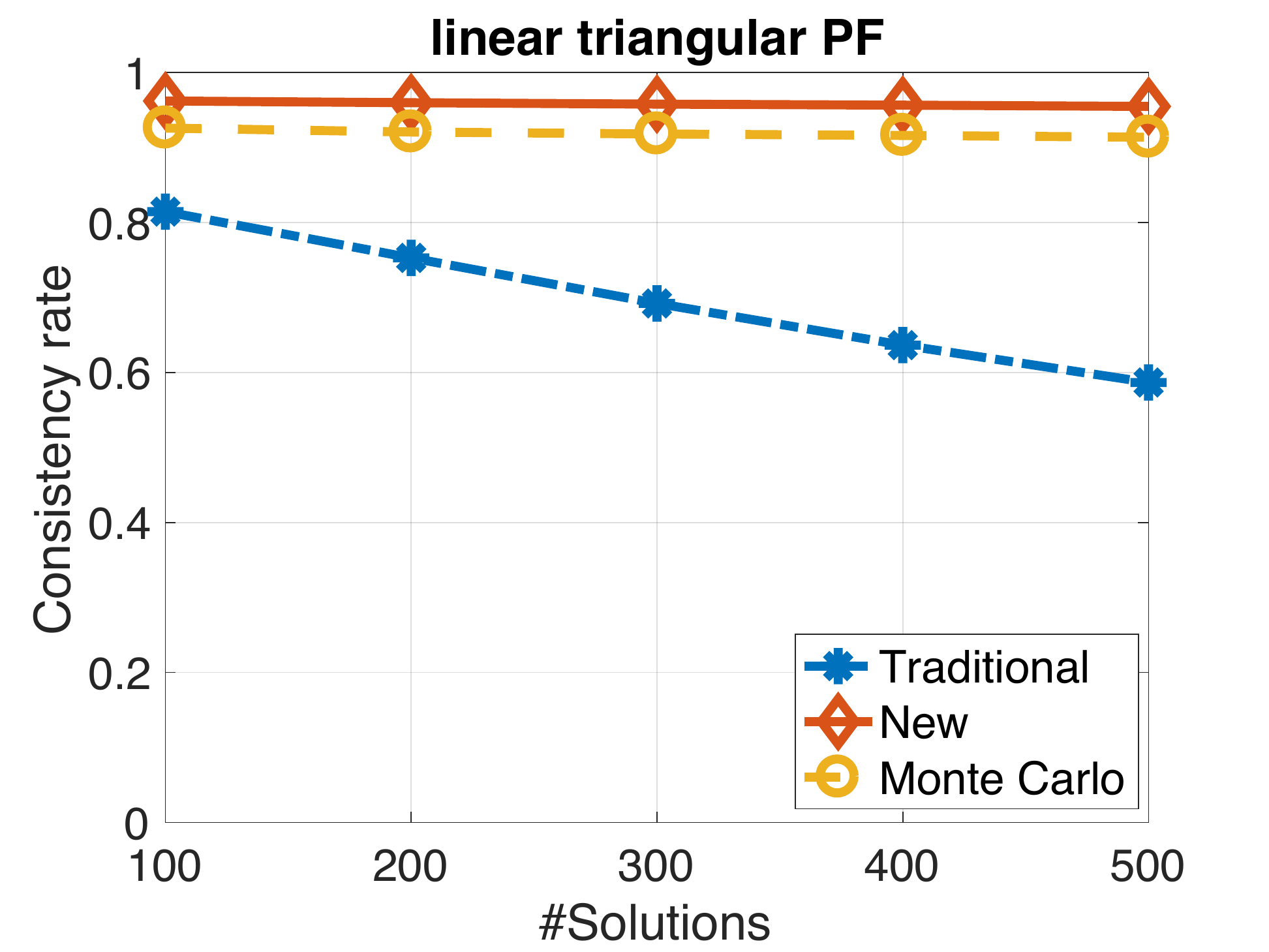}               %以pic.jpg的0.5倍大小输出
\end{minipage}}
\subfigure{                    %第二张子图
\begin{minipage}{0.45\columnwidth}\centering                                                          %子图居中
\includegraphics[scale=0.23]{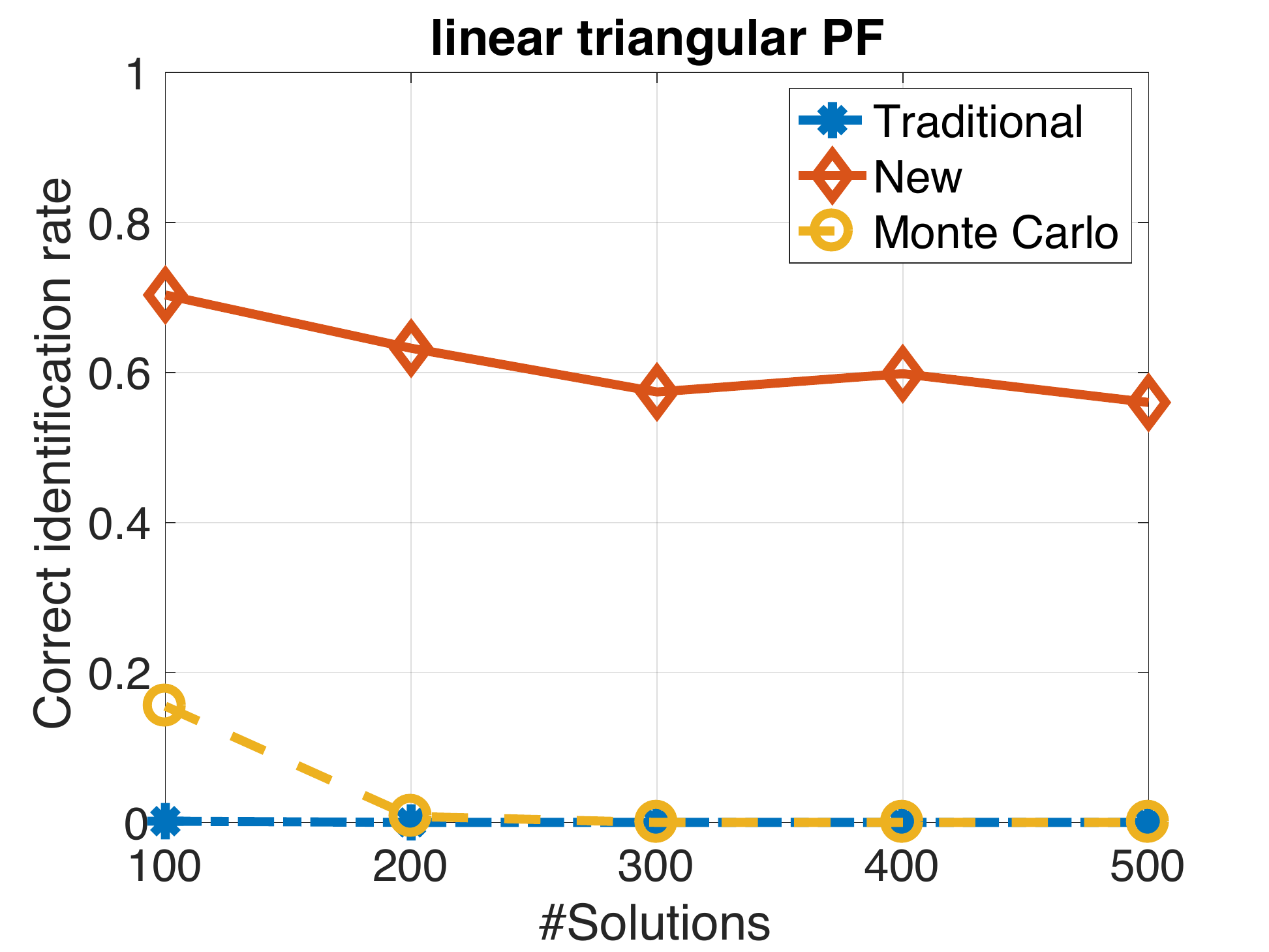}                %以pic.jpg的0.5倍大小输出
\end{minipage}}
\subfigure{        
\hspace{-0.5cm}            %第一张子图
\begin{minipage}{0.45\columnwidth}\centering                                                          %子图居中
\includegraphics[scale=0.23]{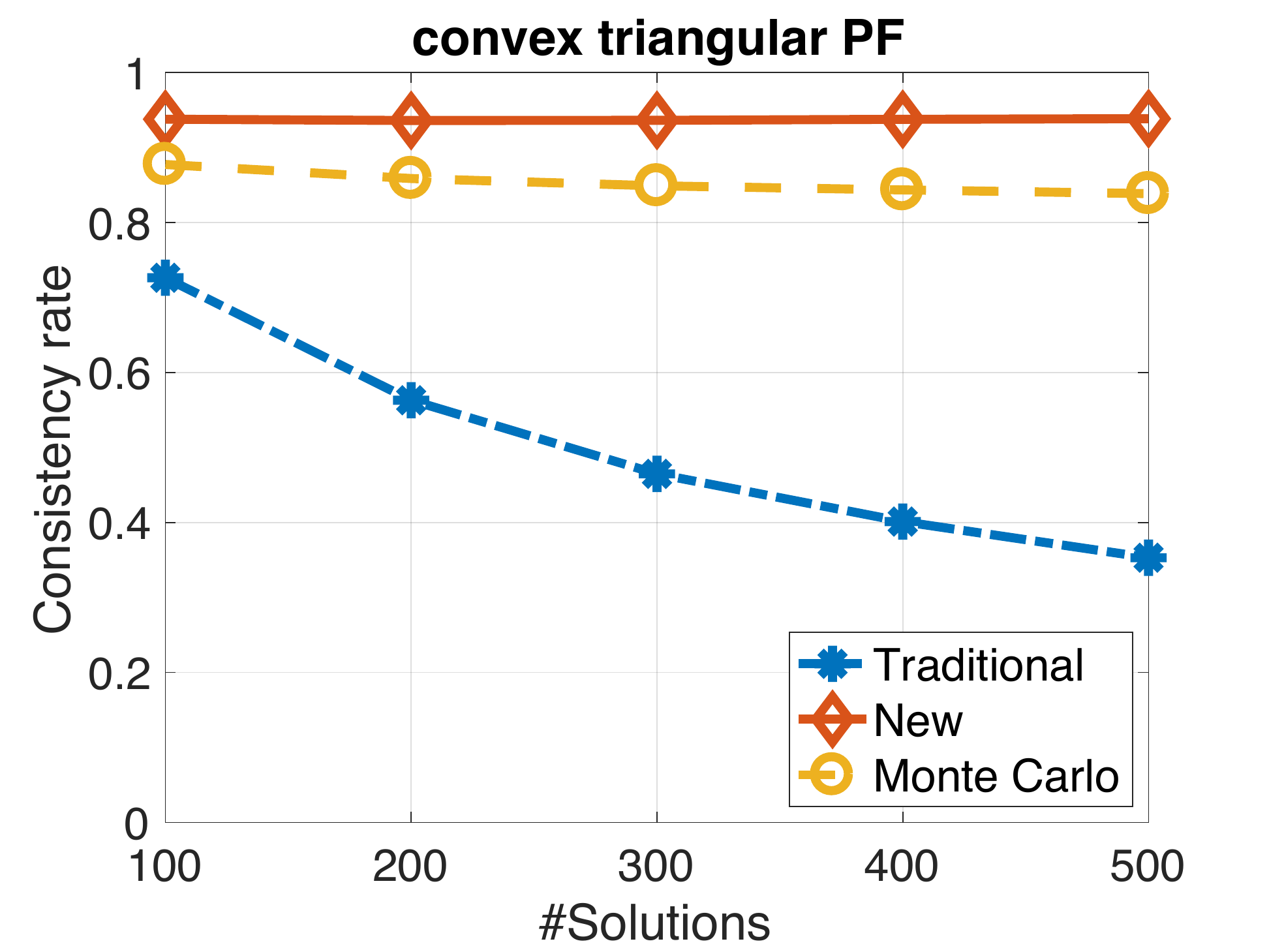}               %以pic.jpg的0.5倍大小输出
\end{minipage}}
\subfigure{                    %第二张子图
\begin{minipage}{0.45\columnwidth}\centering                                                          %子图居中
\includegraphics[scale=0.23]{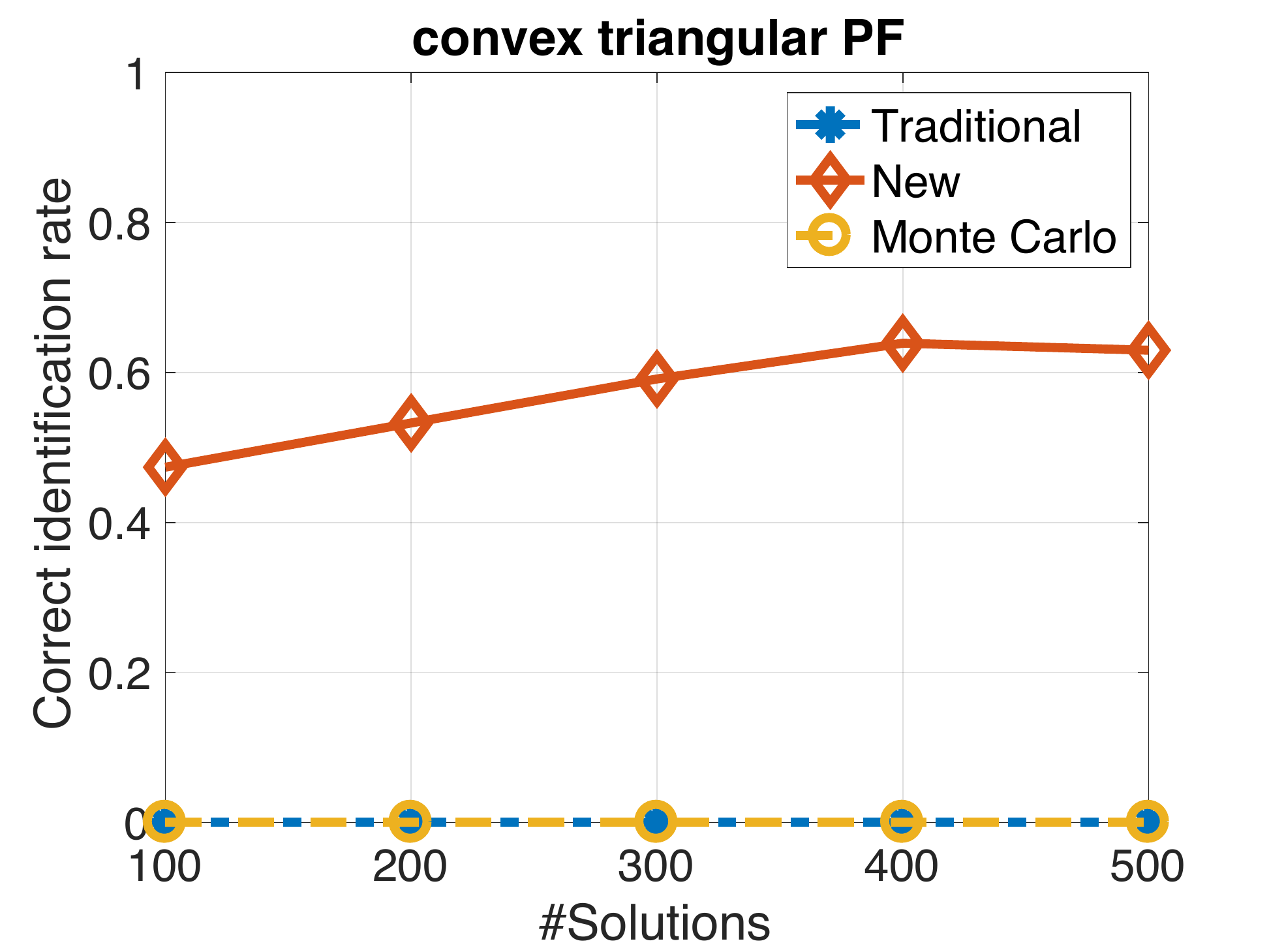}                %以pic.jpg的0.5倍大小输出
\end{minipage}}
\subfigure{ 
\hspace{-0.5cm}                   %第一张子图
\begin{minipage}{0.45\columnwidth}\centering                                                          %子图居中
\includegraphics[scale=0.23]{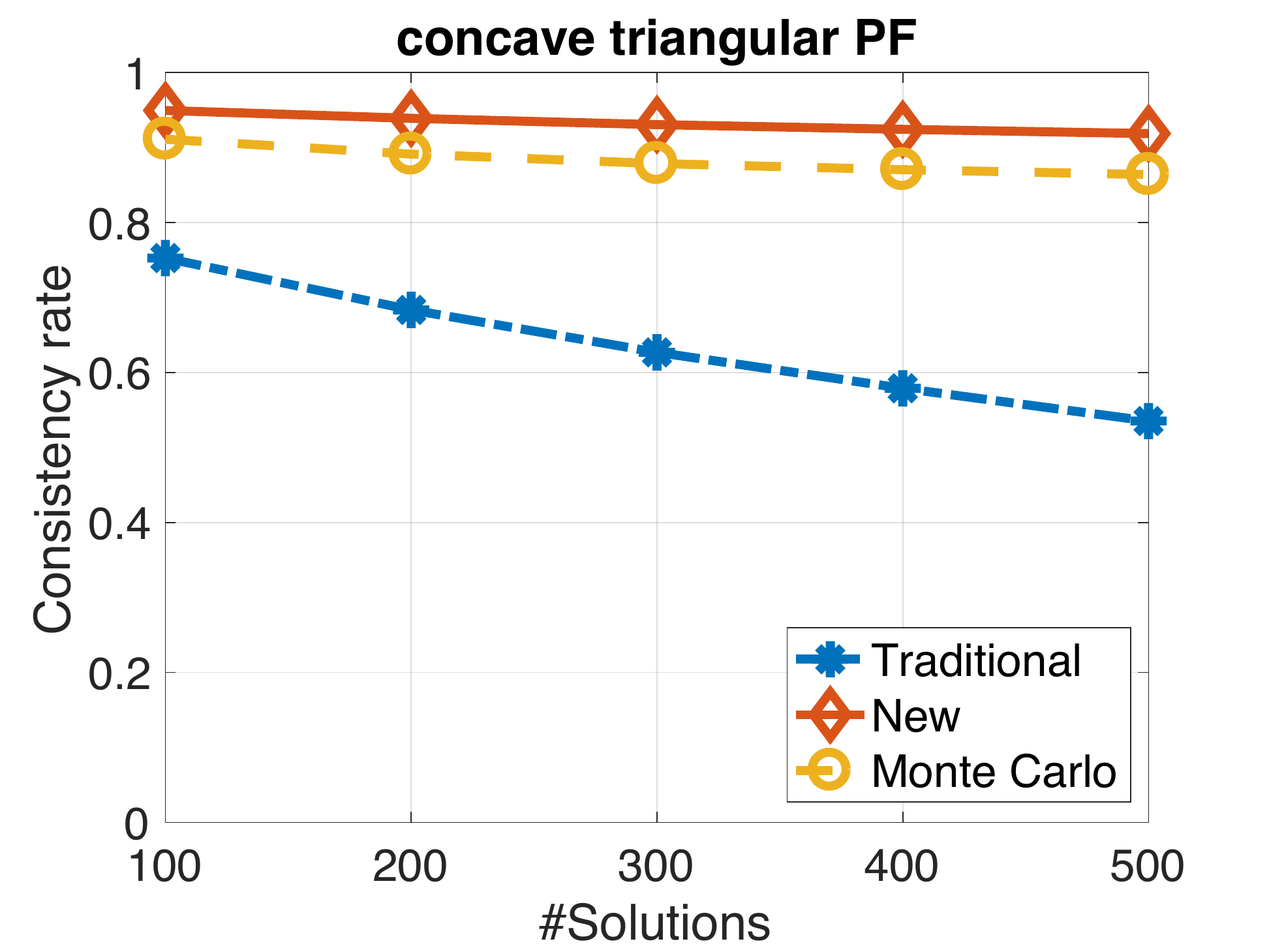}               %以pic.jpg的0.5倍大小输出
\end{minipage}}
\subfigure{                    %第二张子图
\begin{minipage}{0.45\columnwidth}\centering                                                          %子图居中
\includegraphics[scale=0.23]{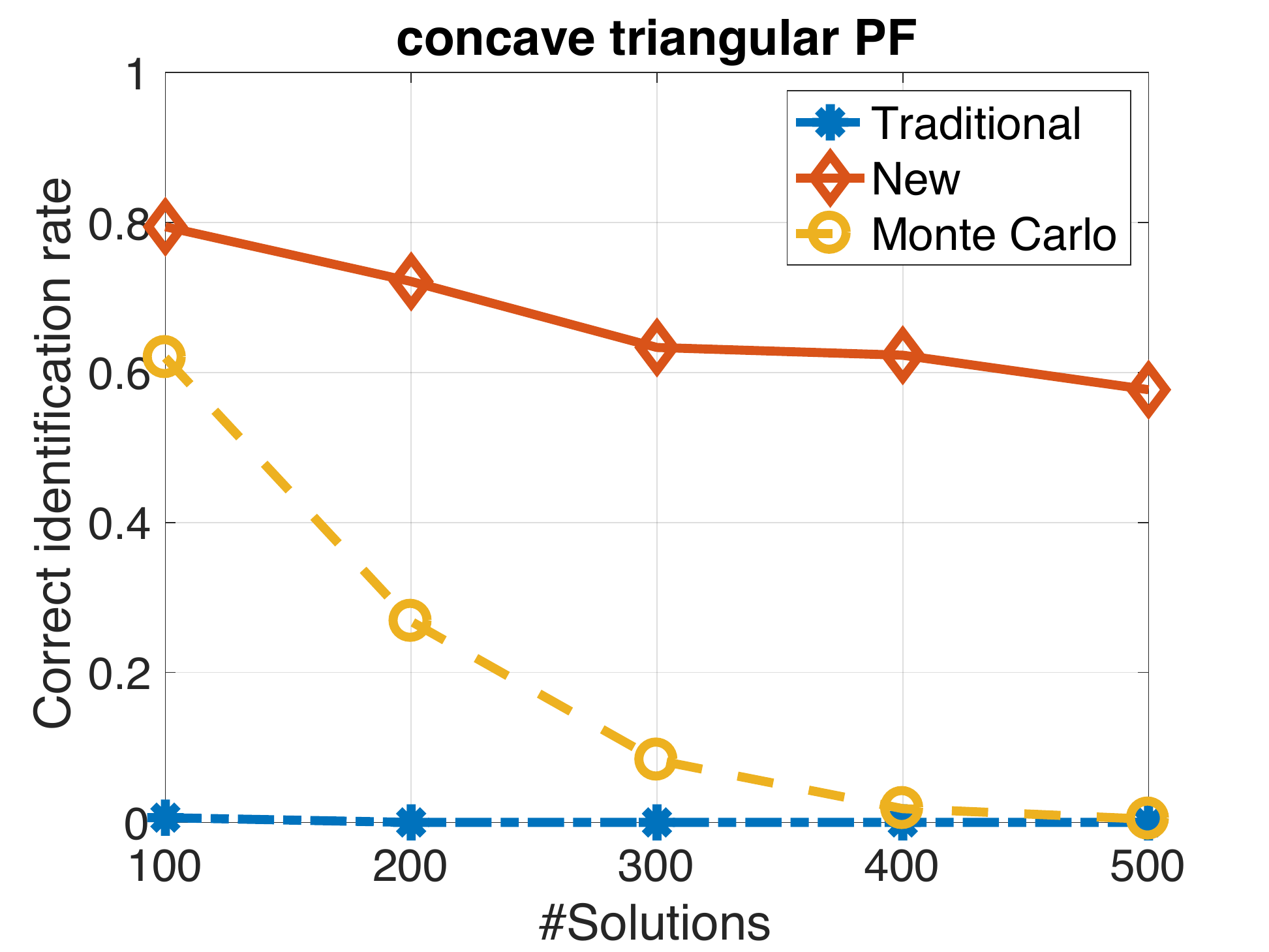}                %以pic.jpg的0.5倍大小输出
\end{minipage}}
\caption{The approximation accuracy with respect to the number of solutions on the triangular PF solution sets.} %                         %大图名称
\label{numsol1}                                                        %图片引用标记
\end{figure}

\begin{figure}[!htbp]
\centering                                                          %居中
\subfigure{                    %第一张子图
\hspace{-0.5cm}
\begin{minipage}{0.45\columnwidth}\centering                                                          %子图居中
\includegraphics[scale=0.23]{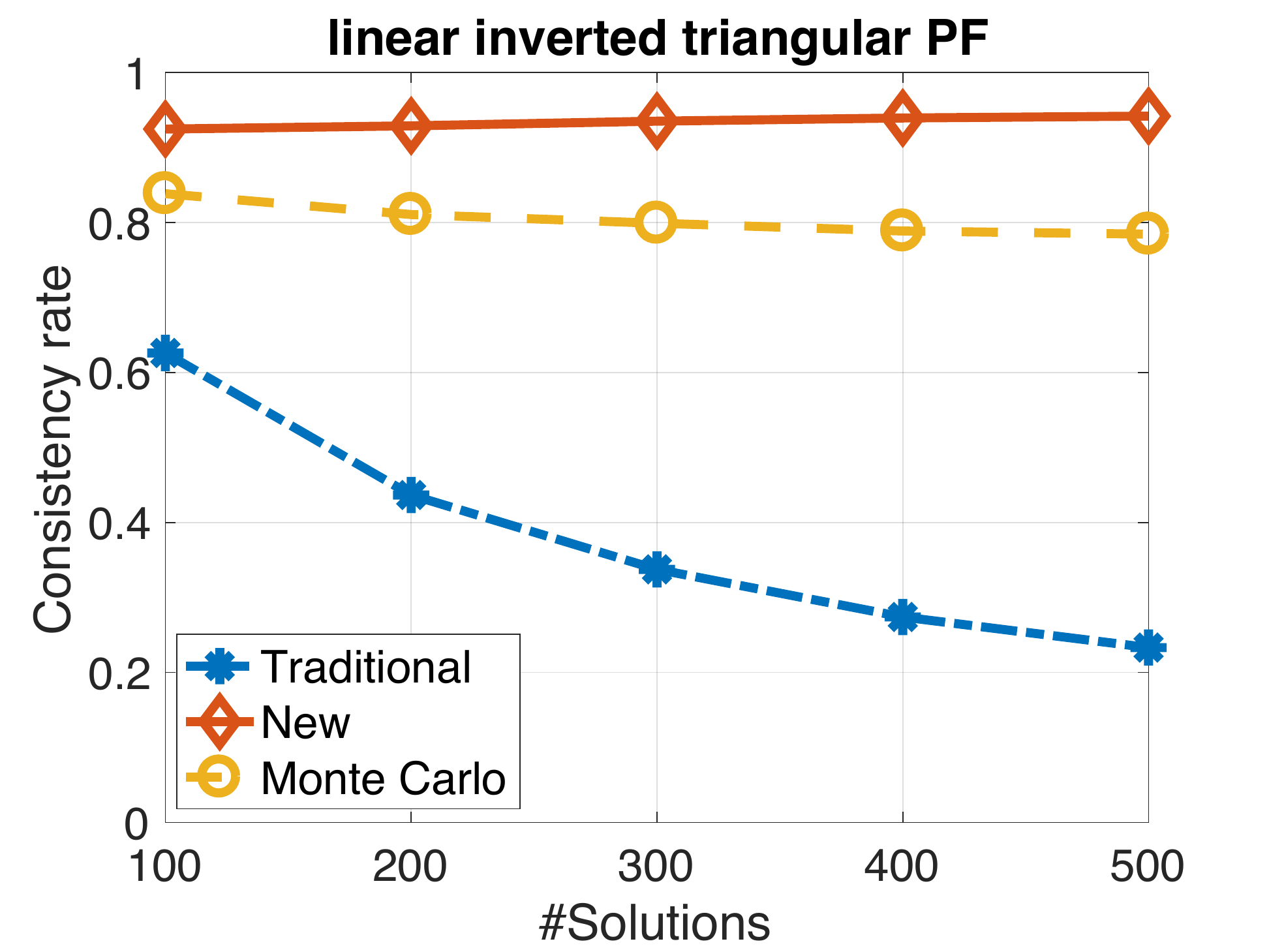}               %以pic.jpg的0.5倍大小输出
\end{minipage}}
\subfigure{                    %第二张子图
\begin{minipage}{0.45\columnwidth}\centering                                                          %子图居中
\includegraphics[scale=0.23]{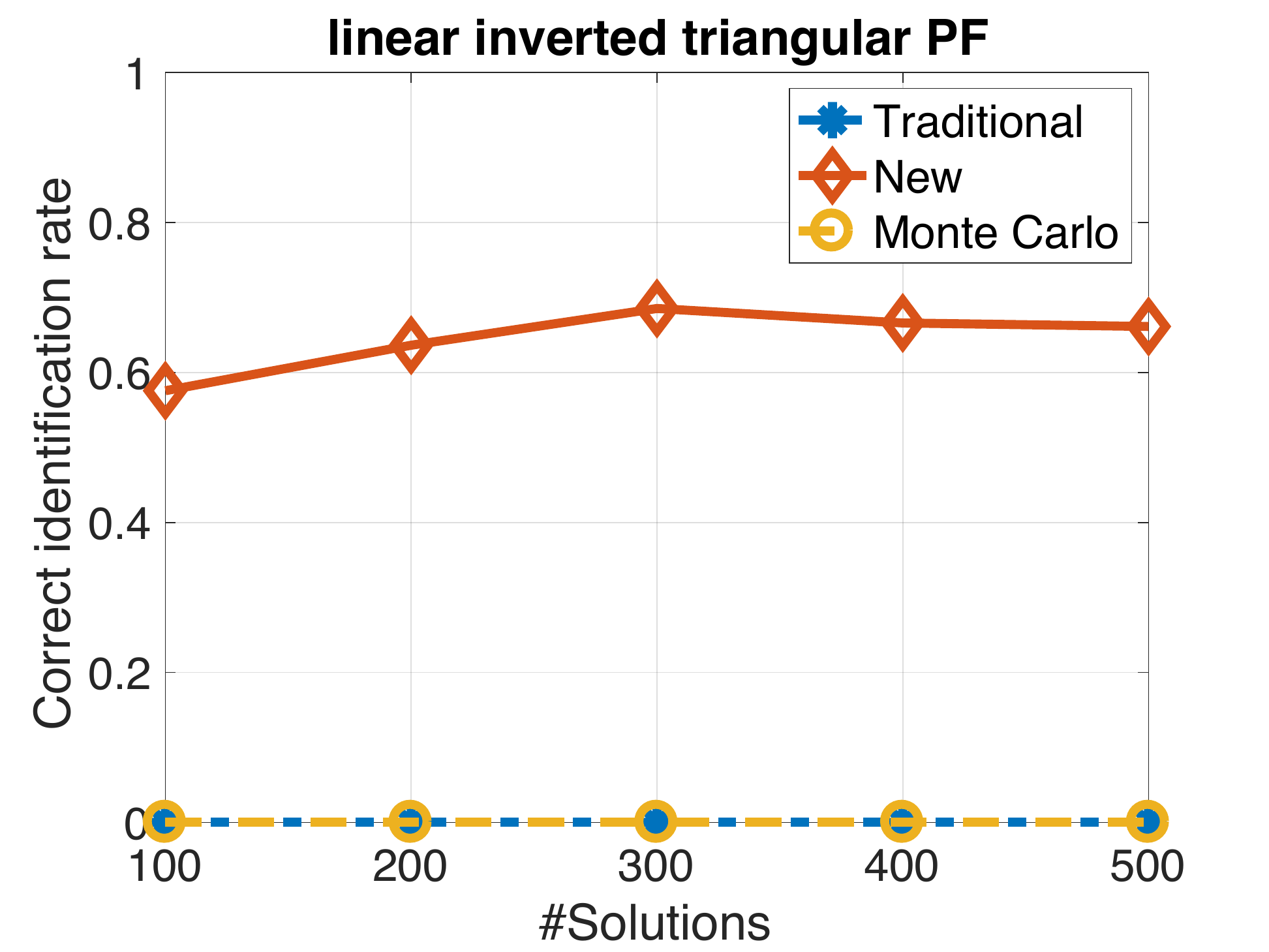}                %以pic.jpg的0.5倍大小输出
\end{minipage}}
\subfigure{      
\hspace{-0.5cm}              %第一张子图
\begin{minipage}{0.45\columnwidth}\centering                                                          %子图居中
\includegraphics[scale=0.23]{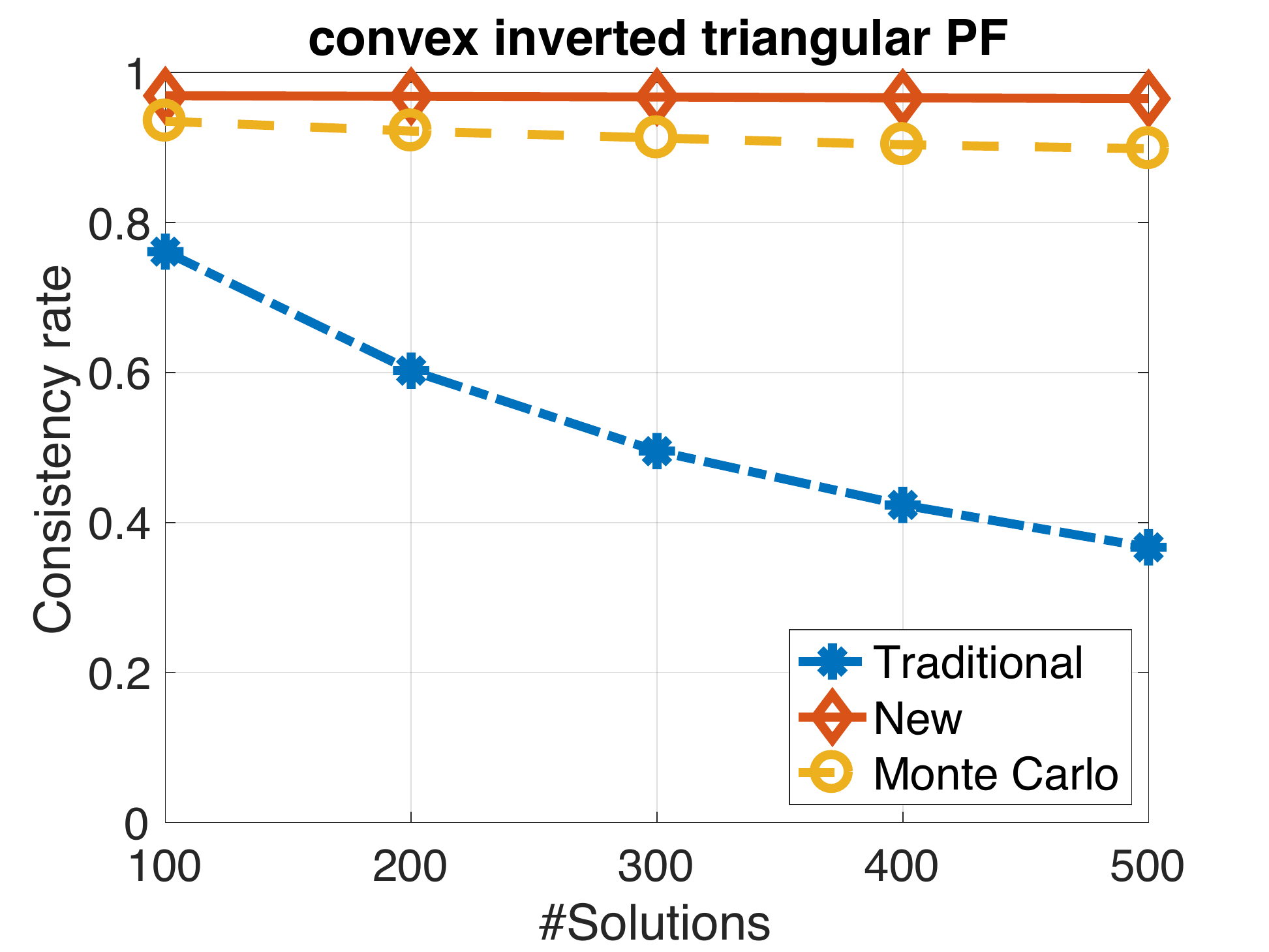}               %以pic.jpg的0.5倍大小输出
\end{minipage}}
\subfigure{                    %第二张子图
\begin{minipage}{0.45\columnwidth}\centering                                                          %子图居中
\includegraphics[scale=0.23]{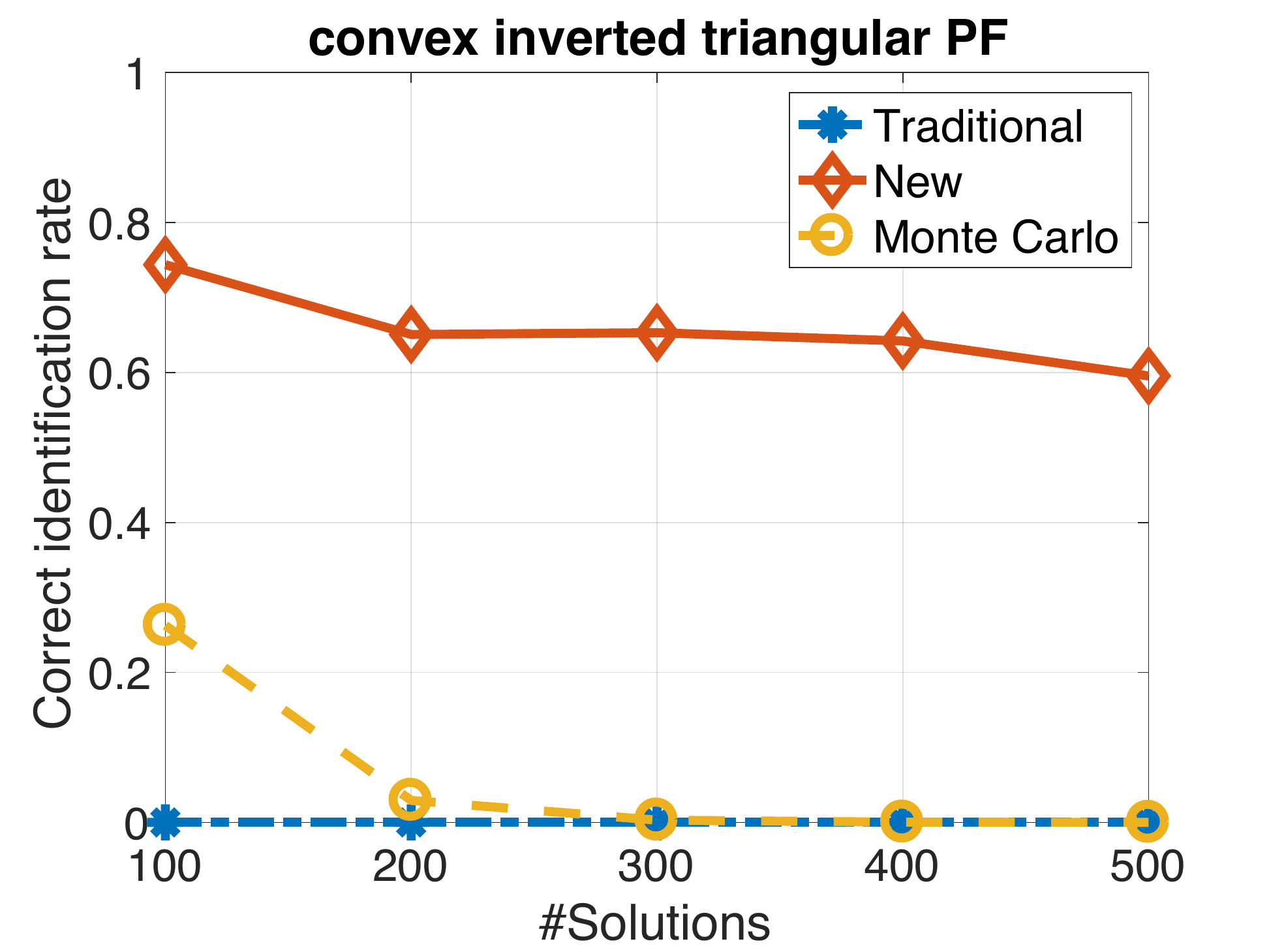}                %以pic.jpg的0.5倍大小输出
\end{minipage}}
\subfigure{      
\hspace{-0.5cm}              %第一张子图
\begin{minipage}{0.45\columnwidth}\centering                                                          %子图居中
\includegraphics[scale=0.23]{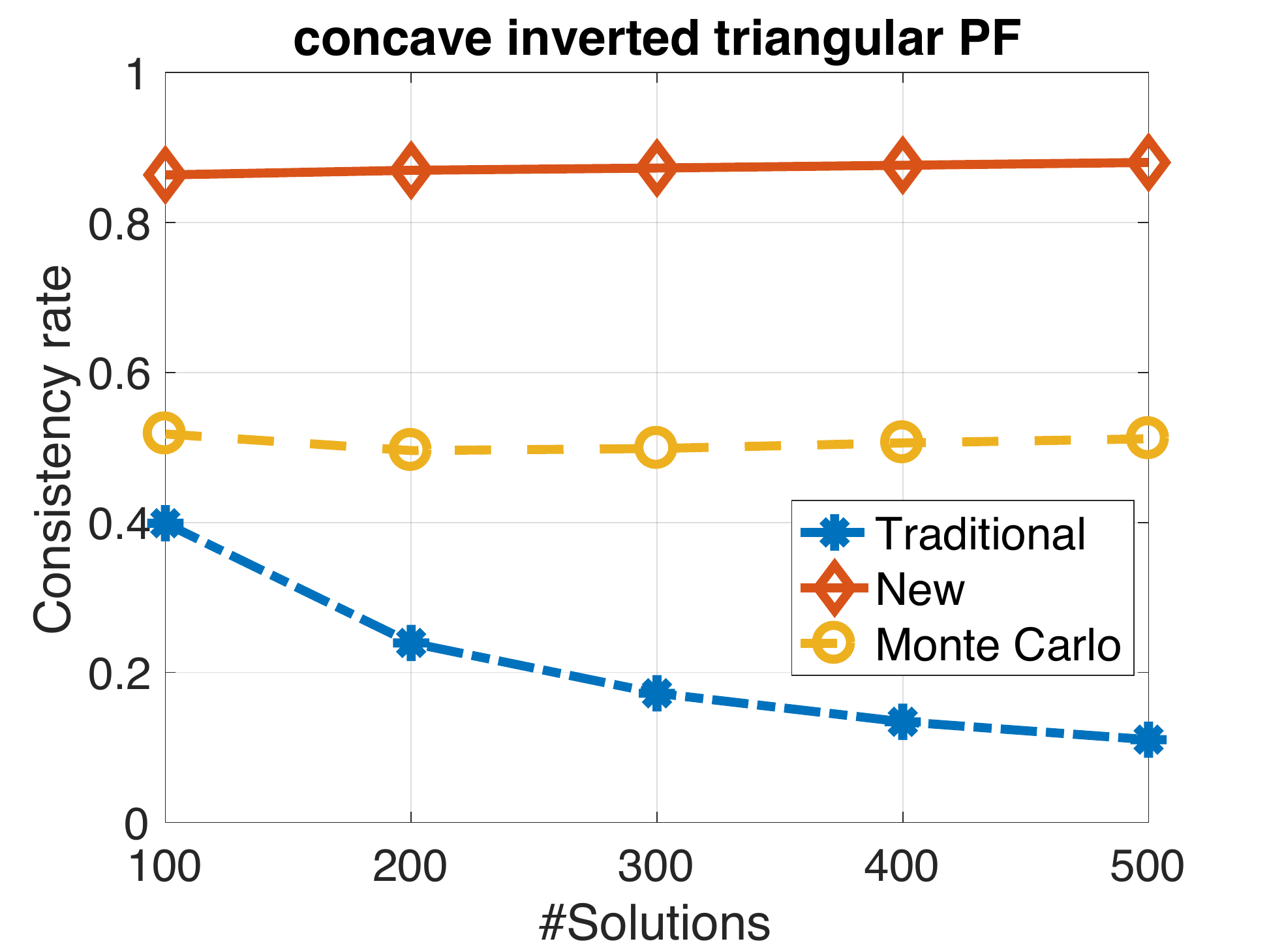}               %以pic.jpg的0.5倍大小输出
\end{minipage}}
\subfigure{                    %第二张子图
\begin{minipage}{0.45\columnwidth}\centering                                                          %子图居中
\includegraphics[scale=0.23]{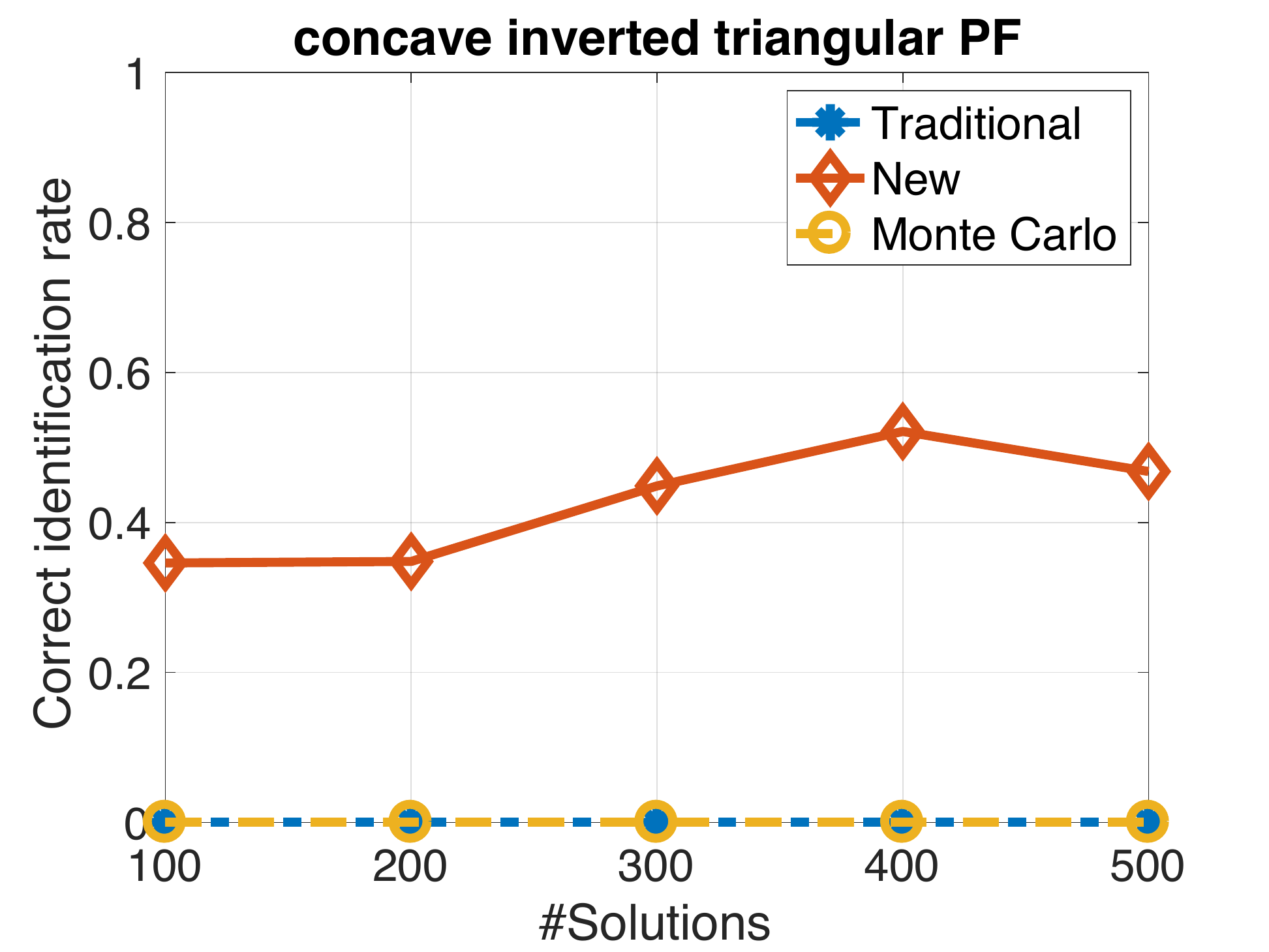}                %以pic.jpg的0.5倍大小输出
\end{minipage}}
\caption{The approximation accuracy with respect to the number of solutions on the inverted triangular PF solution sets.} %                         %大图名称
\label{numsol2}                                                        %图片引用标记
\end{figure}

\textcolor{black}{
Another interesting observation is that the new method is able to achieve a good performance even when the number of the direction vectors is small (e.g., 100). For almost all cases, the new method with 100 direction vectors achieves a comparable or better performance than the other two methods with 1000 direction vectors and sampling points. This clearly shows the advantage of the new method over the other two methods.}
\textcolor{black}{
\subsubsection{The effect of the number of solutions}
Now we examine the effect of the number of solutions in each solution set on the performance of the three approximation methods. We fix the number of the direction vectors and the sampling points to 500. The results on 5-dimension solution sets are shown in Fig.~\ref{numsol1}-\ref{numsol2}. The transverse axis represents the number of solutions in each solution set.}

\textcolor{black}{
We can see that the new method clearly outperforms the other two methods. The traditional method performs the worst as its correct identification rate is always 0 in all cases. In general, with the increase of the number of solutions, the performance of the Monte Carlo method and the traditional method is non-increasing. However, this is not always the case for the new method. As we can observe that for some cases (e.g., convex triangular, linear inverted triangular and concave inverted triangular PFs), its performance can be improved with the number of solutions increases. Moreover, the new method shows a robust performance with respect to the number of solutions. }

\begin{figure}[!htbp]
\centering                                                          %居中
\subfigure{                    %第一张子图
\hspace{-0.5cm}
\begin{minipage}{0.45\columnwidth}\centering                                                          %子图居中
\includegraphics[scale=0.23]{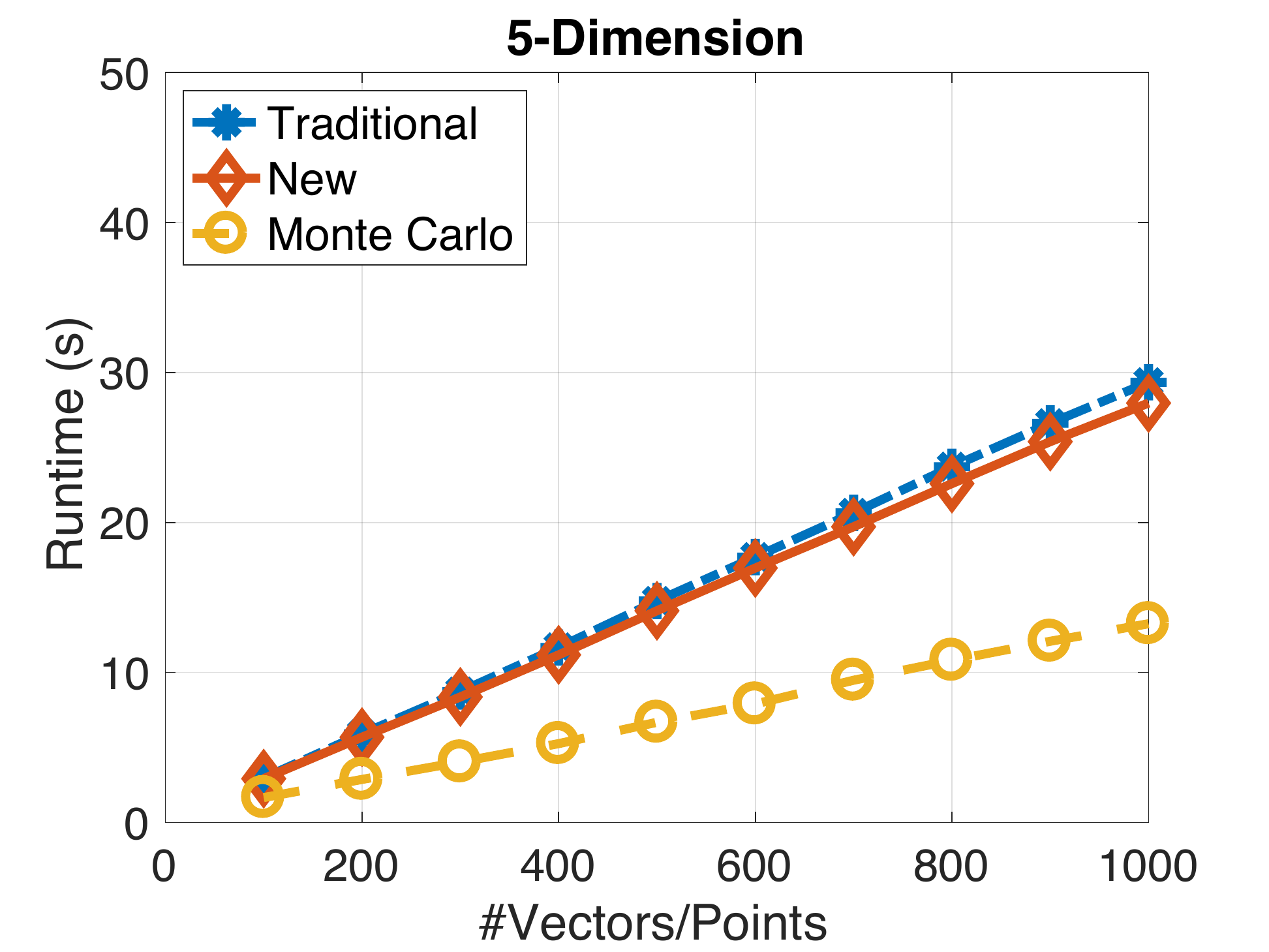}               %以pic.jpg的0.5倍大小输出
\end{minipage}}
\subfigure{                    %第二张子图
\begin{minipage}{0.45\columnwidth}\centering                                                          %子图居中
\includegraphics[scale=0.23]{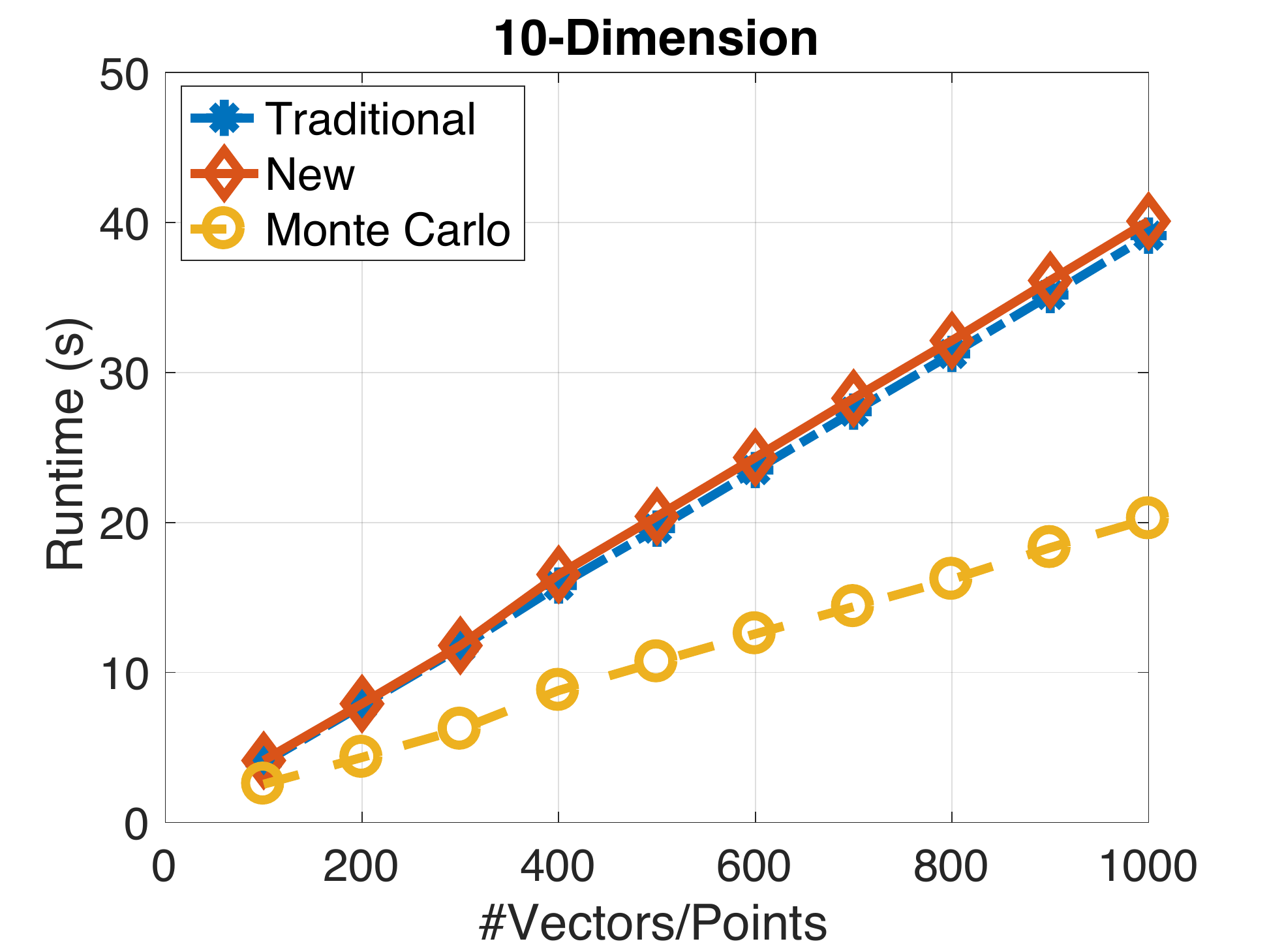}                %以pic.jpg的0.5倍大小输出
\end{minipage}}
\subfigure{                    %第二张子图
\begin{minipage}{0.45\columnwidth}\centering                                                          %子图居中
\includegraphics[scale=0.23]{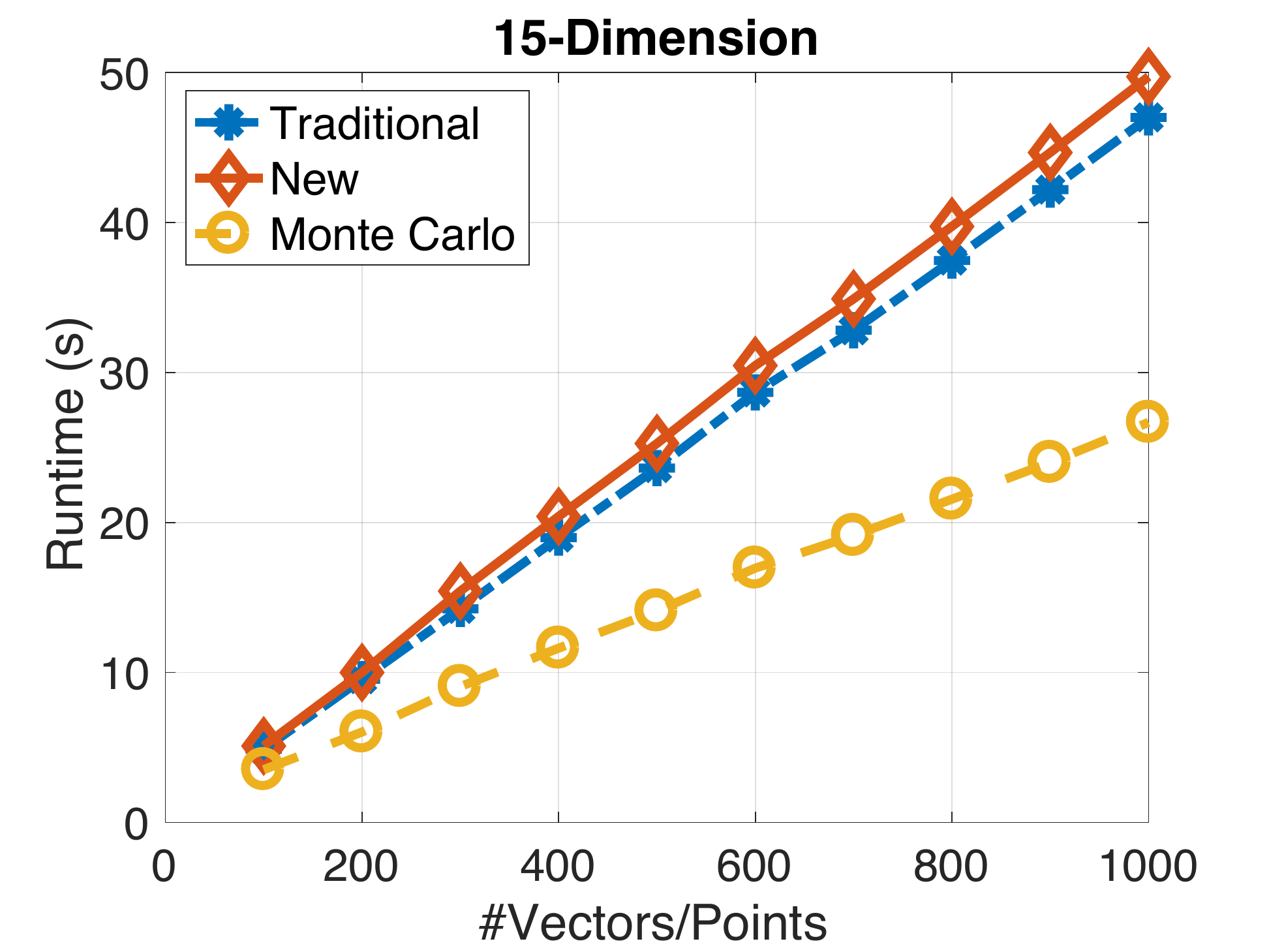}                %以pic.jpg的0.5倍大小输出
\end{minipage}}
\caption{Runtime of the three approximation methods on the linear triangular PF solution sets.} %                         %大图名称
\label{runtime}                                                        %图片引用标记
\end{figure}

As similar observations are obtained from our computational experiments on the 10-dimension solution sets, we do not show their results in this letter. Their results are provided in the supplementary material. 
\textcolor{black}{
\subsection{Runtime comparison}
Now we compare the runtime of the three approximation methods and the two exact methods to evaluate their computational efficiency. We fix the number of solutions $N=100$. The linear triangular PF solution sets are chosen for illustration. Fig.~\ref{runtime} shows the runtime results of the three approximation methods on 5-, 10- and 15-dimension\footnote{To illustrate the computational efficiency of the approximation methods in high-dimensional spaces, we use the method of generating 5- and 10-dimension solution sets to generate 15-dimension solution sets.} solution sets. The runtime of each method is the total time consumed by each method to approximate the hypervolume contributions of all solutions in the 100 solution sets (i.e., $100\times 100$ solutions in total).}

\textcolor{black}{
From the results we can see that the runtime of the new method and the traditional method is comparable to each other. The Monte Carlo method takes about 50\% of the runtime of the new method when the same number of the direction vectors and the sampling points are used. For the three approximation methods, the runtime increases linearly as the number of the direction vectors and the sampling points increases. We can also observe that the runtime increases linearly as the dimension increases for a fixed number of the direction vectors and the sampling points.}

\begin{figure}[!htbp]
\centering                                                          %居中
\subfigure{                    %第一张子图
\hspace{-0.5cm}
\begin{minipage}{0.45\columnwidth}\centering                                                          %子图居中
\includegraphics[scale=0.23]{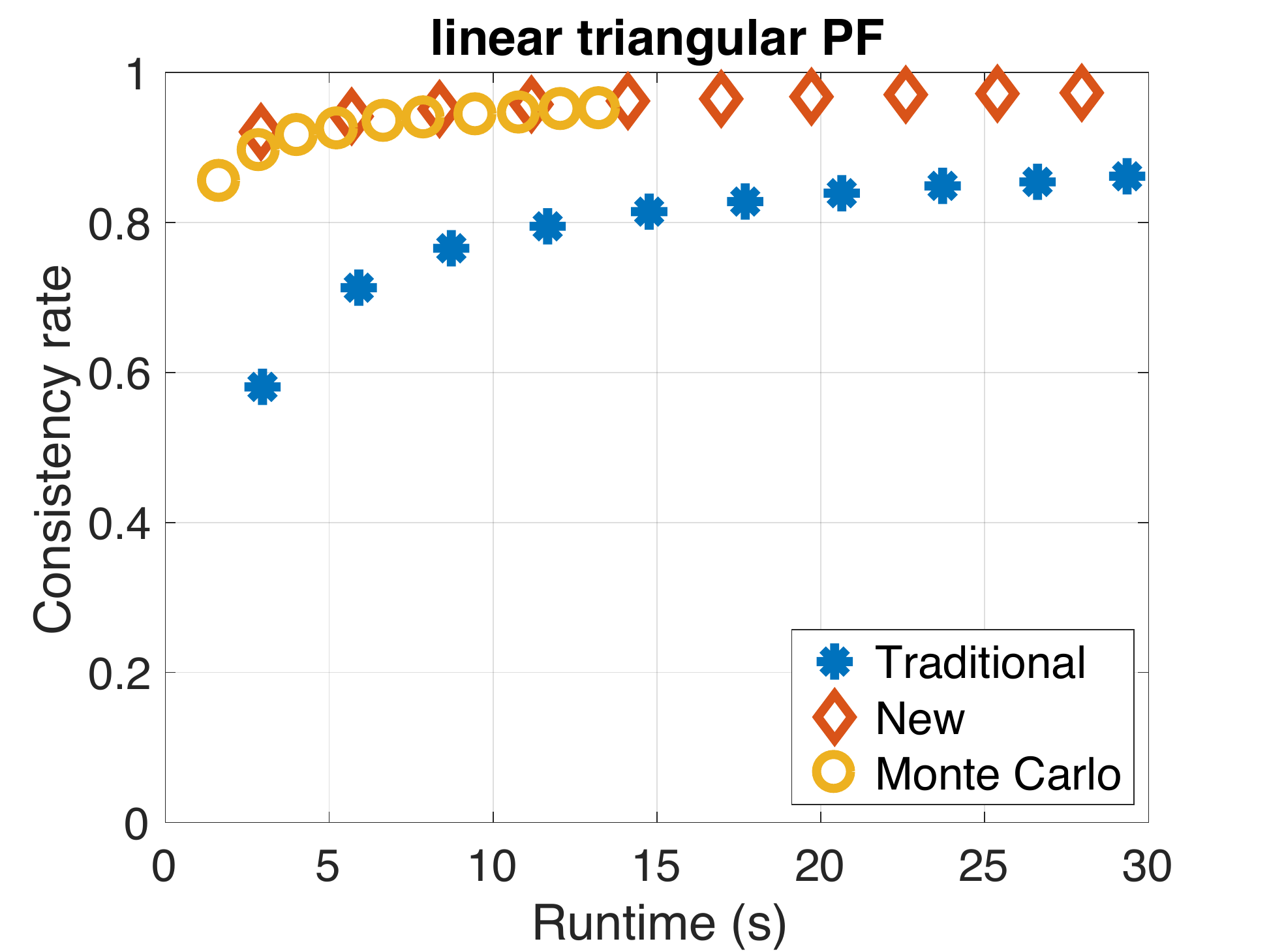}               %以pic.jpg的0.5倍大小输出
\end{minipage}}
\subfigure{                    %第二张子图
\begin{minipage}{0.45\columnwidth}\centering                                                          %子图居中
\includegraphics[scale=0.23]{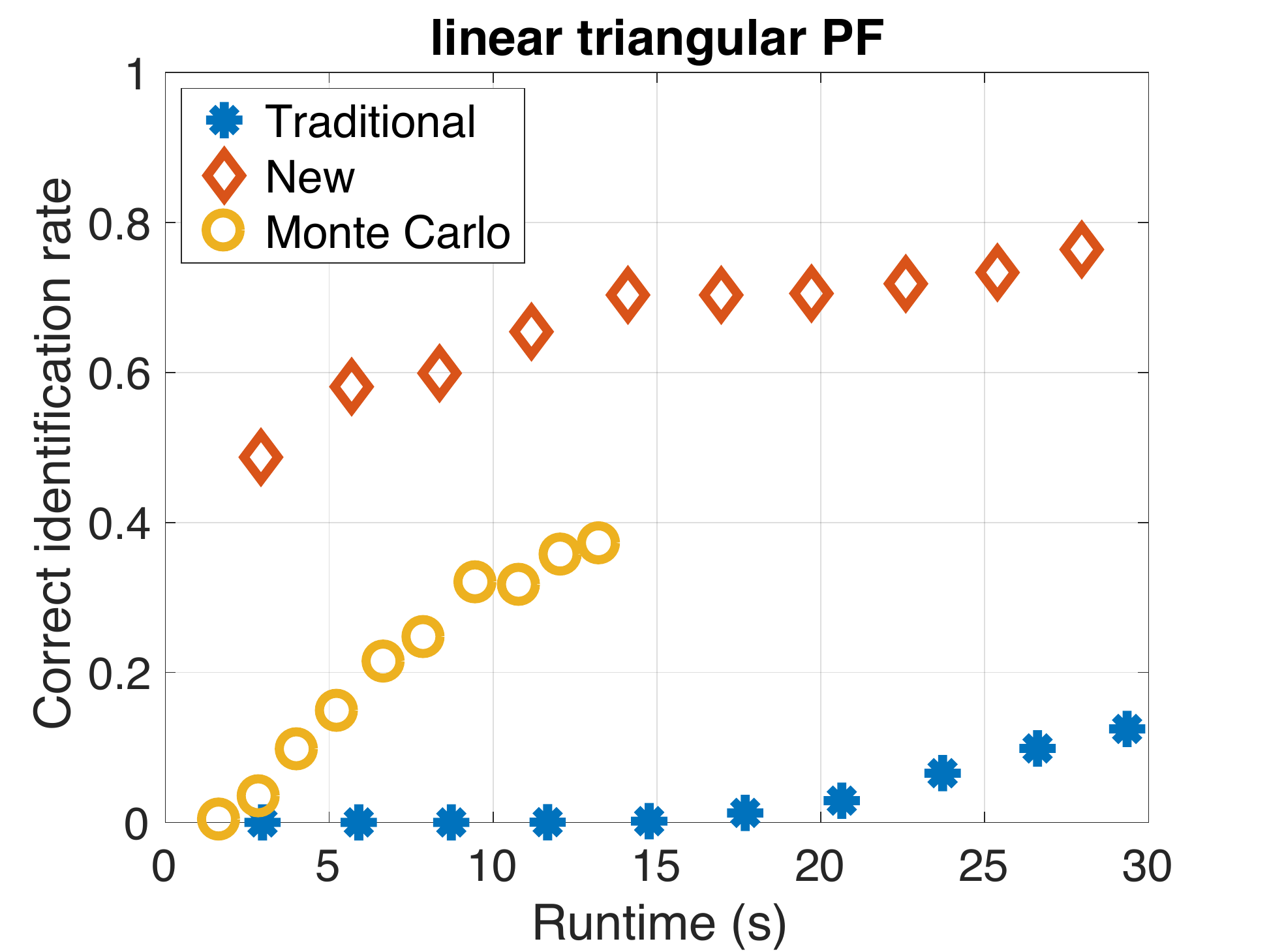}                %以pic.jpg的0.5倍大小输出
\end{minipage}}
\caption{The relation between the approximation accuracy and the runtime of the three approximation methods on the 5-dimension linear triangular PF solution sets.} %                         %大图名称
\label{accuracytime}                                                        %图片引用标记
\end{figure}

\textcolor{black}{
In Section \ref{section4c2}, we compared the approximation accuracy of the three approximation methods, and the new method outperforms the other two methods with the same number of the direction vectors and the sampling points. Here, the Monte Carlo method outperforms the other two methods with respect to the runtime with the same number of the direction vectors and the sampling points. In practice, an approximation method with higher approximation accuracy and faster runtime is preferred. In order to evaluate the three approximation methods from this point of view, we plot the relation between the approximation accuracy and the runtime of the three methods on the 5-dimension linear triangular PF solution sets as shown in Fig.~\ref{accuracytime}.}

\textcolor{black}{
From the results we can clearly observe that the new method is able to run faster and at the same time achieve higher approximation accuracy than the other two methods. The Monte Carlo method is able to run faster than the new method, but it can not achieve better approximation accuracy at the same time. In this sense, the new method clearly outperforms the other two methods.}

\textcolor{black}{
Finally, we compare the runtime of the three approximation methods with the two exact methods. The results on the linear triangular PF solution sets are shown in Table \ref{compareexact}. From the results we can see that the two exact methods are computationally more efficient than the three approximation methods with 100 and 1000 direction vectors and sampling points in 5-dimension case. However, the exact methods become computationally expensive in 10- and 15-dimension cases, especially the exact methods cannot output the results within 2 hours in 15-dimension case. Notice that IWFG method can only find the solution with the smallest hypervolume contribution in a solution set, while the other methods can approximate and calculate the hypervolume contributions of all solutions in a solution set.}

\textcolor{black}{
The results in Table \ref{compareexact} suggest that the approximation methods are more preferable in high-dimensional spaces (e.g., $>$10-dimension) where the exact methods become computationally expensive and impractical to apply.}
\begin{table}[!htbp]
\small
\centering
  \caption{Runtime comparison of five methods on linear triangular PF solution sets (Seconds)}
  \label{compareexact}
  \begin{tabularx}{\linewidth}{@{\extracolsep{\fill}}cc r r r}
    \toprule
    Methods&\#Vectors/Points&5-D& 10-D & 15-D\\
    \midrule
    New  &100 & 2.9167 & 4.1289 & 5.0976\\
    New  &1000 & 27.9614 & 40.0872 &49.7242\\
     Traditional  &100 & 2.9638 & 3.9814 &4.8794\\
     Traditional  &1000 & 29.3457 & 39.1236 &47.0043\\
     Monte Carlo &100 & 1.6563 & 2.5491 & 3.5167\\
     Monte Carlo &1000 & 13.2484 & 20.2383 & 26.6723\\
     IWFG &- & 0.1098 & 71.8063 & $>$7200\\
     exQHV &- & 0.6358 & 665.9809 & $>$7200 \\
  \bottomrule
\end{tabularx}
\end{table}

\section{Conclusions}
\label{conclude}
In this letter, a new method for the hypervolume contribution approximation was proposed. \textcolor{black}{
From the experiment results we obtained the following insights and conclusions:
\begin{enumerate}
\item The Monte Carlo method only performs well when the reference point is the nadir point of the PF. The new method is able to achieve a good and robust performance with respect to the specification of the reference point. From a practical point of view, the new method is more suitable to apply in the hypervolume-based EMOAs.
\item The new method with a small number of the direction vectors is able to achieve a comparable or even better performance than the other two approximation methods with a large number of the direction vectors and the sampling points. In this sense, the new method can consume less time and at the same time achieve higher approximation accuracy. 
\item The new method showed a good and robust performance with respect to the number of solutions while the other two approximation methods can not.
\item The new method is computationally more efficient than the two exact methods in high-dimensional spaces (e.g., $>$10-dimension) where the exact methods become impractical to use. 
\end{enumerate}
}

For our future research, we will develop an indicator-based EMOA based on the new method and compare it with the hypervolume-based EMOAs \cite{beume2007sms,bader2011hype,jiang2015simple} and other R2 indicator-based EMOAs \cite{phan2013r2,hernandez2015improved,brockhoff20152}. Another future research direction is to further improve the performance of the new method \textcolor{black}{in terms of the approximation accuracy and the computational efficiency.}

The solution sets and the source code of our experiments are available at \url{https://github.com/nixizi/R2-HVC}.

\ifCLASSOPTIONcaptionsoff
  \newpage
\fi

\bibliographystyle{IEEEtran}
\bibliography{sample-bibliography} 

\end{document}